\documentclass[11pt]{article}

\usepackage{amssymb}
\usepackage{latexsym}
\usepackage{amsmath}
\usepackage{epsfig}
\usepackage{amscd}
\usepackage{amsthm}
\usepackage{multirow}
\usepackage{float}
\usepackage{footmisc}
\usepackage{color}
\usepackage{enumerate}
\usepackage{bm}
\usepackage{algorithm,algorithmic}

\usepackage{setspace}
\usepackage{adjustbox} 
\usepackage{xparse}
\usepackage{subfigure}
\usepackage{tcolorbox}
\usepackage{lipsum}
\usepackage{blindtext}
\usepackage{mathtools}
\usepackage[titletoc,title]{appendix}
\usepackage[T1]{fontenc}
\usepackage{dsfont}
\usepackage[sort]{natbib}
\usepackage{url}
\usepackage[font=small,labelfont=bf]{caption}
\usepackage{authblk} 
\usepackage{hyperref}[] 
\hypersetup{
	colorlinks = true,
	linkcolor = blue!60!black, 
	filecolor = magenta,
	urlcolor = blue!60!black, 
	citecolor = blue!60!black 
}

\usepackage{sectsty}

\newcommand{\mP}{\mathbb{P}}
\newcommand{\mE}{\mathbb{E}}

\newcommand{\tV}{\text{Var}}

\newcommand{\bW}{\boldsymbol{W}}

\newcommand\independent{\protect\mathpalette{\protect\independenT}{\perp}}
\def\independenT#1#2{\mathrel{\rlap{$#1#2$}\mkern2mu{#1#2}}}

\DeclarePairedDelimiter\ceil{\lceil}{\rceil}

\newtheorem{definition}{Definition}
\newtheorem{theorem}{Theorem}
\newtheorem{lemma}{Lemma}

\newtheorem{proposition}{Proposition}
\newtheorem{corollary}{Corollary}
\newtheorem{remark}{Remark}

\begin{document}

\begin{center}
	\textbf{\LARGE Local permutation tests for conditional independence}	
	
	\vspace*{.2in}

	\begin{author}
	A
	Ilmun Kim$^{\ddagger}$  \quad Matey Neykov$^{\dagger}$ \quad Sivaraman Balakrishnan$^{\dagger}$ \quad Larry Wasserman$^{\dagger}$
	\end{author}

	\vspace*{.2in}
	
	\begin{tabular}{c}
		$^{\dagger}$Department of Statistics and Data Science, Carnegie Mellon University\\
		$^{\ddagger}$Department of Statistics and Data Science, Yonsei University
	\end{tabular}

	\vspace*{.2in}
	
	\today
		
	\vspace*{.2in}
\end{center}

\begin{abstract}
	In this paper, we investigate local permutation tests for testing conditional independence between two random vectors $X$ and $Y$ given $Z$. 
	The local permutation test determines the significance of a test statistic by locally shuffling samples which share similar values of the conditioning variables $Z$,
	and it forms a natural extension of the usual permutation approach for
	unconditional independence testing. 
	Despite its simplicity and empirical support, the theoretical underpinnings of the local permutation test remain unclear. 
	Motivated by this gap, this paper aims to establish theoretical foundations of local permutation tests with a particular focus on binning-based statistics. We start by revisiting the hardness of conditional independence testing and provide an upper bound for the power of any valid conditional independence test, which 
	holds when the probability of observing ``collisions'' in $Z$ is small.  
	This negative result naturally motivates us to impose additional restrictions on the possible distributions under the null and alternate.
	 To this end, we focus our attention on certain classes of smooth distributions and identify provably tight conditions under which the local permutation method is universally valid, i.e.~it is valid when applied to any (binning-based) test statistic. To complement this result on type I error control, we also show that in some cases, a binning-based statistic calibrated via the local permutation method can achieve minimax optimal power.
	 We also introduce a double-binning permutation strategy, which yields a valid test over less smooth null distributions than the typical single-binning method without compromising much power.
	 Finally, we present simulation results to support our theoretical findings.
\end{abstract}

\vskip 2em

\section{Introduction}
Conditional independence (CI) is an important concept in a variety of statistical applications including graphical models~\citep{de2000new,koller2009probabilistic} and causal inference~\citep{spohn1994properties,pearl2014probabilistic,imbens2015causal}. In these applications, the assumption of conditional independence offers significant representational and computational benefits, and helps disentangle causal relationships among variables in an efficient and tractable way.
In a related vein, a problem of essential importance in statistical practice is that of variable selection~\citep{williamson2021nonparametric,dai2021significance}, which 
is concerned with selecting a parsimonious subset of features 
that are predictive of a response variable. In each of these settings, conditional independence tests are an essential tool to validate (or invalidate) critical modeling assumptions, and can lend additional credibility to the conclusions of our data analysis.

The performance of a statistical hypothesis 
test relies not only on the form of the test statistic but also heavily on the method used to ensure type I error control. Indeed, one might argue that a huge part of the practical success and ubiquity
of two-sample and (unconditional) independence tests is the fact that these tests can be tightly calibrated in a black-box fashion using a permutation method. This in turn frees the practitioner to focus on designing
powerful test statistics, without having to further ensure that the distribution of their test statistics are analytically tractable under the null. For two-sample and (unconditional) independence testing,
the permutation method is universal without any additional restrictions, i.e.~it controls the type I error in a non-trivial sense for \emph{any} underlying test statistic. As noted by \cite{shah2020hardness}, part of the hardness of conditional independence testing with a continuous 
conditioning variable $Z$ stems from the fact that it is impossible to control the type I error, via for instance a permutation method, in any non-trivial sense without additional restrictions. 
Our broad goal in this paper is to propose and study natural extensions of the permutation method, namely the local permutation procedure, which are applicable to CI testing. In particular, we aim to investigate restrictions under which these extensions tightly control the type I error for a broad class of test statistics, and further to explore the power of tests calibrated via these methods.

The \emph{local} permutation procedure calibrates a test statistic by locally shuffling samples based on the proximity of their conditioning variables $Z$. When the conditional variable is discrete, the resulting local permutation test has a universal guarantee on type I error control under relatively weak assumptions on the data-generating process. When the conditional variable is continuous, on the other hand, the validity of the local permutation test is far from obvious. While there is a line of work providing some empirical support~\citep{fukumizu2008kernel,doran2014permutation,sen2018mimic,neykov2020minimax}, a rigorous theoretical foundation of the local permutation test has not been fully established in the continuous case. Motivated by this gap, the first aim of this paper is to identify provably tight conditions under which the type I error of the local permutation test is uniformly controlled at least for sufficiently large sample-sizes. To this end, we focus primarily on a binning-based local permutation procedure and determine the size of bins for which the resulting test is asymptotically valid under various smoothness assumptions. 

Once the type I error is under control, our subsequent focus is on power. In contrast to type I error control, which requires the size of bins to be small, the use of bins that are too fine causes a loss of power due to the small sample size in each bin. Our next goal is to balance this trade-off and show that there is a choice of bin-widths which ensures that the local permutation method controls the type I error but still retains minimax optimal power.
We achieve this goal by building on the recent work of \cite{canonne2018testing} and \cite{neykov2020minimax} which study minimax-optimal CI tests and the work of \cite{kim2020minimax} which studies the power of the classical permutation method for two-sample and independence testing.

An interesting aspect of our results is that they elucidate a tension in conditional independence testing 
between ensuring tight control of the type I error, and ensuring high power of the resulting test. 
In many well-studied examples, 
permutation and other simulation 
methods represent an apparent free lunch,
ensuring tight control of the type I error without sacrificing power (see for instance, \cite{kim2020minimax} for precise results on the minimax power of the permutation test in these settings). 
In conditional independence testing, the permutation method is no longer exact and we show that there is a trade-off when using the local permutation method for calibration. In certain cases, ensuring type I error control requires selecting bin-widths which are too small to guarantee high power. In some settings, we are able to mitigate this trade-off by designing 
a careful double-binning strategy
where two resolutions are combined in the permutation method: a finer resolution for permutations which ensures type I error control, and a coarser resolution for computing the test statistic which 
ensures high power (see Section~\ref{Section: Double-binning strategy}). Before we state our contributions in more detail, we briefly review related work.

\subsection{Related work} \label{Section: Related work}
There is an extensive body of literature on CI measures and CI tests. Here, we give a selective review of existing methods, which can be categorized into several groups. 

The first category of methods is based on kernel mean embeddings~\citep[see][for a review]{muandet2016kernel}. The idea of kernel mean embeddings is to represent probability distributions as elements of a reproducing kernel Hilbert space~(RKHS), which enables us to understand properties of these distributions using Hilbert space operations. One of the initial attempts to use kernel mean embeddings for CI testing was made by \cite{fukumizu2008kernel}. In particular, \cite{fukumizu2008kernel} propose a test based on the empirical Hilbert--Schmidt norm of the conditional cross-covariance operator. \cite{zhang2012kernel} introduce another kernel-based test attempting to measure partial correlations, which in turn characterize CI \citep{daudin1980partial}. \cite{strobl2019approximate} use random Fourier features to approximate kernel computations, and propose a more computationally efficient version of the test of \cite{zhang2012kernel}. Other CI measures proposed by \cite{doran2014permutation} and \cite{huang2020kernel} are motivated by the kernel maximum mean discrepancy for two-sample testing~\citep{gretton2012kernel}. In particular, the CI measure introduced by \cite{huang2020kernel} compares whether $Y|X,Z$ and $Y|X$ have the same distribution, and their measure can be viewed as a kernelized version of the CI measure of \cite{azadkia2019simple}. Recently, \cite{sheng2019distance} and \cite{park2020measure} propose kernel CI measures that are closely connected to the Hillbert--Schmidt independence criterion~\citep{gretton2005measuring}. \cite{sheng2019distance} also discuss the connection between their CI measure to the conditional distance correlation proposed by \cite{wang2015conditional}.

Another category of methods relies on estimating regression functions. Consider random variables $X$ and $Y$, and their regression residuals on $Z$, denoted by $\epsilon_{X,Z} :=X - \mE[X|Z]$ and $\epsilon_{Y,Z}:= Y-\mE[Y|Z]$. The underlying idea of regression-based methods is that the expected value of $\epsilon_{X,Z}\epsilon_{Y,Z}$ is zero if $X \independent Y |Z$, and not necessarily zero if $X \not\independent Y|Z$. Thus, one can use an empirical estimate of the expected value of  $\epsilon_{X,Z}\epsilon_{Y,Z}$ as a test statistic for CI. We refer to \cite{zhang2018measuring} for a discussion of the relationship between $\epsilon_{X,Z} \independent \epsilon_{Y,Z}$ and $X \independent Y |Z$. Given that there exist a variety of successful regression algorithms to estimate $\mE[X|Z]$ and $\mE[Y|Z]$, the expected value of $\epsilon_{X,Z}\epsilon_{Y,Z}$ can be accurately estimated as well. This regression-based idea has been exploited by several authors to tackle CI testing. For instance, \cite{shah2020hardness} propose the generalized covariance measure, which has been extended to functional linear models by \cite{lundborg2021conditional}. 
The methods proposed by \cite{zhang2012kernel} and \cite{strobl2019approximate} rest on the regression of a function in a RKHS, thereby belonging to this category as well. 
We also note that there has been a growing interest in estimating the expected conditional covariance in semi-parametric statistics \citep{robins2008higher,li2011higher,newey2018cross} often employing non-parametric 
regression methods followed by adjustments to reduce bias, and this work in turn has implications for the design of CI tests. 
  
Apart from the above two categories, there are many other novel approaches for CI testing developed in recent years. For example, \cite{bellot2019conditional,shi2020double} design nonparametric tests by leveraging the success of generative adversarial networks. \cite{sen2017model,sen2018mimic} convert the CI testing problem into a binary classification problem, which allows one to leverage existing classification algorithms. Approaches based on the partial copula have been examined by \cite{bergsma2004testing,bergsma2010nonparametric,song2009testing,patra2016nonparametric,petersen2021testing}. A metric-based approach is also common in the literature, including tests based on the conditional Hellinger distance~\citep{su2008nonparametric} and conditional mutual information~\citep{runge2018conditional}. The above methods are mainly for continuous data, whereas there are numerous CI tests available for discrete data as well~\citep{agresti1992survey,yao1993exact,kim1997nearly,canonne2018testing,marx2019testing,neykov2020minimax,balakrishnan2018hypothesis}. A more extensive review of CI tests can be found in \cite{li2020nonparametric}.

So far we have mainly reviewed various ways of measuring CI and constructing test statistics. For testing problems, it is also important to determine a reasonable critical value, that results in small type I and type II errors. The current literature usually considers one of the following three approaches for setting critical values. 

\begin{itemize}
	\item \textbf{Asymptotic method.} The first common approach is based on the limiting null distribution of a test statistic. Once the limiting null distribution is known, the critical value is determined by using a quantile of this limiting distribution or a bootstrap procedure. In order to obtain a tractable asymptotic distribution, the test statistic typically has an asymptotically linear or quadratic form. Examples of CI tests based on the asymptotic approach include~\cite{su2008nonparametric,huang2010testing,zhang2012kernel,wang2015conditional,strobl2019approximate,shah2020hardness,zhou2020test}. 
	Due to technical hurdles, this line of work often focuses on a pointwise (rather than uniform) type I error guarantee with a few exceptions~\citep{shah2020hardness,lundborg2021conditional}. 
	\item \textbf{Model-X framework.} Formalized by \cite{candes2016panning}, the model-X framework builds on the assumption that the conditional distribution $P_{X|Z}$ is (approximately) known. In this case, one can compute a set of test statistics, which are exchangeable under the null, by exploiting the knowledge of $P_{X|Z}$ either through direct resampling as in \cite{candes2016panning} or via a permutation method as in \cite{berrett2020conditional}. The critical value is then set to be an empirical quantile of these test statistics, and the resulting test has finite-sample validity. \cite{berrett2020conditional} rigorously characterize the excess type I error when an estimate of $P_{X|Z}$ is considered, and also demonstrate situations where this excess error is asymptotically negligible. Nevertheless, this methodology may not be appropriate for applications where $P_{X|Z}$ is hard to estimate.   
	\item \textbf{Local permutation method.} The third approach is based on local permutations. This method generates a reference distribution by randomly permuting $Y$ within subclasses, which are defined in terms of the proximity of the conditional variable $Z$. Then the critical value is determined as a quantile of this reference distribution. The work of \cite{margaritis2005distribution,fukumizu2008kernel,doran2014permutation,sen2017model} fall into this category. When observing multiple samples with the same value of $Z$ is possible, this method can yield an exact CI test with reasonable power against certain alternatives. However, its validity has not been fully explored beyond discrete settings.
\end{itemize}

As mentioned before, our work heavily builds upon the recent work of \cite{canonne2018testing}, \cite{neykov2020minimax} and \cite{kim2020minimax}. \cite{canonne2018testing} construct tests for CI when $(X,Y,Z)$ are discrete random variables, but with a possibly large number of categories, and establish the optimality of their tests from a minimax perspective, in certain regimes. \cite{neykov2020minimax} extend the work of \cite{canonne2018testing} to the case where $Z$ is a continuous and bounded random variable. However, both tests considered in \cite{canonne2018testing} and \cite{neykov2020minimax} rely on critical values that depend on unspecified constants. In this sense, it has been unknown whether there exists a minimax optimal CI test, which is easily implementable in practice. To address this issue, our work considers the local permutation method, which leads to an explicit critical value. In order to analyze the power of the resulting test, we 
build on results of \cite{kim2020minimax} who provide a sufficient condition under which the permutation test has non-trivial power for (unconditional) independence testing. To verify this sufficient condition, we build on the analysis of U-statistic-based tests from the previous work of \cite{canonne2018testing} and \cite{neykov2020minimax}.

\subsection{Our contributions} \label{Section: Our contributions}
\begin{figure}[t!]
	\begin{center}		
		\includegraphics[width=30em]{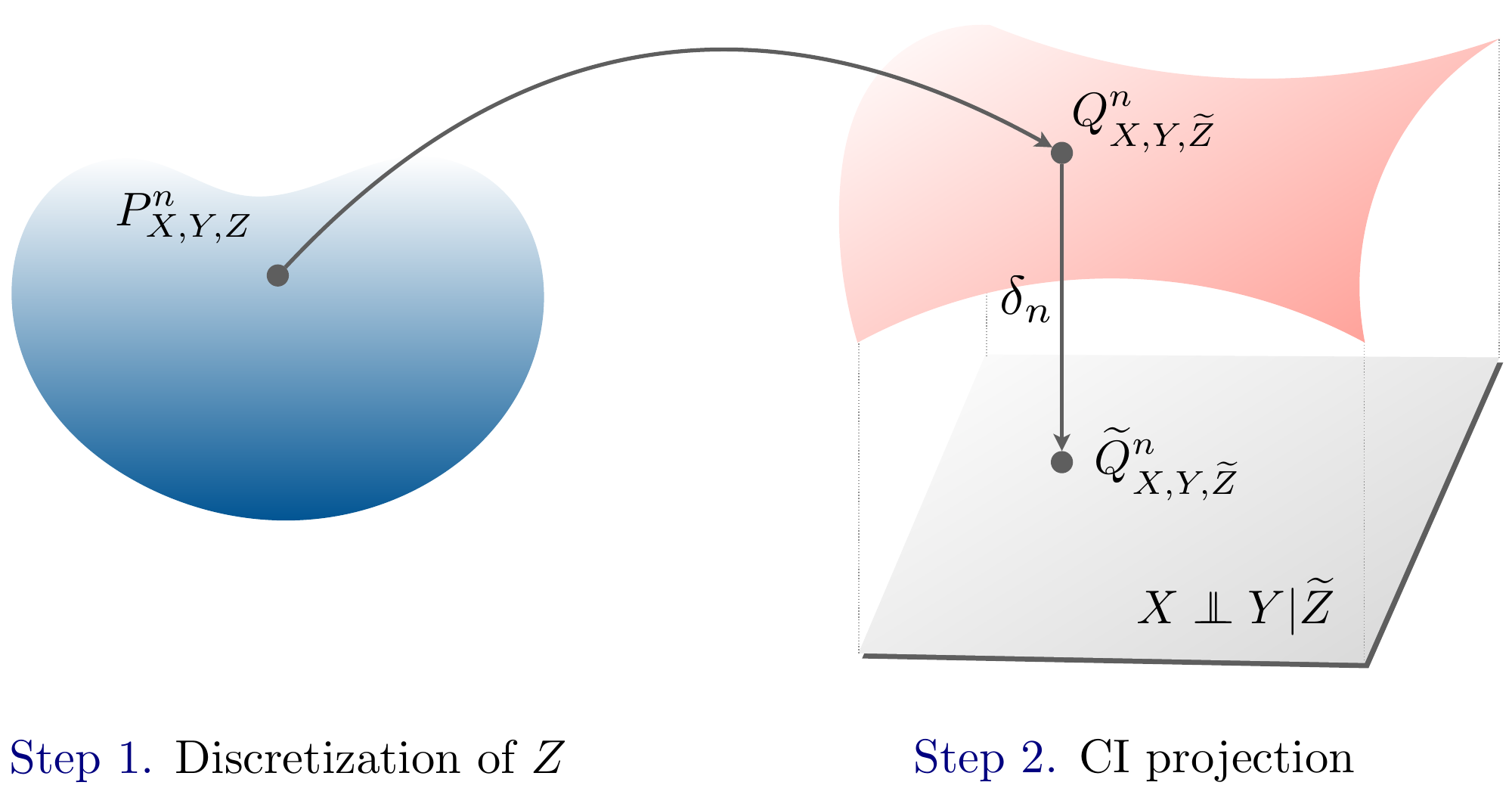} 
		\caption{\small A visualization of our analysis of the local permutation test. To proceed, we first discretize $Z$ into several bins and denote the binned conditional variable by $\widetilde{Z}$. Given a product distribution $P_{X,Y,Z}^n$, the corresponding distribution smoothed over the bins is denoted by $Q_{X,Y,\widetilde{Z}}^n$ (see Section~\ref{Section: Validity} for a more precise description). Next, we consider another distribution $\widetilde{Q}_{X,Y,\widetilde{Z}}^n$, which is the CI projection of $Q_{X,Y,\widetilde{Z}}^n$ onto the space where $X \independent Y |\widetilde{Z}$. The validity of the location permutation test is essentially determined by the TV distance $\delta_n$ (or its upper bound) between $Q_{X,Y,\widetilde{Z}}^n$ and $\widetilde{Q}_{X,Y,\widetilde{Z}}^n$. In particular, the local permutation test is asymptotically valid when $\delta_n \rightarrow 0$ as $n \rightarrow \infty$. In Section~\ref{Section: Validity under smoothness conditions}, we present sharp and tractable conditions under which $\delta_n$ converges to zero, depending on the smoothness of underlying distributions.} \label{Figure: CI projection}
	\end{center}
\end{figure}

We now outline our contributions. 
\begin{itemize}
	\item \textbf{Hardness result of CI testing~(Section~\ref{Section: Fundamental limits of CI testing}).} By leveraging the recent work of \cite{foygel2019limits,barber2020distribution}, our first contribution~(Theorem~\ref{Theorem: Hardness of CI testing}) is to provide a new hardness result for CI testing. For two-sample and unconditional independence testing, one can use the permutation procedure to develop tests that can keep the type I error under control, while having non-trivial power against interesting alternatives~\citep[e.g.][]{kim2020minimax}. However, this is not the case for the continuous CI testing problem. As pointed out by \cite{shah2020hardness}, any valid CI test should have no power against any alternative when $Z$ is a continuous random variable. In Theorem~\ref{Theorem: Hardness of CI testing}, we formalize that the impossibility of CI testing is more fundamentally determined by the probability of observing collisions in $Z$, rather than the type of $Z$. Therefore, even in the discrete or mixture setting, CI testing is difficult or even impossible when the probability of observing the same $Z$ is extremely small. 	
	\item \textbf{Validity of local permutation tests~(Section~\ref{Section: Validity}).} In continuous settings, one typical way to address the replication problem of $Z$ in the design of tests is to hypothesize some notion of smoothness, i.e.~the conditional distribution does not vary too much as a function of $Z$. Under this hypothesis, we have approximate replicates, which we can use to construct our test statistics, and can use to design an (approximate) permutation method. A basic question is to address the validity of the local permutation method. Our preliminary results~(Lemma~\ref{Lemma: Error bound in terms of the TV distance} and Lemma~\ref{Lemma: Error bound in terms of the Hellinger distance}) show that one can control the type I error of binning based CI statistics when the binned distribution is indistinguishable from its CI projection (the distribution we obtain from permutations). See Figure~\ref{Figure: CI projection} for a pictorial illustration. 
	\item \textbf{Tightness of our conditions~(Section~\ref{Section: Tightness}).} We note that, counterintuitively, increasing the sample size $n$ can make type I error control harder to achieve, because ensuring the indistinguishability of the product measures is more challenging as $n$ increases. This forces us to use finer bins as $n$ increases for type I error control. On the other hand, using bins that are too fine may result in a loss of power, which raises the question on a choice of the size of bins. In Theorem~\ref{Theorem: Validity of the permutation test under Hellinger smooth} and Theorem~\ref{Theorem: Validity of the permutation test under Renyi smooth}, we first present concrete upper bounds for the type I error in terms of the size of bins and the sample size $n$, under certain smoothness conditions. These results guide us on the size of bins, which yields rigorous type I error control. As a complementary result, Theorem~\ref{Theorem: Lower bounds} proves that the upper bounds in the previous results are asymptotically tight. In particular, we show that there exists a local permutation test whose type I error is arbitrarily close to one when the given upper bounds are sufficiently far from the significance level. 	
	\item \textbf{Power analysis~(Section~\ref{Section: Power analysis and optimality}).} The next question we address is that of power. We start by revisiting the test statistics for discrete CI testing in \cite{canonne2018testing}. Theorem~\ref{Theorem: Type II error under discrete setting} then shows that the corresponding local permutation tests have the same power guarantee as the tests in \cite{canonne2018testing}. Unlike the discrete CI setting, the replication problem affects both type I and type II error control in the continuous CI setting. As mentioned earlier, taking finer bins helps us for type I error control but not for power. Our next result shows that in some cases, we are able to navigate this trade-off, i.e.~there is a choice of binning which ensures that the local permutation method controls the type I error but still retains minimax power. In particular, we show in Theorem~\ref{Theorem: Type II error under continuous setting I} and Theorem~\ref{Theorem: Type II error under continuous setting II} that the local permutation tests using the same test statistics in \cite{neykov2020minimax} achieve the same minimax power in the total variation (TV) distance. However, this guarantee comes at a cost. Namely, the local permutation test with the optimal choice of binning is valid over a set of null distributions much smoother than those considered in \cite{neykov2020minimax}. 
	\item \textbf{Double-binning strategy~(Section~\ref{Section: Double-binning strategy}).} Finally, we develop and analyze a new double-binning based permutation test, which partly addresses the aforementioned drawback. More specifically, we consider bins of two distinct resolutions where finer bins are used for permutations and coarser bins are used to compute a test statistic. By permuting over the finer bins, our theory in Proposition~\ref{Proposition: Type I error of the double-binned test} guarantees that a double-binning based test has type I error control over a larger class of null distributions than the single-binning counterpart. On the other hand, by computing the test statistic over the coarser bins, Theorem~\ref{Theorem: Type II error of the double-binned test} proves that the power of the resulting test remains the same as the single-binning test, up to a constant factor, under certain regularity conditions. We further demonstrate our theoretical findings in Section~\ref{Section: Experiment 3} through simulations. 
\end{itemize}

\subsection{Organization}
The rest of this paper is organized as follows. We start by explaining the local permutation procedure in Section~\ref{Section: Preliminaries} along with a basic background on probability metrics. Section~\ref{Section: Fundamental limits of CI testing} provides a new hardness result of CI testing that covers the discrete case of $Z$. We then move on to discussing the validity of the local permutation test in Section~\ref{Section: Validity under smoothness conditions}. In particular, we provide upper bounds for the type I error of the local permutation test under certain smoothness conditions. We further prove that these upper bounds are asymptotically tight in some cases. Focusing on the test statistics proposed by \cite{canonne2018testing} and \cite{neykov2020minimax}, we investigate the power property of the local permutation tests in Section~\ref{Section: Power analysis and optimality}. Section~\ref{Section: Double-binning strategy} introduces a double-binning strategy that allows us to choose a smaller binning size without sacrificing power up to a constant factor. Section~\ref{Section: Simulations} includes several illustrative simulation results. Finally, we end the paper with a discussion and future work in Section~\ref{Section: Discussion}. All technical proofs are relegated to the Appendix. 

\section{Preliminaries} \label{Section: Preliminaries}

In this section, we set up the notation and introduce preliminaries including the local permutation procedure and probability metrics.

\subsection{Notation}
Throughout this paper, we mostly follow the notation used in \cite{neykov2020minimax}. Let the triplet $(X,Y,Z)$ have a distribution $P_{X,Y,Z}$ on a measurable space. We denote the conditional distribution of $X,Y|Z=z$ as $P_{X,Y|Z=z}$. We denote the (marginal) conditional distributions of $X|Z=z$ and $Y|Z=z$ by $P_{X|Z=z}$ and $P_{Y|Z=z}$, respectively. In addition, the marginal distributions of $X,Y,Z$ are denoted by $P_X,P_Y,P_Z$ and similarly the joint marginal distributions are denoted by $P_{X,Y},P_{X,Z},P_{Y,Z}$. Moreover, we will use the lowercase $p$ to denote density functions with respect to a base measure. For example, $p_{X,Y|Z}(x,y|z)$, denotes the conditional density (or probability mass) function of $X,Y|Z=z$, 
evaluated at a point $(x,y,z)$. We denote the set of all distributions for which $X \independent Y|Z$ by $\mathcal{P}_0$.

\subsection{Local permutation procedure} \label{Section: Local permutation procedure}
We formalize the local permutation procedure based on i.i.d.~observations $\{(X_i,Y_i,Z_i)\}_{i=1}^n:=(X^n,Y^n,Z^n)$ from $P_{X,Y,Z}$. Throughout this paper, we assume the conditional random variable $Z$ has compact support $\mathcal{Z}$ and briefly discuss an extension to unbounded support in Appendix~\ref{Section: Unbounded support}. Let $\{B_1,\ldots,B_M\}$ denote a partition of $\mathcal{Z}$ such that $\mathcal{Z} = \cup_{m=1}^M B_m$ and $\boldsymbol{\sigma}:=\{\sigma_1,\ldots,\sigma_M\}$ denote the sample sizes within bins $\{B_1,\ldots,B_M\}$. Furthermore, let $\bW_m$ denote the set of the pairs of $(X_i,Y_i)$ that belong to the $m$th bin. More formally, by letting $(X_{i,m},Y_{i,m})$ be the $i$th pair in the $m$th bin, we write $\bW_m := \{(X_{1,m},Y_{1,m}),\ldots, (X_{\sigma_m,m},Y_{\sigma_m,m})\}$ when $\sigma_m \geq 1$ and otherwise $\bW_m= \emptyset$.

Given this binned data, we consider a generic test statistic for CI testing, which maps from $\bW_1,\ldots,\bW_M$ to $\mathbb{R}$, i.e. for some function $f: \bW_1,\ldots,\bW_M \mapsto \mathbb{R}$ 
we compute our test statistic as: 
\begin{align}  \label{Eq: test statistic for CI testing}
T_{\mathrm{CI}} = f(\bW_1,\ldots, \bW_M).
\end{align}
As a concrete example with real-valued data, one can take $f$ to be the average function of arbitrary (unconditional) independence test statistics computed based on $\bW_1,\ldots, \bW_M$, respectively.

In order to determine significance of the statistic $T_{\mathrm{CI}}$, we rely on the local permutation procedure summarized in Algorithm~\ref{alg:local permutation procedure}. To describe the algorithm, consider a permutation $\pi_m=\{\pi_m(1),\ldots,\pi_m(\sigma_m)\}$ of $\{1,\ldots,\sigma_m\}$ and denote $\bW_{m}^{\pi_m} =  \{(X_{1,m},Y_{\pi_m(1),m}),\ldots, (X_{\sigma_m,m},Y_{\pi_m(\sigma_m),m})\}$ for $m=1,\ldots,M$. Notice that when $\sigma_m = 0$, there is nothing to permute and we set $\bW_{m}^{\pi_m}  = \emptyset$. The test statistic computed using the locally permuted data set is denoted by
\begin{align} \label{Eq: permuted statistic}
	T_{\mathrm{CI}}^{\boldsymbol{\pi}} = f( {\bW}_1^{\pi_1},\ldots, {\bW}_M^{\pi_M}).
\end{align}
Let us further denote the set of all possible such local permutations $\boldsymbol{\pi}:=\{\pi_1,\ldots,\pi_M\}$ by $\boldsymbol{\Pi}$ whose cardinality is $K=\prod_{m=1}^M \sigma_m!$. Given this notation, we describe the local permutation procedure in Algorithm~\ref{alg:local permutation procedure}. 

\begin{algorithm}[h]
	\caption{Local permutation procedure} \label{alg:local permutation procedure}
	\begin{algorithmic}
		\ENSURE  data $\{(X_i,Y_i,Z_i)\}_{i=1}^n$, a partition of $\mathcal{Z}$: $\{B_1,\ldots,B_M\}$, a test statistic $T_{\mathrm{CI}}$, a nominal level $\alpha$
	\begin{enumerate}
		\item For each $\boldsymbol{\pi} \in \boldsymbol{\Pi}$, compute $T_{\mathrm{CI}}^{\boldsymbol{\pi}}$ as in (\ref{Eq: permuted statistic}) and denote the resulting statistics by $T_{\mathrm{CI}}^{\boldsymbol{\pi}_1},\ldots,T_{\mathrm{CI}}^{\boldsymbol{\pi}_K}$. 
		\item By comparing the statistic $T_{\mathrm{CI}}$ in~(\ref{Eq: test statistic for CI testing}) with the permuted ones, calculate the $p$-value as
		\begin{align} \label{Eq: permutation p-value}
			p_{\mathrm{perm}} = \frac{1}{K} \sum_{\boldsymbol{\pi}_i \in \boldsymbol{\Pi}} \mathds{1} \big\{ T_{\mathrm{CI}}^{\boldsymbol{\pi}_i} \geq T_{\mathrm{CI}} \big\}.
		\end{align}
		\item Given the nominal level $\alpha \in (0,1)$, define the test function $\phi_{\mathrm{perm},n} = \mathds{1}(p_{\mathrm{perm}} \leq \alpha)$ and reject the null when $\phi_{\mathrm{perm},n}=1$.
	\end{enumerate}
	\end{algorithmic}
\end{algorithm}

The local permutation procedure, like other randomized or permutation procedures \citep[e.g.~Chapter 15 of][]{lehmann2006testing}, can be used with any binning-based test statistic for CI testing. For simplicity, our theoretical results are based on Algorithm~\ref{alg:local permutation procedure} but they can be easily extended to a more practical permutation procedure via Monte Carlo simulations as remarked below.

\begin{remark} \normalfont \leavevmode
	\begin{itemize} \label{Remark: MC approximation/randomization}
		\item \textbf{Monte Carlo approximation.} The permutation $p$-value~(\ref{Eq: permutation p-value}) may be practically unappealing as its computational cost is prohibitively expensive for large $n$. To alleviate this computational issue, it is a common practice to approximate $p_{\mathrm{perm}}$ using Monte Carlo simulations as in (\ref{Eq: approximated p-value}). As noted in \cite{lehmann2006testing}, the difference between $p_{\mathrm{perm}}$ and its Monte Carlo approximation can be made arbitrarily small by taking a sufficiently large number of Monte Carlo samples. This can be formally stated using Dvoretzky--Kiefer--Wolfowitz inequality~\citep{dvoretzky1956asymptotic} and we refer to Corollary 6.1 of \cite{kim2019comparing} or Proposition 4 of \cite{schrab2021mmd} for such argument. 
		\item \textbf{Randomization.} It is well-known that the permutation test can be made exact by introducing randomization. We state the randomized permutation test \citep{hoeffding1952large} for completeness. For a nominal level $\alpha$, we denote $k = K - [K\alpha]$ where $[K\alpha]$ is the largest integer less than or equal to $K\alpha$. In addition let $K^+$ and $K^0$ be the numbers of $T_{\mathrm{CI}}^{\boldsymbol{\pi}_1},\ldots,T_{\mathrm{CI}}^{\boldsymbol{\pi}_K}$, which are greater than or equal to $T_{\mathrm{CI}}$, respectively. Given $a = (K\alpha - K^+) /K^0$ and the $k$th order statistic $T_{\mathrm{CI}}^{(k)}$ of $T_{\mathrm{CI}}^{\boldsymbol{\pi}_1},\ldots,T_{\mathrm{CI}}^{\boldsymbol{\pi}_K}$, we define $\phi_{\mathrm{perm},n,a} = 1,a$ or $0$ depending on whether $T_{\mathrm{CI}}  > T_{\mathrm{CI}}^{(k)}$, $T_{\mathrm{CI}}  = T_{\mathrm{CI}}^{(k)}$ or $T_{\mathrm{CI}}  < T_{\mathrm{CI}}^{(k)}$, respectively. Then under the exchangeability assumption of the permuted statistics, it holds that $\mE[\phi_{\mathrm{perm},n,a}] = \alpha$, whereas $\phi_{\mathrm{perm},n}$ from Algorithm~\ref{alg:local permutation procedure} has a weaker guarantee that $\mE[\phi_{\mathrm{perm},n}] \leq \alpha$ in general.
	\end{itemize}
\end{remark}

In the next subsection, we present several statistical distances between probability measures that we make use  of throughout this paper.

\subsection{Probability metrics}
Let $P$ and $Q$ be two probability measures over a measurable space $(\Omega, \mathcal{F})$ and denote the densities of $P$ and $Q$ with respect to a common dominating measure $\mu$ by $p$ and $q$, respectively. There are two classes of probability metrics that will be considered in this paper. The first class, we call the generalized Hellinger distance~\citep[e.g.][]{kamps1989hellinger}, includes the TV distance and the Hellinger distance as special cases. 
\begin{definition}[Generalized Hellinger distances] \normalfont 
	Given $\gamma \geq 1$, the generalized Hellinger distance with parameter $\gamma$ between $P$ and $Q$ is defined as
	\begin{align*}
		\mathcal{D}_{\gamma\text{,H}}(P,Q) = \bigg( \frac{1}{2} \int  \big| p^{1/\gamma} - q^{1/\gamma} \big|^\gamma d\mu \bigg)^{1/\gamma}.
	\end{align*}
\end{definition}
From the definition, it is clear that the above distance becomes the TV distance when $\gamma = 1$ and the Hellinger distance when $\gamma = 2$. Since these two values deserve special attention, we denote the corresponding TV distance and Hellinger distance by $\mathcal{D}_{\mathrm{TV}}(P,Q)$ and $\mathcal{D}_{\mathrm{H}}(P,Q)$, respectively. The generalized Hellinger distance has the monotonicity property that $\mathcal{D}_{\gamma_2\text{,H}}^{\gamma_2}(P,Q) \leq \mathcal{D}_{\gamma_1 \text{,H}}^{\gamma_1}(P,Q)$ for $1 \leq \gamma_1 \leq \gamma_2$ \citep[Corollary 3 of][]{kamps1989hellinger}. This monotonic relationship generalizes the well-known inequality between the TV and Hellinger distances, namely $\mathcal{D}_{\mathrm{H}}^2(P,Q) \leq \mathcal{D}_{\mathrm{TV}}(P,Q)$ \citep[e.g.~Chapter 4 of][]{le2012asymptotic}.

Another class of probability metrics that we consider is R{\'e}nyi divergence defined as follows. 
\begin{definition}[R{\'e}nyi divergences] \normalfont
	For $\gamma \in (0,\infty)$, R{\'e}nyi divergence of order $\gamma$ of $P$ from $Q$ is defined as
	\begin{align*}
		\mathcal{D}_{\gamma \text{,R}}(P \| Q) = 
		\begin{cases} \displaystyle
			\frac{1}{\gamma - 1} \log \bigg\{  \int \left( \frac{p}{q}\right)^\gamma q d\mu \bigg\}, \quad & \text{if $\gamma \neq 1$,} \\[1.2em] \displaystyle
			\int p \log \left(\frac{p}{q}\right) d\mu, \quad & \text{if $\gamma = 1$.}
		\end{cases}
	\end{align*}
\end{definition}
In the above definition, some notable values of $\gamma$ include $\gamma \in \{1/2, 1,2\}$ and the corresponding R{\'e}nyi divergence is directly or indirectly associated with the Hellinger distance $(\gamma=1/2)$, Kullback--Leibler (KL) divergence $(\gamma=1)$ and $\chi^2$ divergence $(\gamma=2)$ as stated in Appendix~\ref{Section: Auxiliary lemmas}. We refer the reader to \cite{van2014renyi} and \cite{sason2016f} for more information on R{\'e}nyi divergences.

\section{Fundamental limits of CI testing} \label{Section: Fundamental limits of CI testing}
Before we start analyzing local permutation tests, we provide a new hardness result for CI testing. In view of the recent hardness result of \cite{shah2020hardness}, further revisited by \cite{neykov2020minimax}, CI testing is intrinsically difficult in the following sense. Let $\mathcal{P}_{\mathbb{R}^d}$ be the set of all distributions for $(X,Y,Z)$ on $\mathbb{R}^{d_x+d_y+d_z}$, and let $\mathcal{P}_{K} \subset \mathcal{P}_{\mathbb{R}^d}$ be the subset of $\mathcal{P}_{\mathbb{R}^d}$ whose support is defined within a $L_\infty$ ball of radius $K$. We also assume that the distributions in $\mathcal{P}_K$ are absolutely continuous with respect to the Lebesgue measure. Let $\mathcal{P}_{0,K} \subset \mathcal{P}_K$ be the subset of distributions such that $X \independent Y|Z$ and denote its complement by $\mathcal{P}_{1,K} = \mathcal{P}_K \setminus \mathcal{P}_{0,K}$. By denoting the joint distribution of $n$ i.i.d.~random vectors from $P_{X,Y,Z}$ by $P_{X,Y,Z}^n$, the result of \cite{shah2020hardness} states that any valid CI test $\phi$ for the class of null distributions $\mathcal{P}_{0,K}$ should satisfy $\mE_{P^n_{X,Y,Z}}[\phi] \leq \alpha$ for any $P_{X,Y,Z} \in \mathcal{P}_{1,K}$. In words, no CI test, which has valid type I error control for all absolutely continuous conditionally independent distributions, can have meaningful power against any single alternative distribution in $\mathcal{P}_{1,K}$. It therefore emphasizes that one should consider smaller sets of null and alternative distributions in order to make the CI testing problem feasible. 

The story, on the other hand, is different when $Z$ has a discrete or a mixture distribution where one can observe the same value of $Z$. In this case, by permuting the samples within groups having the same value of $Z$, the local permutation test can be valid, while possessing non-trivial power against certain alternatives. However, even in this case, there exists an intrinsic difficulty of CI testing when the probability of observing the same $Z$ is extremely small. We precisely characterize this challenge in the following theorem. To describe the result, let $\mathcal{P}_{0} \subset \mathcal{P}_{\mathbb{R}^d}$ be the subset of distributions such that $X \independent Y|Z$, and define $\mathcal{P}_1 = \mathcal{P}_{\mathbb{R}^d} \setminus \mathcal{P}_{0}$. Further, let $\mathcal{P}_{0,\mathrm{disc}} \subset \mathcal{P}_0$ be the subset of null distributions where $Z$ is supported on a countable set. Then our result is stated as follows. 

\begin{theorem}[Hardness of CI testing] \label{Theorem: Hardness of CI testing}
	For an arbitrary integer $J \geq n(n-1)$, let us define $\rho_{J,P} := \mP\{Z_1,\ldots,Z_J~\mathrm{\text{are distinct}}\}$, where $Z_1,\ldots,Z_J$ are i.i.d.~samples from the marginal distribution of $Z$. Suppose that a test $\phi$ satisfies $\sup_{P_{X,Y,Z} \in \mathcal{P}_{0,\mathrm{disc}}}\mE_{P_{X,Y,Z}^n}[\phi] \leq \alpha$ for $\alpha \in (0,1)$. Then for any $P_{X,Y,Z} \in \mathcal{P}_1$, the power of $\phi$ is bounded above by
	\begin{align} \label{Eq: hardness result}
		\mE_{P_{X,Y,Z}^n}[\phi] \leq \alpha \times \rho_{J,P} + (1- \rho_{J,P}) + \frac{n(n-1)}{J}.
	\end{align}
\end{theorem}

A few remarks on this result are given below.

\begin{remark} \leavevmode \normalfont
	\begin{itemize}
		\item Theorem~\ref{Theorem: Hardness of CI testing} states that what makes the CI problem hard is not just whether $Z$ is discrete or continuous, but whether one can observe the same value of $Z$ repeatedly with high probability. This difficulty is precisely captured by the quantity $\rho_{J,P}$. As an illustration, suppose that $Z$ has a multinomial distribution with equal probabilities over bins. Intuitively, when the number of bins is much larger than the sample size, one cannot expect to see the same $Z$ even twice with high probability, and thus the bound~(\ref{Eq: hardness result}) becomes close to $\alpha$. We make this intuition more precise in Remark~\ref{Remark: Results on impossibility} where we demonstrate that any valid CI test becomes (asymptotically) powerless if $Z$ is uniformly distributed over bins and the number of bins increases much faster than $n^4$. We also refer to Section~\ref{Section: Experiment 1} for a numerical illustration. 
		\item Our result is not restricted to the case of discrete $Z$. Suppose that $\phi$ is valid over a subset of $\mathcal{P}_0$ that contains $\mathcal{P}_{0,\mathrm{disc}}$. 
		Then the same result trivially follows. Let $\mathcal{P}_{1,\text{no-atom}} \subset \mathcal{P}_1$ be the subset of alternative distributions where the marginal distribution of $Z$ has no atoms, i.e.~$\rho_{J,P}=1$. Then, as a corollary of Theorem~\ref{Theorem: Hardness of CI testing}, it holds that $\sup_{P_{X,Y,Z} \in \mathcal{P}_{1,\text{no-atom}}} \mE_{P^n_{X,Y,Z}}[\phi] \leq \alpha$ for any $\phi$ such that $\sup_{P_{X,Y,Z} \in \mathcal{P}_{0}}\mE_{P_{X,Y,Z}^n}[\phi] \leq \alpha$. However, in this argument, it is crucial to assume that $\mathcal{P}_{0,\mathrm{disc}} \subset \mathcal{P}_0$. In contrast, the hardness result of \cite{shah2020hardness} does not require $\mathcal{P}_{0,\mathrm{disc}} \subset \mathcal{P}_0$. In this sense, Theorem~\ref{Theorem: Hardness of CI testing} is a weaker result than Theorem 2 of \cite{shah2020hardness}. 
		\item On the other hand, the proof of Theorem~\ref{Theorem: Hardness of CI testing} is much simpler than that of the hardness result of \cite{shah2020hardness}, being highly motivated by the recent impossibility results in distribution-free conditional predictive inference~e.g.,~Lemma A.1 of \cite{foygel2019limits} and Lemma 1 of \cite{barber2020distribution}. The key idea of the proof is to introduce ghost samples $\{X_i,Y_i,Z_i\}_{i=1}^J$ and express the type II error of $\phi$ as the iterative expectations associated with sampling without replacement from these ghost samples. When $Z_1,\ldots,Z_J$ are distinct, random draws from these ghost samples with replacement can be viewed as random draws under the null. We note that such argument via sampling with replacement can be traced back to \cite{gretton2012kernel} (see Example 1) who provide a negative result for two-sample testing. Given this key observation, we can connect the type II error with the significance level $\alpha$, and the result follows by a union bound along with the total variation distance between sampling with and without replacement. The details can be found in Appendix~\ref{Section: Proof of Theorem: Hardness of CI testing}.  
	\end{itemize}
\end{remark}

In summary, the hardness result of \cite{shah2020hardness} and Theorem~\ref{Theorem: Hardness of CI testing} indicate that CI testing is a difficult task without further assumptions. This negative result naturally motivates us to explore reasonable conditions under which CI testing, especially based on the local permutation procedure, is feasible. This is the main topic of the next section.

\section{Validity under smoothness conditions} \label{Section: Validity under smoothness conditions}
The goal of this section is two-fold: one is to identify universal conditions under which any local permutation test based on a binned statistic is asymptotically valid; the other is to show that these conditions are tight under certain smoothness conditions in the sense that there exists a permutation test whose type I error is not controlled even asymptotically when our conditions are violated. We start by introducing our smoothness assumptions (Section~\ref{Section: Smoothness conditions}) and then state the main results (Section~\ref{Section: Validity} and Section~\ref{Section: Tightness}).

\subsection{Smoothness conditions} \label{Section: Smoothness conditions}
The validity of the local permutation test crucially relies on the smoothness assumption on conditional distributions. For instance, suppose that $p_{X|Z}(x|z)$ and $p_{Y|Z}(y|z)$ are constant with respect to $z$ for all $z \in \mathcal{Z}$. In this case, $X$ and $Y$ are independent of $Z$ and, consequently, CI testing is the same as unconditional independence testing for which the (local) permutation test has finite-sample validity. Motivated by this observation, we consider similar (but more general) smoothness classes to those in \cite{neykov2020minimax} defined as follows. 
\begin{definition}[$\gamma$-Hellinger Lipschitzness]  \label{Definition: Hellinger Lipschitzness}
	Let $\mathcal{P}_{0,\text{\emph{H}},\gamma,\delta}(L) \subset \mathcal{P}_0$ be the collection of distributions $P_{x,y,z}$ such that for all $z,z' \in \mathcal{Z}$,
	\begin{align*}
		\mathcal{D}_{\gamma\text{,\emph{H}}}(P_{X|Z=z},P_{X|Z=z'}) \leq L \delta (z,z') \quad \text{and} \quad \mathcal{D}_{\gamma\text{,\emph{H}}}(P_{Y|Z=z},P_{Y|Z=z'}) \leq L \delta (z,z'),
	\end{align*}
	where $\delta(z,z')$ is a distance between $z$ and $z'$ in $\mathcal{Z}$. 
\end{definition}

Another smoothness assumption is made based on R{\'e}nyi divergence. 

\begin{definition}[$\gamma$-R{\'e}nyi Lipschitzness]
	Let $\mathcal{P}_{0,\text{\emph{R}},\gamma,\delta}(L) \subset \mathcal{P}_0$ be the collection of distributions $P_{x,y,z}$ such that for all $z,z' \in \mathcal{Z}$,
\begin{align*}
	\mathcal{D}_{\gamma\text{,\emph{R}}}^{1/2}(P_{X|Z=z} \| P_{X|Z=z'}) \leq L \delta (z,z') \quad \text{and} \quad \mathcal{D}_{\gamma\text{,\emph{R}}}^{1/2}(P_{Y|Z=z} \| P_{Y|Z=z'}) \leq L \delta (z,z'),
\end{align*}
where $\delta(z,z')$ is a distance between $z$ and $z'$ in $\mathcal{Z}$. 
\end{definition}

For both Lipschitz conditions, we let $\delta(\cdot,\cdot)$ be the Euclidean distance by default when $Z \in \mathbb{R}^d$. With these smoothness conditions in place, the next section studies the asymptotic validity of the local permutation test.

\subsection{Validity} \label{Section: Validity}
To state the validity result, we begin with additional notation. Let $(X^n,Y^n,\widetilde{Z}^n)$ be the binned version of $(X^n,Y^n,Z^n)$ where $\widetilde{Z} \in \{1,\ldots,M\}$ is a discrete random variable with probability $\mP(\widetilde{Z} = m) = \mP(Z \in B_m)$ for $m=1,\ldots,M$. For simplicity, let us write $q_{\widetilde{Z}}(m) = \mP(\widetilde{Z} = m)$ and denote the conditional distribution of $Z|Z \in B_m$ by $\widetilde{P}_{Z|Z \in B_m}$. Then $(X,Y,\widetilde{Z})$ has its density function $q_{X,Y,\widetilde{Z}}(x,y,m) = q_{XY|\widetilde{Z}}(x,y|m) q_{\widetilde{Z}}(m)$ where 
\begin{align*}
	q_{XY|\widetilde{Z}}(x,y|m) = \int_{B_m} p_{X,Y|Z}(x,y|z) d\widetilde{P}_{Z|Z \in B_m}(z),
\end{align*}
and we denote the corresponding joint distribution by $Q_{X,Y,\widetilde{Z}}$.

Under the considered binning scheme, we make a key observation that the test statistic $T_{\mathrm{CI}}$ is defined only through the binned data, which means that $T_{\mathrm{CI}}$ based on $(X^n,Y^n,Z^n)$ is equal in distribution to that based on $(X^n,Y^n,\widetilde{Z}^n)$. Furthermore, since our test function $\phi_{\mathrm{perm},n}$ is computed only through the binned data, we observe that 
\begin{align} \label{Eq: equivalence of two permutation tests}
	\mE_{P_{X,Y,Z}^n}[\phi_{\mathrm{perm},n}] = \mE_{Q_{X,Y,\widetilde{Z}}^n}[\phi_{\mathrm{perm},n}].
\end{align}
To proceed further, we consider a product density $\widetilde{q}_{X,Y,\widetilde{Z}}(x,y,m) = q_{X\cdot |\widetilde{Z}}(x|m) q_{\cdot Y|\widetilde{Z}}(y|m)q_{\widetilde{Z}}(m)$ where $q_{X\cdot |\widetilde{Z}}(x|m)$ and $q_{\cdot Y|\widetilde{Z}}(y|m)$ are the marginals of $q_{XY|\widetilde{Z}}(x,y|m)$, and denote the corresponding joint distribution by $\widetilde{Q}_{X,Y,\widetilde{Z}}$. By construction, $(X^n,Y^n,\widetilde{Z}^n)$ from $\widetilde{Q}_{X,Y,\widetilde{Z}}$ satisfies $X \independent Y | \widetilde{Z}$ and therefore exchangeability of $Y$ within each bin yields
\begin{align*}
	\mE_{\widetilde{Q}_{X,Y,\widetilde{Z}}^n}[\phi_{\mathrm{perm},n}] \leq \alpha.
\end{align*}
Combining the above inequality with the identity~(\ref{Eq: equivalence of two permutation tests}) yields a generic type I error bound for the local permutation test in terms of the total variation distance. 
\begin{lemma}[Type I error bound in terms of the TV distance] \label{Lemma: Error bound in terms of the TV distance}
	Suppose that the distribution $P_{X,Y,Z}$ belongs to $\mathcal{P}_0$. Then for any $\alpha \in (0,1)$, the type I error of $\phi_{\text{\emph{perm}},n}$ is bounded above by
	\begin{align*}
		& \mE_{P_{X,Y,Z}^n}[\phi_{\text{\emph{perm}},n}] \leq \alpha + \mathcal{D}_{\emph{TV}}\big(Q_{X,Y,\widetilde{Z}}^n,\widetilde{Q}_{X,Y,\widetilde{Z}}^n\big).
	\end{align*}
\end{lemma}
The given bound implies that the local permutation test is valid when the binned null distribution is indistinguishable from its CI projection. The proof of this result follows by the definition of the TV distance. It is worth pointing out that since the randomized permutation test $\phi_{\mathrm{perm},n,a}$ (Remark~\ref{Remark: MC approximation/randomization}) is exact under the law of $\widetilde{Q}_{X,Y,\widetilde{Z}}^n$, one can establish a stronger result that 
\begin{align*}
	\big| \mE_{P_{X,Y,Z}^n,R}[\phi_{\mathrm{perm},n,a}] - \alpha \big| \leq \mathcal{D}_{\mathrm{TV}}\big(Q_{X,Y,\widetilde{Z}}^n,\widetilde{Q}_{X,Y,\widetilde{Z}}^n\big),
\end{align*}
where $R \sim \text{Uniform}[0,1]$ is for randomization and independent of the data. In either randomized or non-randomized test, our main task boils down to identifying reasonable conditions under which the TV distance between $Q_{X,Y,\widetilde{Z}}^n$ and $\widetilde{Q}_{X,Y,\widetilde{Z}}^n$ tends to zero asymptotically. In general, however, it is challenging to directly work with the TV distance between two product measures. Instead we upper bound the TV distance by the Hellinger distance as follows.

\begin{lemma}[Type I error bound in terms of the Hellinger distance] \label{Lemma: Error bound in terms of the Hellinger distance}
	Suppose that the distribution $P_{X,Y,Z}$ belongs to $\mathcal{P}_0$. Then for any $\alpha \in (0,1)$, the type I error of $\phi_{\text{\emph{perm}},n}$ is bounded above by
	\begin{align} \label{Eq: Error bound in terms of the Hellinger distance}
		& \mE_{P_{X,Y,Z}^n}[\phi_{\text{\emph{perm}},n}] \leq \alpha + \Bigg\{2n \sum_{m=1}^M q_{\widetilde{Z}}(m)  \times \mathcal{D}_{\emph{H}}^2\big(Q_{X,Y|\widetilde{Z}=m},\widetilde{Q}_{X,Y|\widetilde{Z}=m}\big) \Bigg\}^{1/2},
	\end{align}
	where the second term on the right-hand side is simply $\sqrt{2n}\mathcal{D}_{\mathrm{H}}\big(Q_{X,Y,\widetilde{Z}},\widetilde{Q}_{X,Y,\widetilde{Z}}\big)$.
\end{lemma}
We mention that, while the Hellinger bound~(\ref{Eq: Error bound in terms of the Hellinger distance}) is provably looser than that based on the TV distance, it has the same characterization as the TV bound in terms of when the local permutation is asymptotically valid. More specifically, using well-known bounds relating the TV distance and the Hellinger distance, it can be verified that 
\begin{align*}
	\text{$\mathcal{D}_{\mathrm{TV}}\big(Q_{X,Y,\widetilde{Z}}^n,\widetilde{Q}_{X,Y,\widetilde{Z}}^n\big) \rightarrow 0$ if and only if $n \mathcal{D}_{\mathrm{H}}^2\big(Q_{X,Y,\widetilde{Z}},\widetilde{Q}_{X,Y,\widetilde{Z}}\big) \rightarrow 0$.}
\end{align*}
We are now ready to state the main results of this section. To this end, we denote the maximum diameter of bins $B_1,\ldots,B_M$ by $h$, that is
\begin{align} \label{Eq: maximum diameter}
	h := \max_{1 \leq m \leq M} \sup_{z,z' \in B_m} \delta(z,z').
\end{align}
We first state the result under $\gamma$-Hellinger Lipschitzness in the next theorem and then consider $\gamma$-R{\'e}nyi Lipschitzness in Theorem~\ref{Theorem: Validity of the permutation test under Renyi smooth}. 
\begin{theorem}[Validity of $\phi_{\mathrm{perm},n}$ under $\gamma$-Hellinger Lipschitzness] \label{Theorem: Validity of the permutation test under Hellinger smooth}
	For any $\alpha \in (0,1)$, the type I error of $\phi_{\text{\emph{perm}},n}$ under $\gamma$-Hellinger Lipschitzness is bounded above by 
	\begin{equation}
	 \begin{aligned} \label{Eq: upper bounds under Hellinger Lipschitzness}
	 	\sup_{P_{X,Y,Z} \in \mathcal{P}_{0, \text{\emph{H}},\gamma,\delta}(L)}  \mE_{P_{X,Y,Z}^n}[\phi_{\text{\emph{perm}},n}] ~\leq
	 	\begin{cases}
	 		\alpha + C_{\gamma}n^{1/2} L^{\gamma} h^{\gamma}, \quad & \text{if $\gamma \in [1,2]$,}\\[.5em]	
	 		\alpha + C_{\gamma} n^{1/2} L^{2} h^{2}, \quad & \text{if $\gamma >2$,} 	 	
 		\end{cases}
	 \end{aligned}
 	\end{equation}
 where $C_{\gamma}$ is a constant that only depends on $\gamma$. 
\end{theorem}

A few remarks are in order.
\begin{remark} \leavevmode \normalfont
	\begin{itemize}
		\item The result of Theorem~\ref{Theorem: Validity of the permutation test under Hellinger smooth} shows that once we assume that $L$ and $\gamma$ are fixed in the sample size, a sufficient condition for the validity of $\phi_{\mathrm{perm},n}$ is $h = o(n^{-1/2\gamma})$ if $\gamma \in [1,2]$ and $h = o(n^{-1/4})$ if $\gamma > 2$ under $\gamma$-Hellinger Lipschitzness. Since users have control over $h$, Theorem~\ref{Theorem: Validity of the permutation test under Hellinger smooth} provides a guideline for the choice of binning that ensures type I error control of the local permutation test. We also note that when $h$ is too small, most of the bins are empty, which adversely affects the power performance. In other words, there is an intriguing trade-off between the type I error and power in terms of the choice of $h$. We discuss this trade-off more in Section~\ref{Section: Power analysis and optimality}.
		\item In most practical applications, the smoothness parameter $\gamma$ is unknown. In such case, one can choose $h$ such that $n^{1/2} h \rightarrow 0$ as $n \rightarrow \infty$, which leads to an asymptotically valid permutation test for all $\gamma \geq 1$, assuming other parameters are fixed. However, as mentioned before, choosing small $h$ comes at a price of low power under the alternative. It would be interesting to see whether an adaptive way of choosing $h$ is possible without much loss of power. We leave this direction to future work.
		\item We observe an interesting phenomenon that there exists a sharp transition at $\gamma = 2$, which corresponds to the Hellinger distance. In particular, the result illustrates that the smoothness condition beyond $\gamma > 2$ does not really help to improve the convergence rate of the type I error. Importantly, the given upper bound~(\ref{Eq: upper bounds under Hellinger Lipschitzness}) is tight in terms of $n$ and $h$ in some cases. More specifically, we show in Section~\ref{Section: Tightness} that there exists a local permutation test whose type I error rate can be made arbitrarily large unless the upper bound~(\ref{Eq: upper bounds under Hellinger Lipschitzness}) converges to $\alpha$. 
		\item The proof of Theorem~\ref{Theorem: Validity of the permutation test under Hellinger smooth} builds on Lemma~\ref{Lemma: Error bound in terms of the Hellinger distance} and the monotonicity property of the generalized Hellinger distance. We note that Lemma~\ref{Lemma: Error bound in terms of the Hellinger distance} has a bound in terms of smoothed distributions over partitions, whereas $\gamma$-Hellinger Lipschitzness is stated in terms of original distributions. The bulk of the effort in proving Theorem~\ref{Theorem: Validity of the permutation test under Hellinger smooth} lies in connecting the Hellinger distance between $Q_{X,Y|\widetilde{Z}=m}$ and $\widetilde{Q}_{X,Y|\widetilde{Z}=m}$ to the $\gamma$-Hellinger distance between $P_{X|Z=z}$ and $P_{X|Z=z'}$ (and also between $P_{Y|Z=z}$ and $P_{Y|Z=z'}$). The details can be found in Appendix~\ref{Section: Proof of Theorem: Validity of the permutation test under Hellinger smooth}. 
	\end{itemize}
\end{remark}

Next we present a similar result under $\gamma$-R{\'e}nyi Lipschitzness.

\begin{theorem}[Validity of $\phi_{\mathrm{perm},n}$ under $\gamma$-R{\'e}nyi Lipschitzness] \label{Theorem: Validity of the permutation test under Renyi smooth}
	For any $\alpha \in (0,1)$ and $\gamma>0$, the type I error of $\phi_{\text{\emph{perm}},n}$ under $\gamma$-R{\'e}nyi Lipschitzness is bounded above by 
	\begin{align} \label{Eq: upper bounds under Renyi Lipschitzness}
		\sup_{P_{X,Y,Z} \in \mathcal{P}_{0, \text{\emph{R}},\gamma,\delta}(L)}  \mE_{P_{X,Y,Z}^n}[\phi_{\text{\emph{perm}},n}] ~\leq \alpha + C_{\gamma} n^{1/2} L^{2} h^{2},
	\end{align}
	where $C_{\gamma}$ is a constant that only depends on $\gamma$. 
\end{theorem}
In contrast to Theorem~\ref{Theorem: Validity of the permutation test under Hellinger smooth}, the above result indicates that the smoothness parameter $\gamma$ in R{\'e}nyi Lipschitzness does not affect the type I error of $\phi_{\mathrm{perm},n}$ by more than a constant factor. At a high-level, we observe this phenomenon because R{\'e}nyi divergence is lower bounded by the squared Hellinger distance, up to a constant factor, for any $\gamma > 0$ (see Lemma~\ref{Lemma: Renyi divergence}). In other words, the conditional distribution in $\mathcal{P}_{0, \text{R},\gamma,\delta}(L)$ is at least as smooth as the one in $\mathcal{P}_{0, \mathrm{H},\gamma=2,\delta}(L)$, which means that we are essentially in the second regime of Theorem~\ref{Theorem: Validity of the permutation test under Hellinger smooth} for $\gamma \geq 2$. Indeed, it should be clear from the proof that the same upper bound in Theorem~\ref{Theorem: Validity of the permutation test under Renyi smooth} holds for any Lipschitzness condition whose underlying divergence is lower bounded by the squared Hellinger distance such as KL divergence and $\chi^2$ divergence.

\begin{remark}[Poissonization] \label{Remark: Poissonization} \normalfont
To analyze the power of binning-based tests, it is often convenient to assume that the sample size has a Poisson distribution \citep[e.g.][]{canonne2018testing,balakrishnan2019hypothesis,neykov2020minimax}. This, so-called, Poissonization trick allows us to bypass the difficulty in dealing with the dependence between different bins. In fact, as we proved in Proposition~\ref{Proposition: Poissonization} in Appendix~\ref{Section: Poissonization}, the local permutation test under Poissonization has the same validity as before in Theorem~\ref{Theorem: Validity of the permutation test under Hellinger smooth} and Theorem~\ref{Theorem: Validity of the permutation test under Renyi smooth}. To explain it briefly, we shall use the convenient notation
\begin{align*}
	\mE_{P_{X,Y,Z}^N,N}[\cdot]  := \sum_{k=0}^\infty \mP(N = k)  \mE_{P_{X,Y,Z}^k} [\cdot] 
\end{align*}
to denote the expectation operator with respect to $P_{X,Y,Z}^N$ where $N$ is a random sample size. We similarly write $\mP_{P_{X,Y,Z}^N,N}[\cdot]$ to denote $\mE_{P_{X,Y,Z}^N,N}[\mathds{1}(\cdot)]$. Suppose now that $\sup_{P_{X,Y,Z} \in \mathcal{P}_{0}'} \mE_{P_{X,Y,Z}^n}[\phi_{\mathrm{perm},n}] \leq \alpha + Cn^{-\epsilon}$ for $\mathcal{P}_{0}' \subset \mathcal{P}_0$ and some constants $C,\epsilon>0$. Then under the same condition, the permutation test under Poissonization (i.e.~$N\sim \mathrm{Pois}(n)$) satisfies $\sup_{P_{X,Y,Z} \in \mathcal{P}_{0}'} \mE_{P_{X,Y,Z}^N,N}[\phi_{\mathrm{perm},N}] \leq \alpha + C_\epsilon n^{-\epsilon}$ where $C_\epsilon$ is a constant that only depends on $\epsilon$. See Proposition~\ref{Proposition: Poissonization} for a more rigorous statement. 
\end{remark}

We now move to the next section where we provide complementary results of this section.

\subsection{Lower bounds} \label{Section: Tightness}
In this section, we demonstrate that the upper bounds for the type I error established in Section~\ref{Section: Validity} cannot be improved further in some cases. In particular, we prove that there exists a local permutation test whose type I error cannot be controlled if one chooses $h$ in such a way that the upper bounds~(\ref{Eq: upper bounds under Hellinger Lipschitzness}) or (\ref{Eq: upper bounds under Renyi Lipschitzness}) diverge. In order to simplify our presentation, we focus on the case where $X$ and $Y$ are discrete random variables whereas $Z$ is continuous and bounded between $[0,1]$. Other cases such as multivariate continuous $(X,Y,Z)$ will be discussed in Remark~\ref{Remark: Lower bounds}.

Suppose that $X$ and $Y$ are discrete random variables supported on $[\ell_1] \times [\ell_2]$ for some positive integers $\ell_1$ and $\ell_2$. By convention, $[\ell_1]$ denotes the set of integers $\{1,\ldots,\ell_1\}$ and $[\ell_2]$ is similarly defined. Let $\{B_1,\ldots,B_M\}$ be an equi-partition of $[0,1]$ so that the length of each bin is $h = M^{-1}$. The given partition yields the binned data sets $\bW_1,\ldots,\bW_M$ defined in Section~\ref{Section: Local permutation procedure}. To study lower bounds, we work with the weighted sum of U-statistics proposed by \cite{canonne2018testing} and \cite{neykov2020minimax}. Let 
\begin{align*}
	\psi_{ij}^m(x,y) = \mathds{1}(X_{i,m}=x,Y_{i,m}=y) -  \mathds{1}(X_{i,m}=x) \mathds{1}(Y_{j,m}=y),
\end{align*}
and define a kernel function as
\begin{align*} 
	h_{i_1,i_2,i_3,i_4}^m= \frac{1}{4!} \sum_{\pi \in \boldsymbol{\Pi}_4} \sum_{x \in [\ell_1], y \in [\ell_2]} \psi^m_{\pi(1) \pi(2)}(x,y) \psi^m_{\pi(3)\pi(4)}(x,y),
\end{align*}
where $\boldsymbol{\Pi}_4$ is the set of all permutations of $\{i_1,i_2,i_3,i_4\}$. By linearity of expectations, it is seen that $h_{i_1,i_2,i_3,i_4}^m$ is an unbiased estimator of the squared $L_2$ norm between $Q_{X,Y|\widetilde{Z}=m}$ and $\widetilde{Q}_{X,Y|\widetilde{Z}=m}$. Given this kernel and by recalling $\bW_m = \{(X_{1,m},Y_{1,m}),\ldots, (X_{\sigma_m,m},Y_{\sigma_m,m})\}$, the resulting U-statistic is calculated as
\begin{align} \label{Eq: fourth order U-statistic}
	U(\bW_m):= \frac{1}{\binom{\sigma_m}{4}} \sum_{i_1<i_2<i_3<i_4: (i_1,i_2,i_3,i_4) \in[\sigma_m]} h_{i_1,i_2,i_3,i_4}^m.
\end{align}
The final statistic is a weighted sum of $U(\bW_1), \ldots, U(\bW_M)$ given by
\begin{align} \label{Eq: unweighted statistic}
	T_{\mathrm{CI}} := \sum_{m \in [M]} \mathds{1}(\sigma_m \geq 4) \sigma_m U(\bW_m).
\end{align}
Several properties of $T_{\mathrm{CI}}$, such as minimax power optimality, have been studied under Poissonization by \cite{canonne2018testing} and \cite{neykov2020minimax}. To fully benefit from their results, we work with a modified local permutation test: First draw $N \sim \text{Pois}(n/2)$ and accept the null when $N > n$. If $N \leq n$, we carry out a local permutation test with $N$ samples randomly chosen from $(X^n,Y^n,Z^n)$. Formally, we define the modified local permutation test as
\begin{align} \label{Eq: modified local permutation test}
	\phi_{\mathrm{perm},n}^\dagger := \phi_{\mathrm{perm},N} \times \mathds{1}(N \leq n),
\end{align}
where $\phi_{\mathrm{perm},N}$ denotes the local permutation test using $T_{\mathrm{CI}}$~(\ref{Eq: unweighted statistic}) computed based on the $N$ samples chosen before. By Proposition~\ref{Proposition: Poissonization} along with the inequality $\phi_{\mathrm{perm},n}^\dagger \leq \phi_{\mathrm{perm},N} $, it is clear that $\phi_{\mathrm{perm},n}^\dagger$ is asymptotically valid whenever the upper bounds in~(\ref{Eq: upper bounds under Hellinger Lipschitzness}) and (\ref{Eq: upper bounds under Renyi Lipschitzness}) converge to $\alpha$. The next theorem provides complementary results establishing lower bounds for the type I error of $\phi_{\mathrm{perm},n}^\dagger$.
\begin{theorem}[Lower bounds] \label{Theorem: Lower bounds}
	For arbitrarily small but fixed $\alpha, \beta >0$, choose $h_n \rightarrow 0$ such that $\sqrt{n} h_n^\gamma \rightarrow \infty$ for $\gamma \in [1,2]$ and $\sqrt{n} h_n^2 \rightarrow \infty$ for $\gamma > 2$ under $\gamma$-Hellinger Lipschitzness. On the other hand, choose $h_n \rightarrow 0$ such that $\sqrt{n} h_n^2 \rightarrow \infty$ for $\gamma > 0$ under $\gamma$-R{\'e}nyi Lipschitzness.
	Consider $\phi_{\text{\emph{perm}},n}^\dagger$ based on $T_{\text{\emph{CI}}}$ in~(\ref{Eq: unweighted statistic}). Then there exist constants $n_0,L>0$, where $L$ depends only on $\gamma$, such that for all $n \geq n_0$, the following two inequalities hold:
	\begin{align*}
		& \sup_{P_{X,Y,Z} \in \mathcal{P}_{0, \text{\emph{H}},\gamma,\delta}(L)}  \mE_{P_{X,Y,Z}^N,N}[\phi_{\text{\emph{perm}},n}^\dagger] \geq  1 - \beta \quad \text{and} \\[.5em]
		& \sup_{P_{X,Y,Z} \in \mathcal{P}_{0, \text{\emph{R}},\gamma,\delta}(L)}  \mE_{P_{X,Y,Z}^N,N}[\phi_{\text{\emph{perm}},n}^\dagger] \geq  1 - \beta.
	\end{align*}
\end{theorem}
Let us provide several comments on Theorem~\ref{Theorem: Lower bounds}. 
\begin{remark} \leavevmode \normalfont \label{Remark: Lower bounds}
	\begin{itemize}
		\item A crucial observation is that the type I error of $\phi_{\mathrm{perm},n}^\dagger$ for CI testing corresponds to its power for testing
		\begin{align} \label{Eq: modified hypothesis}
			\widetilde{H}_0: Q_{X,Y,\widetilde{Z}} = \widetilde{Q}_{X,Y,\widetilde{Z}} \quad \text{versus} \quad 	\widetilde{H}_1: Q_{X,Y,\widetilde{Z}} \neq \widetilde{Q}_{X,Y,\widetilde{Z}}.
		\end{align}
		This means that, in order to verify that the type I error of $\phi_{\mathrm{perm},n}^\dagger$ is inflated, it suffices to show that $\phi_{\mathrm{perm},n}^\dagger$ is asymptotically powerful against the above alternative~(\ref{Eq: modified hypothesis}). With this observation in place, our main task is to construct a distributional setting where $\phi_{\mathrm{perm},n}^\dagger$ is able to distinguish the binned distribution and its CI projection with high probability. In fact, the conditions of Theorem~\ref{Theorem: Lower bounds} guarantee that these two distributions are far enough for  $\phi_{\mathrm{perm},n}^\dagger$ to be asymptotically powerful. 
		\item In order to ease our analysis, we carefully design the distribution of $Z$ such that most samples are observed in one of the partitions with high probability. In this case, the test statistic $T_{\mathrm{CI}}$ approximates $\sigma_1 U(\bW_1)$, which is much easier to handle. It is then sufficient to study the permutation test based on $\sigma_1 U(\bW_1)$ and prove that it is asymptotically powerful under the given conditions. We show that this is indeed the case by building on the results of \cite{kim2020minimax} where the authors investigate the permutation test based on $U(\bW_1)$ for unconditional independence testing. 
		\item We expect that $\phi_{\mathrm{perm},n}$ (i.e.~without Poissionization) also achieves the same error bounds in Theorem~\ref{Theorem: Lower bounds} as it always uses more samples than $\phi_{\mathrm{perm},n}^\dagger$ by definition. See empirical evidence in Figure~\ref{Figure: Double-binning}. However, we found it challenging to analyze $T_{\mathrm{CI}}$ without Poissonization, especially its variance, due to a non-trivial dependence between the summands of $T_{\mathrm{CI}}$. Due to this technical difficulty, we focus on the Poisson-sampling scheme as in \cite{canonne2018testing} and \cite{neykov2020minimax} and leave the detailed analysis of $\phi_{\mathrm{perm},n}$ to future work. Nevertheless, the concentration property of a Poisson random variable allows us to say a certain negative result on $\phi_{\mathrm{perm},n}$ without Poissonization. Specifically, note that a Poisson random variable $N$ with parameter $n/2$ is bound between $cn \leq N \leq Cn$ with high probability where $c,C>0$ some positive constants~\citep[e.g.][]{canonne2019note}. Thus it is guaranteed that one can find a fixed sample size $\widetilde{n}$ such that $c n \leq \widetilde{n} \leq C n$ and Theorem~\ref{Theorem: Lower bounds} holds for $\phi_{\mathrm{perm},\widetilde{n}}$ without Poissonization.
		\item For simplicity, we prove Theorem~\ref{Theorem: Lower bounds} using the example where $(X,Y)$ are discrete and $Z$ is a univariate continuous random variable. Nevertheless, our proof can be extended to the case of multivariate continuous $(X,Y,Z)$ as follows. Consider piecewise constant densities for $(X,Y)$ and assume that the components $Z^{(2)},\ldots,Z^{(d_Z)}$ of $Z \in \mathbb{R}^{d_Z}$ are independent of the rest of variables. In this case, we are essentially in the setting where $X,Y$ are discrete and $Z$ is univariate. Therefore, the same proof carries through except now that the maximum diameter $h$ depends on $d_Z$. When an equi-partition is considered, we note that $d_Z$ only affects the scaling factor in $h$ and hence the statement of Theorem~\ref{Theorem: Lower bounds} remains true, provided that $d_Z$ is fixed. We also note that the upper bound results~(Theorem~\ref{Theorem: Validity of the permutation test under Hellinger smooth} and Theorem~\ref{Theorem: Validity of the permutation test under Renyi smooth}) are stated in terms of the maximum diameter $h$; thereby the upper bounds remain the same for both univariate and multivariate cases of $Z$.
	\end{itemize}
\end{remark}

So far we have explored type I error control of the local permutation test. Next we turn our attention to the power and study its optimality in certain regimes.

\section{Power analysis} \label{Section: Power analysis and optimality}

This section considers both discrete and continuous cases of the conditional variable $Z$ and investigates the power property of local permutation tests. In order to achieve meaningful power, we focus on a subset of alternatives, which are at least $\varepsilon$ far away from the null in terms of the TV distance. Our main interest is then to characterize $\varepsilon$ for which local permutation tests can be powerful.

\subsection{Discrete $(X,Y,Z)$}
To start with the discrete case where $(X,Y,Z) \in [\ell_1] \times [\ell_2] \times [M]$, we revisit the test statistics proposed by \cite{canonne2018testing} and demonstrate the power property of the local permutation procedure based on the same test statistics. \cite{canonne2018testing} propose two test statistics for CI testing. The first one is $T_{\mathrm{CI}}$ given in (\ref{Eq: unweighted statistic}), which is defined as the sum of unweighted U-statistics. While $T_{\mathrm{CI}}$ is simple and performs well in certain regimes, it may suffer from a large variance especially when the dimensions $\ell_1$ and $\ell_2$ of $X$ and $Y$ are large. To mitigate this issue, \cite{canonne2018testing} propose another statistic building on the flattening idea of \cite{diakonikolas2016new}. The latter statistic can be viewed as the sum of weighted U-statistics as observed by \cite{neykov2020minimax}.

\paragraph{Weighted U-statistic.} To proceed, let us formally write down the weighted U-statistic. First recall that $\bW_m =\{(X_{1,m},Y_{1,m}),\ldots,(X_{\sigma_m,m},Y_{\sigma_m,m})\}$ is the set of pairs of $(X_i,Y_i)$ with $Z_i = m$. Suppose that the sample size of $\bW_m$ is $\sigma_m \geq 4$ and $\sigma_m = 4 + 4t_m$ for some $t_m \in \mathbb{N}$. Following the notation in \cite{neykov2020minimax}, let $t_{1,m} := \min\{t_m, \ell_1\}$ and $t_{2,m} := \min\{t_m, \ell_2\}$. We then randomly split the data $\bW_m$ into three sets $\bW_{X,m}$, $\bW_{Y,m}$ and $\bW_{XY,m}$ of size $t_{1,m}$, $t_{2,m}$ and $2{t_m}+4$, respectively, where $\bW_{X,m}:=\{X_{i,m}:i \in [t_{1,m}]\}$, $\bW_{Y,m} :=\{Y_{i,m}: t_{1,m} + 1\leq i \leq t_{1,m} + t_{2,m}\}$ and $\bW_{XY,m}:=\{(X_{i,m},Y_{i,m}: 2t_m + 1 \leq i \leq \sigma_m)\}$. The purposes of these splits are different: the first two will be used to compute weights and the last one will be used to compute the U-statistic. In particular, as a weight function, consider a positive integer $1 + a_{xy,m} = (1 + a_{x,m}) (1 + a_{y,m})$ where $a_{x,m}$ is the number of occurrences $x$ in $\bW_{X,m}$ and, similarly, $a_{y,m}$ is the number of occurrences $y$ in $\bW_{Y,m}$. Next let $h_{i_1,i_2,i_3,i_4}^{m,\boldsymbol{a}}$ denote a weighted kernel function defined as
\begin{align*}
	h_{i_1,i_2,i_3,i_4}^{m,\boldsymbol{a}}:= \frac{1}{4!} \sum_{\pi \in \boldsymbol{\Pi}_4} \sum_{x \in [\ell_1], y \in [\ell_2]} \frac{\psi^m_{\pi(1) \pi(2)}(x,y) \psi^m_{\pi(3)\pi(4)}(x,y)}{1+a_{xy,m}}.
\end{align*}
Given this kernel, we compute the weighted U-statistic for each $1 \leq m \leq M$ as
\begin{align*}
	U_W(\bW_m) := \frac{1}{\binom{2t_m+4}{2}} \sum_{i_1<i_2<i_3<i_4: (i_1,i_2,i_3,i_4) \in \bW_{XY,m}} h_{i_1,i_2,i_3,i_4}^{m,\boldsymbol{a}}, 
\end{align*}
where $(i_1,i_2,i_3,i_4) \in \bW_{XY,m}$ stands for taking four observations from $\bW_{XY,m}$. Now, by letting $\omega_m:=\sqrt{\min(\sigma_m,\ell_1)\min(\sigma_m,\ell_2)}$, the final test statistic for CI is defined as a weighted sum of $U_W(\bW_1),\ldots,U_W(\bW_M)$ given by
\begin{align} \label{Eq: weighted statistic}
	T_{\mathrm{CI},W} := \sum_{m \in [M]} \mathds{1}(\sigma_m \geq 4) \sigma_m \omega_m U_W(\bW_m).
\end{align}

\paragraph{Tests of \cite{canonne2018testing}.} For both test statistics $T_{\mathrm{CI}}$ in~(\ref{Eq: unweighted statistic}) and $T_{\mathrm{CI},W}$ in~(\ref{Eq: weighted statistic}), \cite{canonne2018testing} suggest that one rejects the null when the test statistic is larger than $\zeta \sqrt{\min(n,M)}$ where $\zeta$ is a sufficiently large (but unspecified) constant. This cutoff value can be roughly understood as an upper bound of the standard deviation of the test statistic under the null. For ease of reference, we let $\phi_{\mathrm{CDKS},1} := \mathds{1}(T_{\mathrm{CI}} \geq \zeta  \sqrt{\min(n,M)})$ denote the test based on the unweighted statistic and similarly let $\phi_{\mathrm{CDKS},2} := \mathds{1}(T_{\mathrm{CI},W} \geq \zeta  \sqrt{\min(n,M)})$ denote the test based on the weighted test statistic computed based on $N \sim \mathrm{Pois}(n)$ samples. To describe their power results, let $\mathcal{P}_{[M]}$ be the set of discrete distributions defined on the support $[\ell_1] \times [\ell_2] \times [M]$. Moreover, let $\mathcal{P}_{0,[M]} \subset \mathcal{P}_{[M]}$ where $X \independent Y|Z$ and $\mathcal{P}_{1,[M]} := \mathcal{P}_{[M]} \backslash \mathcal{P}_{0,[M]}$. \cite{canonne2018testing} first consider the regime where $\ell_1$ and $\ell_2$ are fixed and $\varepsilon$ satisfies a certain condition recalled in (\ref{Eq: condition on varepsilon 2}) of the Appendix, and then show that 
\begin{align} \label{Eq: CKDS power result 1}
	\sup_{P_{X,Y,Z} \in \mathcal{P}_{1,[M]}: \inf_{Q \in \mathcal{P}_{0,[M]}} \mathcal{D}_{\mathrm{TV}}(P_{X,Y,Z},Q) \geq \varepsilon} \mE_{P_{X,Y,Z}^N,N}[1 - \phi_{\mathrm{CDKS},1}] \leq \frac{1}{100}.
\end{align}
For example, the condition for $\varepsilon$ is fulfilled when $\varepsilon \geq c \max\{M^{1/4}/n^{1/2}, M^{7/8}/n, M^{3/4}/n^{7/8}\}$ for some large constant $c>0$. In the second regime where $\ell_1$ and $\ell_2$ can vary, the authors consider a more involved condition for $\varepsilon$ depending on $(\ell_1,\ell_2,M,n)$, and show that 
\begin{align} \label{Eq: CKDS power result 2}
	\sup_{P_{X,Y,Z} \in \mathcal{P}_{1,[M]}: \inf_{Q \in \mathcal{P}_{0,[M]}} \mathcal{D}_{\mathrm{TV}}(P_{X,Y,Z},Q) \geq \varepsilon} \mE_{P_{X,Y,Z}^N,N}[1 - \phi_{\mathrm{CDKS},2}] \leq \frac{1}{100}.
\end{align}
The condition for $\varepsilon$ in this second regime is recalled in (\ref{Eq: condition on varepsilon}) of the Appendix for completeness.

\paragraph{Main results for the discrete case.} In contrast to \cite{canonne2018testing}, we consider relatively more practical tests calibrated by the local permutation procedure, which do not rely on unspecified constants. We then argue that the permutation-based tests have the same theoretical guarantee as the tests of \cite{canonne2018testing}. As mentioned earlier, the local permutation test can control the type I error rate in the discrete setting without any further assumptions. Therefore our focus is on the power of the test. In the next remark, we explain a modified local permutation procedure, which we refer to as the ``half-permutation'' procedure, that facilitates the power analysis of the test based on $T_{\mathrm{CI},W}$. 

\begin{remark}[Full- versus.~half-permutation] \leavevmode \normalfont \label{Remark: full- vs half-permutation}
	For the weighted test statistic $T_{\mathrm{CI},W}$, there are two possible ways of calibrating the test via the local permutation procedure. The first one, we call ``\emph{full-permutation}'', computes the $p$-value by permuting all $Y$ labels within $\bW_m$, independently, for each $m$. This is equivalent to the procedure described in Algorithm~\ref{alg:local permutation procedure}. The second one, we call ``\emph{half-permutation}'', only permutes the $Y$ labels within $\bW_{XY,m}$, independently, for each $m$. Both approaches have finite-sample validity but the power of the first approach is intrinsically more difficult to analyze since each permutation destroys the independence structure among $\bW_{X,m}$, $\bW_{Y,m}$ and $\bW_{XY,m}$. On the other hand, the half-permutation approach preserves the independence between $\bW_{XY,m}$ and $\{ \bW_{X,m}, \bW_{Y,m}\}$ even after permutations. Moreover, it has computational advantage over the full-permutation test since we do not need to recompute weights $1 + a_{xy,m}$ for each permutation. A similar strategy has been used in \cite{kim2020minimax,kim2021classification} to analyze two-sample and (unconditional) independence tests.
\end{remark}

We are now ready to state the main results of this subsection. As in \cite{canonne2018testing}, suppose that we draw $N$ i.i.d.~samples from $P_{X,Y,Z} \in \mathcal{P}_{[M]}$ where $N \sim \mathrm{Pois}(n)$. Given these samples, let $\phi_{\mathrm{perm},1}$ be the local permutation test based on the unweighted test statistic $T_{\mathrm{CI}}$~(\ref{Eq: unweighted statistic}) through the full-permutation procedure described in Remark~\ref{Remark: full- vs half-permutation}. Similarly, we let $\phi_{\mathrm{perm},2}$ be the local permutation test base on the weighted test statistic $T_{\mathrm{CI},W}$~(\ref{Eq: weighted statistic}) through the half-permutation procedure described in Remark~\ref{Remark: full- vs half-permutation}. For both tests, we set the significance level $\alpha = 0.01$ for simplicity. These tests have the following guarantee on the type II error rate.

\begin{theorem}[Type II error for discrete $X,Y,Z$] \label{Theorem: Type II error under discrete setting}
	Consider the local permutation tests $\phi_{\mathrm{perm},1}$ and $\phi_{\mathrm{perm},2}$, described above. In the setting of discrete ($X,Y,Z$), $\phi_{\mathrm{perm},1}$ and $\phi_{\mathrm{perm},2}$ have the same type II error guarantee as in (\ref{Eq: CKDS power result 1}) and (\ref{Eq: CKDS power result 2}), respectively. 
\end{theorem}

The implications of Theorem~\ref{Theorem: Type II error under discrete setting} are as follows. 
\begin{remark} \leavevmode \normalfont 
	\begin{itemize}
		\item In Theorem~\ref{Theorem: Type II error under discrete setting}, we set the type I error and the type II error by $1/100$ for simplicity. In fact, $1/100$ can be replaced with an arbitrarily small number by adjusting the constant factor in the condition for $\varepsilon$ given in (\ref{Eq: condition on varepsilon}) of the Appendix.
		\item To make the given tests feasible for a fixed sample size, one can apply the truncation trick as in (\ref{Eq: modified local permutation test}) and consider $\phi_{\mathrm{perm},1}^\dagger = \phi_{\mathrm{perm},1} \times \mathds{1}(N \leq n)$ and $\phi_{\mathrm{perm},2}^\dagger = \phi_{\mathrm{perm},2} \times \mathds{1}(N \leq n)$ where $N \sim \mathrm{Pois}(n/2)$. As discussed before, these modified tests have smaller type I errors than $\phi_{\mathrm{perm},1}$ and $\phi_{\mathrm{perm},2}$ based on $N \sim \mathrm{Pois}(n/2)$, respectively, and have the same power guarantee up to $e^{-n/8}$ factor. See	equation~(\ref{Eq: bound on type II error}) in the Appendix for more details.
		\item \cite{canonne2018testing} further prove that the condition~(\ref{Eq: condition on varepsilon}) in terms of $n$ cannot be improved in certain regimes (depending on $\ell_1,\ell_2,M,\varepsilon$) by providing matching lower bounds. This together with Theorem~\ref{Theorem: Type II error under discrete setting} implies that the corresponding permutation test shares the same rate optimality as \cite{canonne2018testing} whenever the tests of \cite{canonne2018testing} are rate optimal in terms of the sample complexity. Despite the same optimality property, the permutation test may be more attractive than the corresponding test of \cite{canonne2018testing} as it does not depend on an unspecified constant and it tightly controls the type I error rate in finite-sample settings.
		\item A major difficulty of proving Theorem~\ref{Theorem: Type II error under discrete setting} is in controlling randomness arising from the permutation procedure. We tackle this difficulty by building on the recent work of \cite{kim2020minimax}. In particular, we derive an upper bound for the $1-\alpha$ quantile of the permutation distribution of the test statistic, which holds with high probability. More details can be found in Appendix~\ref{Section: Proof of Theorem: Type II error under discrete setting}.
	\end{itemize}
\end{remark}

Next we switch gear to the continuous case of $Z$ and develop similar results as in the discrete case.

\subsection{Continuous $Z$} \label{Section: Continuous Z}
In this subsection, we build on the recent work of \cite{neykov2020minimax} and investigate the power of local permutation tests for continuous data. The idea of \cite{neykov2020minimax} is to carefully discretize $Z$ into several bins and apply the tests based on $T_{\mathrm{CI}}$ and $T_{\mathrm{CI},W}$ as if the original data were discrete. \cite{neykov2020minimax} investigate the type I and II errors of these tests and prove that they are minimax optimal under certain smoothness conditions. However, their tests depend on unspecified constants and, in their simulations, the authors instead use the local permutation procedure to determine critical values. Hence there is a gap between theory and practice. The goal of this subsection is to close this gap by showing that the local permutation tests have the same power property as the tests considered in \cite{neykov2020minimax}.

\paragraph{Discrete $X,Y$ and continuous $Z$.} 
First recall the setting described in Section~\ref{Section: Tightness} where $X$ and $Y$ are discrete random variables supported on $[\ell_1] \times [\ell_2]$ and $Z$ is continuous and bounded between $[0,1]$. Denote by $\mathcal{P}_{[0,1]}$ the collection of distributions $P_{X,Y,Z}$ of such random variables. Let $\mathcal{P}_{0,[0,1]} \subset \mathcal{P}_{[0,1]}$ where $X \independent Y | Z$ and $\mathcal{P}_{1,[0,1]} := \mathcal{P}_{[0,1]} \setminus \mathcal{P}_{0,[0,1]}$. Furthermore, let $\mathcal{P}_{1,[0,1],\mathrm{TV}}(L) \subset \mathcal{P}_{1,[0,1]}$ be the subset of alternative distributions that satisfies the TV smoothness:
\begin{align} \label{Eq: TV smoothness for Q}
	\mathcal{D}_{\mathrm{TV}}(P_{X,Y|Z=z}, P_{X,Y|Z=z'}) \leq L |z - z'| \quad \text{for all $z,z' \in [0,1]$.}
\end{align}
To describe the result of \cite{neykov2020minimax}, draw $N \sim \text{Pois}(n/2)$. If $N \leq n$, we take a random subset of size $N$ from $\{(X_i,Y_i,Z_i)\}_{i=1}^n$ and otherwise accept the null. Given that $N \leq n$, compute the unweighted test statistic $T_{\mathrm{CI}}$~(\ref{Eq: unweighted statistic}) using the binned data set $\bW_1,\ldots,\bW_M$ where $M = \ceil{n^{2/5}}$ and $\sum_{i=1}^M \sigma_i = N$. For a sufficiently large constant $\zeta$ depending on $L$, the corresponding test of \cite{neykov2020minimax} is defined as
\begin{align*}
	\phi_{\mathrm{NBW},1} := \mathds{1}(T_{\mathrm{CI}} \geq \zeta n^{1/5})\times \mathds{1}(N \leq n).
\end{align*} 
In terms of the type II error, the authors show that there exists a sufficiently large constant $c$ depending on $(\zeta, L, \ell_1, \ell_2)$, and for $\varepsilon \geq c n^{-2/5}$,
\begin{align} \label{Eq: type II error 1}
	\sup_{P_{X,Y,Z} \in \mathcal{P}_{1,[0,1],\mathrm{TV}}(L): \inf_{Q \in \mathcal{P}_{0,[0,1]}} \mathcal{D}_{\mathrm{TV}}(P_{X,Y,Z},Q) \geq \varepsilon} \mE_{P_{X,Y,Z}^N,N} [1 - \phi_{\mathrm{NBW},1}] \leq \frac{1}{100} + e^{-n/8}.
\end{align}
The type I error of $\phi_{\mathrm{NBW},1}$ is also guaranteed over a class of null distributions determined by certain smoothness conditions. 

As mentioned before, a test based on the unweighted U-statistic may not perform well when $\ell_1$ and $\ell_2$ potentially increase with $n$. To avoid this issue, \cite{neykov2020minimax} follow the idea of \cite{canonne2018testing} and propose another test based on the weighted test statistic $T_{\mathrm{CI},W}$~(\ref{Eq: weighted statistic}) where $M = \ceil{\frac{n^{2/5}}{(\ell_1\ell_2)^{1/5}}}$. More formally, for a sufficiently large $\zeta$ depending on $L$, the second test is defined as
\begin{align} \label{Eq: NBW 2}
	\phi_{\mathrm{NBW},2} := \mathds{1}(T_{\mathrm{CI},W} \geq \sqrt{\zeta M}) \times \mathds{1}(N \leq n).
\end{align}
Again, let $c$ be a sufficiently large constant depending on $(\zeta, L)$ such that $\varepsilon \geq c \frac{(\ell_1 \ell_2)^{1/5}}{n^{2/5}}$. Under this condition for $\varepsilon$ and an extra condition that $M \ell_1 \lesssim n$ for $\ell_1 \geq \ell_2$, Theorem 5.5 of \cite{neykov2020minimax} guarantees that
\begin{align} \label{Eq: type II error 2}
	\sup_{P_{X,Y,Z} \in \mathcal{P}_{1,[0,1],\mathrm{TV}}(L): \inf_{Q \in \mathcal{P}_{0,[0,1]}} \mathcal{D}_{\mathrm{TV}}(P_{X,Y,Z},Q) \geq \varepsilon} \mE_{P_{X,Y,Z}^N,N} [1 - \phi_{\mathrm{NBW},2}] \leq \frac{1}{100} + e^{-n/8}.
\end{align}
Furthermore, the authors prove that no CI test can be uniformly powerful in the TV distance when $\varepsilon$ is much less than $\frac{(\ell_1 \ell_2)^{1/5}}{n^{2/5}}$. That means, $\phi_{\mathrm{NBW},1}$ and $\phi_{\mathrm{NBW},2}$ achieve minimax optimal rate, while the optimality of $\phi_{\mathrm{NBW},1}$ is only guaranteed when $\ell_1$ and $\ell_2$ are bounded. However the optimal power of $\phi_{\mathrm{NBW},2}$ over a broader regime comes at the cost of decreasing the size of null distributions. In fact, the type I error of $\phi_{\mathrm{NBW},2}$ is guaranteed over a $\chi^2$-smooth class of null distributions, which is smaller than the class of null distributions considered for $\phi_{\mathrm{NBW},1}$. See \cite{neykov2020minimax} for more details. 

Having described the results of \cite{neykov2020minimax}, our aim is to reproduce the type II error guarantees (\ref{Eq: type II error 1}) and (\ref{Eq: type II error 2}) based on the same test statistics but with (explicit) cutoff values determined by the local permutation procedure. Given the data set of size $N \leq n$ where $N \sim \text{Pois}(n/2)$, let us denote the local permutation tests based on $T_{\mathrm{CI}}$ and $T_{\mathrm{CI},W}$ by $\phi_{\mathrm{perm},1}$ and $\phi_{\mathrm{perm},2}$, respectively. As described in Remark~\ref{Remark: full- vs half-permutation}, it greatly simplifies the power analysis for $T_{\mathrm{CI},W}$ when we consider the half-permutation method. Hence, unlike $\phi_{\mathrm{perm},1}$ that builds on Algorithm~\ref{alg:local permutation procedure}, the $p$-value of $\phi_{\mathrm{perm},2}$ is determined by the half-permutation method. Now in order to take into account the random sample size $N$, the final test functions are defined as $\phi_{\mathrm{perm},1}^\dagger:= \phi_{\mathrm{perm},1} \times \mathds{1}(N \leq n)$ and $\phi_{\mathrm{perm},2}^\dagger:= \phi_{\mathrm{perm},2} \times \mathds{1}(N \leq n)$. These tests have the following guarantee on the type II error rate. 

\begin{theorem}[Type II error for discrete $X,Y$ and continuous $Z$] \label{Theorem: Type II error under continuous setting I}
	Consider the local permutation tests $\phi_{\text{\emph{perm}},1}^\dagger$ and $\phi_{\text{\emph{perm}},2}^\dagger$, described above. In the setting of discrete ($X,Y$) and continuous $Z$, $\phi_{\text{\emph{perm}},1}^\dagger$ and $\phi_{\text{\emph{perm}},2}^\dagger$ have the same type II error guarantee as in (\ref{Eq: type II error 1}) and (\ref{Eq: type II error 2}), respectively. 
\end{theorem}

Several remarks are provided below. 

\begin{remark} \leavevmode \normalfont \label{Remark: Theorem: Type II error under continuous setting I}
	\begin{itemize}
		\item We highlight once again that the local permutation tests $\phi_{\mathrm{perm},1}^\dagger$ and $\phi_{\mathrm{perm},2}^\dagger$ do not require the knowledge on unspecified constant $\zeta$ in $\phi_{\mathrm{NBW},1}$ and $\phi_{\mathrm{NBW},2}$, while achieving the same power (in terms of rate) as stated in Theorem~\ref{Theorem: Type II error under continuous setting I}.
		\item However, such a nice property does not come for free. In general, the local permutation test requires a stronger condition on the class of null distributions than the corresponding theoretical test for type I error control. For instance, under the class of null distributions~$\mathcal{P}_{0, \mathrm{H},\gamma,\delta}(L)$ in Definition~\ref{Definition: Hellinger Lipschitzness}, we require that $n^{1/2}h^{\gamma} \rightarrow 0$ for $\gamma \in [1,2]$ and $n^{1/2}h^2 \rightarrow 0$ for $\gamma \geq 2$ for the local permutation test to be valid~(Theorem~\ref{Theorem: Validity of the permutation test under Hellinger smooth}). With the choice of $M = \ceil{n^{2/5}}$ (recall that $M = h^{-1}$), for instance, 
		we have only shown that $\phi_{\mathrm{perm},1}^\dagger$ is valid when $\gamma > 5/4$. On the other hand, $\phi_{\mathrm{NBW},1}$ can control the type I error even when $\gamma = 1$ \citep[see Theorem 5.2 of][]{neykov2020minimax}. Similarly, in order for $\phi_{\mathrm{perm},2}^\dagger$ to be valid under the same $\chi^2$-smoothness condition in Theorem 5.5 of \cite{neykov2020minimax}, we require that $n^{1/2}h^2 = \frac{(\ell_1\ell_2)^{2/5}}{n^{3/10}} \rightarrow 0$, which was not needed for $\phi_{\mathrm{NBW},2}$. In Section~\ref{Section: Double-binning strategy}, we attempt to partly address this drawback by introducing a novel double-binning strategy, which requires less stringent conditions for type I error control.
		\item The proof of Theorem~\ref{Theorem: Type II error under continuous setting I} is similar to that of Theorem~\ref{Theorem: Type II error under discrete setting}, that is, we show that the $1-\alpha$ quantile of the permutation distribution of $T_{\mathrm{CI}}$ is upper bounded by $\zeta n^{1/5}$ with high probability under the alternative. This result yields that the type II error of $\phi_{\mathrm{perm},1}^\dagger$ is upper bounded by that of $\phi_{\mathrm{NBW},1}$, up to a small error term. From here, we can directly benefit the previous bounds~(\ref{Eq: type II error 1}) and (\ref{Eq: type II error 2}) and show that the type II error of $\phi_{\mathrm{perm},1}^\dagger$ is small. The proof for $\phi_{\mathrm{perm},2}^\dagger$ follows similarly. The details can be found in Appendix~\ref{Section: Proof of Theorem: Type II error under continuous setting I}.
	\end{itemize}
\end{remark}

\paragraph{Continuous $X,Y,Z$.} Next we develop a similar result for the case where $(X,Y,Z)$ is supported on $[0,1]^3$ with a joint distribution absolutely continuous with respect to the Lebesgue measure. Let $\mathcal{P}_{[0,1]^3}$ be the set of distributions of such random variables and $\mathcal{P}_{0,[0,1]^3} \subset \mathcal{P}_{[0,1]^3}$ for which $X \independent Y |Z$. Let $\mathcal{P}_{1,[0,1]^3,\mathrm{TV}}$ be the subset of $\mathcal{P}_{1,[0,1]^3}:= \mathcal{P}_{[0,1]^3} \setminus \mathcal{P}_{0,[0,1]^3}$, which satisfies the TV smoothness condition in~(\ref{Eq: TV smoothness for Q}). In addition, we assume that, for any $P_{X,Y,Z} \in \mathcal{P}_{1,[0,1]^3,\mathrm{TV}}$, the corresponding conditional density function $p_{X,Y|Z}$ given any $z \in [0,1]$ belongs to $ \mathcal{H}^{2,s}(L)$, where $ \mathcal{H}^{2,s}(L)$ is the class of H\"{o}lder smooth functions $[0,1]^2 \mapsto \mathbb{R}$ with parameter $s$ defined in Definition~\ref{Definition: Holder smoothness} of the Appendix. 

In order to apply $T_{\mathrm{CI},W}$ to continuous data, we need to further discretize $X,Y$ into several bins. For this purpose, for a given $M > 0$ (specified in the sequel), consider a partition of $[0,1]$ into $M' = \ceil{M^{1/s}}$ bins of equal size denoted by $\{B_m'\}_{m=1}^{M'}$. We then transform $X$ (and similarly $Y$) through the map $g: [0,1] \mapsto \{1,\ldots,M'\}$ by defining $g(x) = m$ if and only if $x \in B_m$. On the other hand, we partition the conditional variable $Z$ into $M$ bins of equal size. Given this binned data set, we implement $\phi_{\mathrm{NBW},2}$ in~(\ref{Eq: NBW 2}), but with a different choice of $M = \ceil{n^{2s/(5s+2)}}$. According to Theorem 5.6 of \cite{neykov2020minimax}, the resulting test achieves minimax optimal power. In particular, they show that there exists a sufficiently large constant $c$ depending on $\zeta, L$, and if $\varepsilon \geq c n^{-2s/(5s+2)}$, then
\begin{align} \label{Eq: type II error 3}
	\sup_{P_{X,Y,Z} \in \mathcal{P}_{1,[0,1]^3,\mathrm{TV}}(L): \inf_{Q \in \mathcal{P}_{0,[0,1]^3}} \mathcal{D}_{\mathrm{TV}}(P_{X,Y,Z},Q) \geq \varepsilon} \mE_{P_{X,Y,Z}^N,N} [1 - \phi_{\mathrm{NBW},2}] \leq \frac{1}{100} + e^{-n/8}.
\end{align}
We now show that the corresponding local permutation test $\phi_{\mathrm{perm},2}^\dagger$ has the same power property. 

\begin{theorem}[Type II error for continuous $X,Y,Z$] \label{Theorem: Type II error under continuous setting II}
	Consider the local permutation test $\phi_{\text{\emph{perm}},2}^\dagger$ applied to the discretized data set described above. In the setting of continuous ($X,Y,Z$), $\phi_{\text{\emph{perm}},2}^\dagger$ has the same type II error guarantee as in (\ref{Eq: type II error 3}) and thereby shares the same optimal power as $\phi_{\text{\emph{NBW}},2}$.
\end{theorem}

The same points in Remark~\ref{Remark: Theorem: Type II error under continuous setting I} apply to Theorem~\ref{Theorem: Type II error under continuous setting II}. While $\phi_{\mathrm{perm},2}^\dagger$ has the same optimality as $\phi_{\mathrm{NBW},2}$ in terms of power, we need to restrict the class of null distributions further to rigorously control the type I error. In particular, under $\gamma$-Hellinger Lipschitzness with $\gamma \geq 2$, the underlying conditional density function should be smooth enough to 
ensure that $n^{1/2}h^2 = n^{-4s/(5s+2) + 1/2} \rightarrow 0$, equivalently $s > 2/3$, for type I error control. 
As we discussed in more detail in the introduction, the tension in CI testing between tightly controlling the type I error (requiring narrow bins) and ensuring high power (requiring, in some cases, wider bins) is a unique feature
of CI testing. 
In the next section, we introduce a novel double-binning permutation test that allows us to consider less smooth null distributions, while maintaining the power (up to a constant factor).

\section{Double-binning strategy} \label{Section: Double-binning strategy}
As we discussed above, the type I error of the local permutation test is guaranteed to be small over a smaller set of null distributions than the conservatively calibrated U-statistic test used by \cite{neykov2020minimax}.
 The reason for this gap is that the permutation approach relies on an additional condition for its validity. In particular, it requires that the binned distribution~$Q_{X,Y,\widetilde{Z}}^n$ and its CI projection $\widetilde{Q}_{X,Y,\widetilde{Z}}^n$ be close in the TV distance. The goal of this section is to mitigate this issue via double-binning. The idea is to consider bins of two distinct resolutions where a test statistic is computed over coarser bins, whereas the permutation procedure is implemented over finer bins. This double-binning strategy allows us to keep the TV distance smaller than the single-binning approach, while maintaining similar power under certain conditions. 

To elaborate on the idea, recall that $\{B_1,\ldots,B_M\}$ is a partition of $\mathcal{Z}$. For $m=1,\ldots,M$, let us further partition $B_m$ into $b$ bins, which results in smaller bins $\{B_{m,1},\ldots,B_{m,b}\}$. For $(m,k) \in[M] \times [b]$, let $\sigma_{m,k}$ be the sample size within $B_{m,k}$, and $\bW_{m,k}$ denote the set of the pairs of $(X_i,Y_i)$ that belong to $B_{m,k}$. More formally, by letting $(X_{i,m,k},Y_{i,m,k})$ be the $i$th pair in $B_{m,k}$, we write $\bW_{m,k} := \{ (X_{1,m,k},Y_{1,m,k}),\ldots,(X_{\sigma_{m,k},m,k},Y_{\sigma_{m,k},m,k}) \}$ when $\sigma_{m,k} \geq 1$ and otherwise $\boldsymbol{W}_{m,k} = \emptyset$. Under this setting, we compute the test statistic $T_{\mathrm{CI}}$ as in (\ref{Eq: test statistic for CI testing}) based on the larger bins $\{B_1,\ldots,B_M\}$. The permuted test statistic $T_{\mathrm{CI}}^{\boldsymbol{\pi}}$ is computed similarly as before except now that $Y$ values are permuted within the smaller bins. We then compute the permutation $p$-value as in (\ref{Eq: permutation p-value}) by counting how many permuted statistics are larger than or equal to $T_{\mathrm{CI}}$. 

\begin{figure}[t!]
	\begin{center}		
		\includegraphics[width=40em]{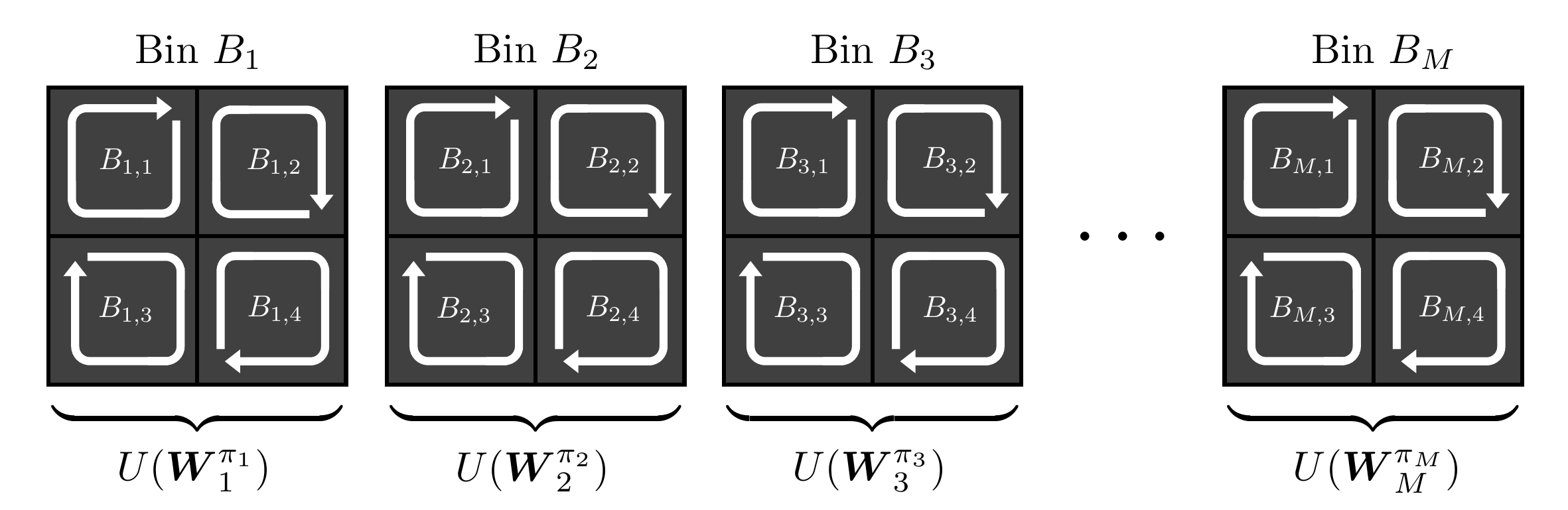} 
		\caption{\small Illustration of the double-binning strategy where we consider bins of two distinct resolutions. The purpose of coarser bins $\{B_1,\ldots,B_M\}$ is to compute a test statistic, whereas the purpose of finer bins $\{B_{1,1},\ldots,B_{M,4}\}$ is to perform the permutation procedure. In particular, we only permute the observations within the finer bins. This strategy allows us to keep the type I error under control over a larger class of null distributions, while maintaining the competitive power under certain conditions.} \label{Figure: doublebinning}
	\end{center}
\end{figure}

Here, to simplify our theoretical analysis, we focus our attention on a subset of all possible local permutations. In particular, within each small bin $B_{m,k}$ for $(m,k) \in [M] \times [b]$, we only consider a set of cyclic permutations of $\{1,\ldots, \sigma_{m,k}\}$. As an illustration, when $ \sigma_{m,k} = 4$, we have four distinct cyclic permutations of $\{1,2,3,4\}$ as $\{1,2,3,4\}, \{2,3,4,1\}, \{3,4,1,2\}, \{4,1,2,3\}$. A formal definition of a cyclic permutation can be found in Definition~\ref{Definition: cyclic permutation}. This cyclic restriction results in $K_\ast := \prod_{(m,k) \in [M] \times [b]} \max\{\sigma_{m,k},1\}$ number of local permutations, and the set of all possible such cyclic local permutations is denoted by $\boldsymbol{\Pi}_{\mathrm{cyclic}}$. Given this notation, we describe the local permutation procedure via double-binning in Algorithm~\ref{alg:double-binning permutation procedure}. We also refer to Figure~\ref{Figure: doublebinning} for an illustration of the procedure.

\begin{algorithm}[h]
	\caption{Local permutation procedure via double-binning} \label{alg:double-binning permutation procedure}
	\begin{algorithmic}
		\ENSURE  data $\{(X_i,Y_i,Z_i)\}_{i=1}^n$, a super-partition of $\mathcal{Z}$: $\big\{B_{1},\ldots,B_{M} \big\}$, sub-partitions of each $B_m: \{B_{m,1},\ldots,B_{m,b}\}$ for $m=1,\ldots,M$, a test statistic $T_{\mathrm{CI}}$, a nominal level $\alpha$
		\begin{enumerate}
			\item For each $\boldsymbol{\pi} \in \boldsymbol{\Pi}_{\mathrm{cyclic}}$, compute $T_{\mathrm{CI}}^{\boldsymbol{\pi}}$ as in (\ref{Eq: permuted statistic}) and denote the resulting statistics by $T_{\mathrm{CI}}^{\boldsymbol{\pi}_1},\ldots,T_{\mathrm{CI}}^{\boldsymbol{\pi}_{K_\ast}}$. 
			\item By comparing the statistic $T_{\mathrm{CI}}$ in~(\ref{Eq: test statistic for CI testing}) with the permuted ones, calculate the $p$-value as
			\begin{align*}
				p_{\mathrm{perm}} = \frac{1}{K_\ast} \sum_{\boldsymbol{\pi}_i \in \boldsymbol{\Pi}_{\mathrm{cyclic}}} \mathds{1} \big\{ T_{\mathrm{CI}}^{\boldsymbol{\pi}_i} \geq T_{\mathrm{CI}} \big\}.
			\end{align*}
			\item Given the nominal level $\alpha \in (0,1)$, define the test function $\phi_{\mathrm{perm},n} = \mathds{1}(p_{\mathrm{perm}} \leq \alpha)$ and reject the null when $\phi_{\mathrm{perm},n}=1$.
		\end{enumerate}
	\end{algorithmic}
\end{algorithm}

Next, we discuss the type I and type II errors of the local permutation test via double-binning. The analysis of the type I error is relatively straightforward. Indeed, all of the results in Section~\ref{Section: Validity} continue to hold for this new approach but now the maximum diameter~(\ref{Eq: maximum diameter}) becomes smaller as it is defined over the finer bins. Therefore, in view of Theorem~\ref{Theorem: Validity of the permutation test under Hellinger smooth} and Theorem~\ref{Theorem: Validity of the permutation test under Renyi smooth}, an upper bound for the type I error can be much tighter than the single-binning approach based on $\{B_1,\ldots,B_M\}$. The only caveat, here, is that we consider the set of cyclic local permutations $\boldsymbol{\Pi}_{\mathrm{cyclic}}$, rather than all possible local permutations. However, this change does not affect the validity.

\begin{proposition}[Type I error of the double-binning test] \label{Proposition: Type I error of the double-binned test}
	Let $\phi_{\mathrm{perm},n}$ be the local permutation test via double-binning in Algorithm~\ref{alg:double-binning permutation procedure}. Then the same bounds~(\ref{Eq: upper bounds under Hellinger Lipschitzness}) and (\ref{Eq: upper bounds under Renyi Lipschitzness}) hold for $\phi_{\mathrm{perm},n}$ with the maximum diameter $h$ defined as 
	\begin{align*} 
		h = \max_{(m,k) \in [M] \times [b]} \sup_{z,z' \in B_{m,k}} \delta(z,z').
	\end{align*}
\end{proposition}

We note that the above result holds under Poissonization in a similar fashion to Proposition~\ref{Proposition: Poissonization} in the Appendix. Indeed, once the bounds~(\ref{Eq: upper bounds under Hellinger Lipschitzness}) and (\ref{Eq: upper bounds under Renyi Lipschitzness}) are given, the validity of the double-binning test under Poissonization can be proved along the same lines of the proof of Proposition~\ref{Proposition: Poissonization}. Given the above proposition, our main concern is the type II error. Here, to illustrate ideas, we only focus on the case where $(X,Y,Z) \in [\ell_1] \times [\ell_2] \times [0,1]$ and $\ell_1,\ell_2$ are fixed. In addition, we recall the definition of $\mathcal{P}_{1,[0,1],\mathrm{TV}}(L)$ from Section~\ref{Section: Continuous Z} and consider a subset of $\mathcal{P}_{1,[0,1],\mathrm{TV}}(L)$, denoted by $\mathcal{P}'_{1,[0,1],\mathrm{TV}}(L)$, such that the marginal density of $Z$ is bounded below by some fixed constant $c_{\mathrm{low}}>0$, i.e.~$p_Z(z) \geq c_{\mathrm{low}}$ for all $z \in [0,1]$. Furthermore, we assume that any distribution in $\mathcal{P}'_{1,[0,1],\mathrm{TV}}(L)$ satisfies
\begin{align*}
	\mathcal{D}_{\mathrm{TV}}(P_{X|Z=z}, P_{X|Z=z'}) \leq L |z - z'| \quad \text{and} \quad \mathcal{D}_{\mathrm{TV}}(P_{Y|Z=z}, P_{Y|Z=z'}) \leq L |z - z'|, 
\end{align*}
for all $z,z' \in [0,1]$. 

To describe type II error results, we consider the test statistic $T_{\mathrm{CI}}$ used in $\phi_{\mathrm{NBW},1}$. Moreover, assume that each of the coarser bins and the finer bins has the same length of an interval $1/M$ and $1/(Mb)$, respectively. Now, given a sample $\{(X_i,Y_i,Z_i)\}_{i=1}^N$ where $N \sim \mathrm{Pois}(n/2)$, we compute $\phi_{\mathrm{perm},N}$ based on $T_\mathrm{CI}$ using Algorithm~\ref{alg:double-binning permutation procedure} at significance level $\alpha = 0.01$ (for simplicity) and define $\phi_{\mathrm{perm},n}^\dagger := \phi_{\mathrm{perm},N} \times \mathds{1}\{N \leq n\}$ as in (\ref{Eq: modified local permutation test}). Under this setting, the resulting test has the following type II error guarantee.

\begin{theorem}[Type II error of the double-binning test] \label{Theorem: Type II error of the double-binned test}
	Consider the test $\phi_{\mathrm{perm},n}^\dagger$ defined above with the number of larger bins $M = \ceil{n^{2/5}}$. Suppose that we choose $b$ such that $n /(Mb)^2 > 400 c_{\mathrm{low}}^{-1}$ and $Mb \rightarrow \infty$ as $n \rightarrow \infty$. Suppose further that $\varepsilon \geq cn^{-2/5}$ for a sufficiently large $c$ depending on ($L, \ell_1, \ell_2$). Then 
	\begin{align*}
		\sup_{P_{X,Y,Z} \in \mathcal{P}'_{1,[0,1],\mathrm{TV}}(L): \inf_{Q \in \mathcal{P}_{0,[0,1]}} \mathcal{D}_{\mathrm{TV}}(P_{X,Y,Z},Q) \geq \varepsilon} \mE_{P_{X,Y,Z}^N,N} [1 - \phi_{\mathrm{perm},n}^\dagger] \leq \frac{1}{100} + e^{-n/8} + \rho_n,
	\end{align*}
	where $\rho_n$ is a positive sequence converging to zero as $n \rightarrow \infty$.
\end{theorem}

We first note that an explicit form of $\rho_n$ in the above result can be found in (\ref{Eq: delta}) given in the Appendix. Next, observe that the above double-binning test achieves the same minimax separation rate as $\phi_{\mathrm{NBW},1}$ without being dependent on an unspecified constant in its critical value. Its guarantee holds over a smaller set of alternative distributions, namely $ \mathcal{P}'_{1,[0,1],\mathrm{TV}}(L)$ (in contrast, the result of 
\cite{neykov2020minimax} does not require any lower bound on the density of $Z$). 

On the other hand, compared to the corresponding permutation test via single-binning, the double-binning method controls the type I error rate over a larger class of null distributions without sacrificing power up to a constant factor. In particular, with an optimal choice of $M = \ceil{n^{2/5}}$, the single-binning test requires $\gamma > 5/4$ in order to control the type I error under $\gamma$-Hellinger smoothness. In contrast, the double-binning test is valid for $\gamma > 1$ as long as $b$ is chosen appropriately such that $n/(Mb)^{2\gamma} \rightarrow 0$ but $n/(Mb)^2$ is bounded below by a large constant.

\begin{remark}[Extension to other settings] \normalfont
	It is worth highlighting that the validity result of the double-binning test in Proposition~\ref{Proposition: Type I error of the double-binned test} holds universally for any binning-based test statistic. Theorem~\ref{Theorem: Type II error of the double-binned test} regarding the power, on the other hand, is based specifically on the test statistic $T_{\mathrm{CI}}$ used in $\phi_{\mathrm{NBW},1}$. While we believe that a similar power result can be developed based on other test statistics including $T_{\mathrm{CI},W}$ used in $\phi_{\mathrm{NBW},2}$, a detailed treatment of this direction is beyond the scope of this paper. Another important direction one can pursue is to see whether the lower bound restriction on the density of $Z$ can be removed in Theorem~\ref{Theorem: Type II error of the double-binned test}. We leave these topics to future work.
\end{remark}

\section{Simulations} \label{Section: Simulations}

In this section, we illustrate the numerical performance of the local permutation test through Monte Carlo simulations. Throughout our experiments, we use Monte Carlo simulations to compute the permutation $p$-value, defined as $p_{\mathrm{perm}}$ in (\ref{Eq: permutation p-value}). In particular, we draw	 $B$ permutations, denoted by $\boldsymbol{\pi}_1',\ldots,\boldsymbol{\pi}_B'$, from $\boldsymbol{\Pi}$ with replacement. Then the permutation $p$-value using the test statistic $T_{\mathrm{CI}}$ is defined as
\begin{align} \label{Eq: approximated p-value}
	\widehat{p}_{\mathrm{perm}} := \frac{1}{B+1} \Bigg[ \sum_{i=1}^B \mathds{1}\{ T_{\mathrm{CI}}^{\boldsymbol{\pi}_i'} \geq T_{\mathrm{CI}}  \} + 1 \Bigg].
\end{align}
It is well-known that $\mathds{1}(\widehat{p}_{\mathrm{perm}} \leq \alpha)$ is a valid test in finite-sample scenarios, whenever $T_{\mathrm{CI}},T_{\mathrm{CI}}^{\boldsymbol{\pi}_1},\ldots,T_{\mathrm{CI}}^{\boldsymbol{\pi}_K}$ are exchangeable under the null. Furthermore, $\widehat{p}_{\mathrm{perm}}$ can be arbitrarily close to $p_{\mathrm{perm}}$ for a sufficiently large $B$, and we take $B=100$ in our simulations. We also note that the type I error and the power presented in this section are approximated by repeating simulations $1000$ times at significance level $\alpha = 0.05$. 

\vskip 1em

\subsection{Experiment 1}  \label{Section: Experiment 1}
In our first experiment, we demonstrate Theorem~\ref{Theorem: Hardness of CI testing} in a discrete setting of $(X,Y,Z)$ where $(X,Y,Z)$ is distributed over $[\ell_1] \times [\ell_2] \times [M]$ and $\ell_1 = \ell_2 = 2$ and $M \in \{10,20,\ldots,110,120\}$. In particular, we let $Z$ have a multinomial distribution with equal probabilities over the bins. Similarly, we let $X$ have a multinomial distribution with equal probabilities, and set $X=Y$ and $(X,Y) \independent Z$. We are under the alternative hypothesis where $X$ and $Y$ are perfectly correlated conditional on $Z$. In this setting, we compute the empirical power of the local permutation test based on the test statistic $T_{\mathrm{CI}}$ considered in Theorem~\ref{Theorem: Type II error under discrete setting} by varying $(M,n)$. The result can be found in the left panel of Figure~\ref{Figure: Discrete Case}. As predicted by Theorem~\ref{Theorem: Hardness of CI testing}, we see that the power of the test degrades quickly as $M$ increases for any given $n$. This in turn illustrates that CI testing is impossible unless the probability observing the same value of $Z$ is properly controlled.  

\vskip 1em

\subsection{Experiment 2}  \label{Section: Experiment 2}
\begin{figure}[t!]
	\begin{center}		
		\begin{minipage}[b]{0.46\textwidth}
			\includegraphics[width=\textwidth]{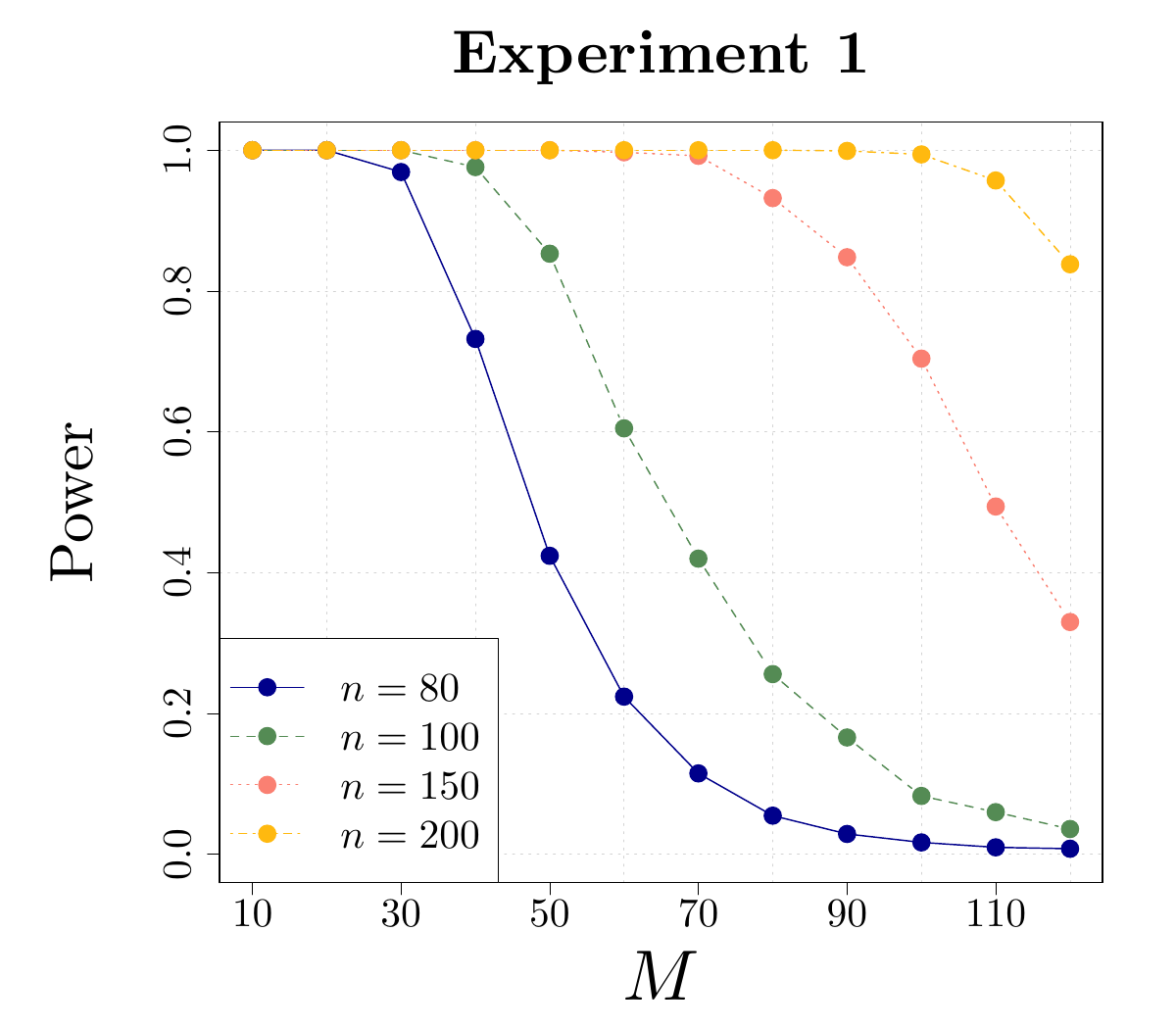}
		\end{minipage} 
		\hskip 2em
		\begin{minipage}[b]{0.46\textwidth}
			\includegraphics[width=\textwidth]{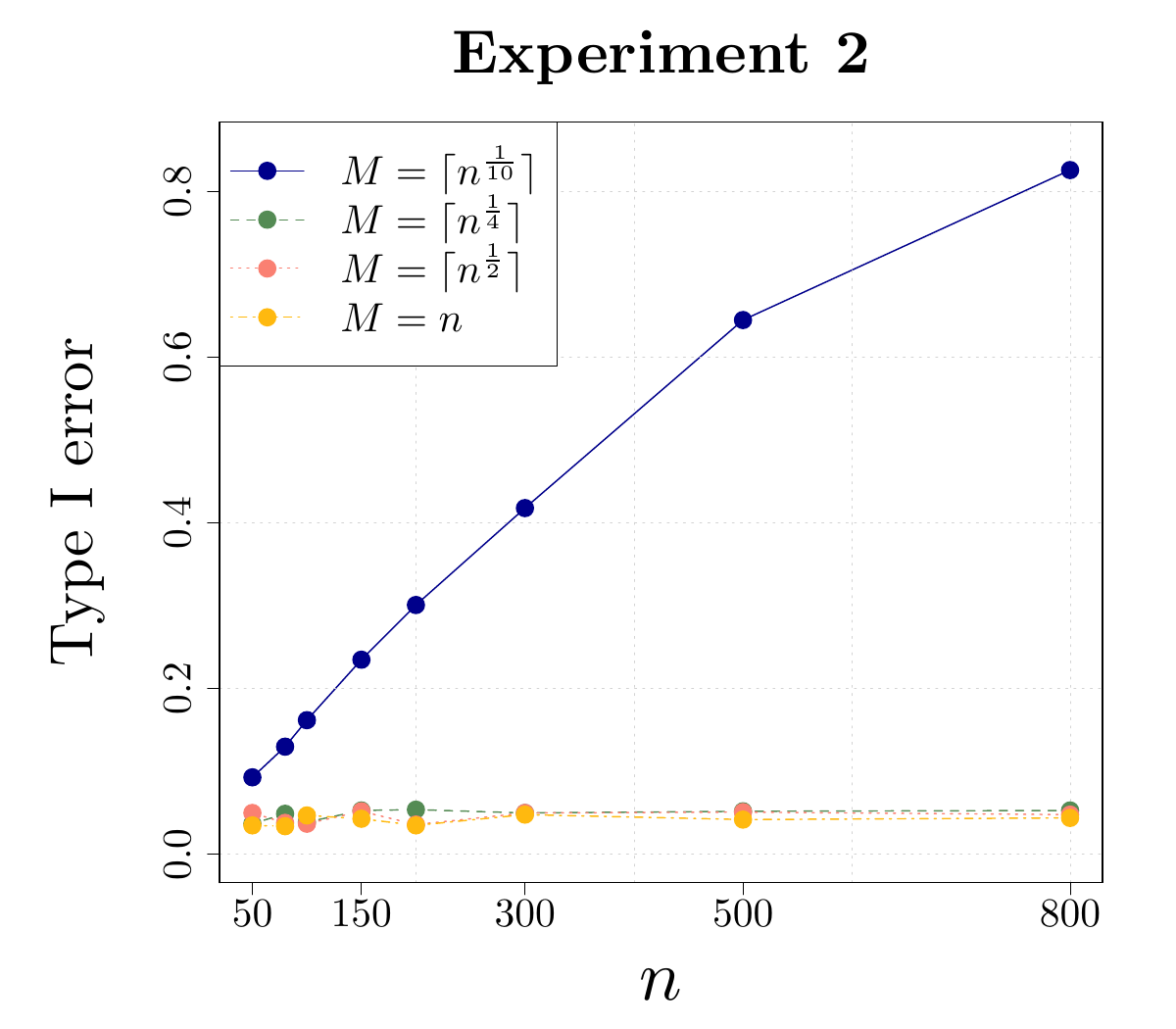}
		\end{minipage}
		\caption{\small Illustration of the power and type I error of the local permutation test using $T_{\mathrm{CI}}$. The left panel demonstrates that the power of the permutation test keeps decreasing as $M$ increases even though $X$ and $Y$ are perfectly correlated conditional on $Z$; thereby, confirming Theorem~\ref{Theorem: Hardness of CI testing}. The right panel demonstrates that the type I error of the local permutation test is well-controlled under the smoothness condition, described in Section~\ref{Section: Experiment 2}, unless $M$ is chosen such that $nM^{-4} \rightarrow \infty$.} \label{Figure: Discrete Case}
	\end{center}
\end{figure}

In our second experiment, we demonstrate the validity result of the local permutation test in Section~\ref{Section: Validity under smoothness conditions}. To generate the data, we consider a distributional setup used in construction of our lower bound result (Theorem~\ref{Theorem: Lower bounds}). In particular, we consider the marginal density of $Z$ in~(\ref{Eq: density of z}) with $\epsilon = n^{-1}$, and the conditional density of $X|Z=z$ in (\ref{Eq: Example of the distribution}). We further let $Y|Z=z$ have the same conditional density as $X|Z=z$ for all $z \in [0,1]$, while satisfying $X\independent Y|Z$. As proved in Appendix~\ref{Section: Hellinger case 1} and Appendix~\ref{Section: Renyi}, the considered distribution satisfies $\gamma$-Hellinger Lipschitzness with any fixed $\gamma \geq 1$ as well as $\gamma$-R{\'e}nyi Lipschitzness with any fixed $\gamma >0$. Therefore, by Theorem~\ref{Theorem: Validity of the permutation test under Hellinger smooth} and Theorem~\ref{Theorem: Validity of the permutation test under Renyi smooth}, the type I error of the local permutation test based on any test statistic is approximately $\alpha$ as long as $nM^{-4} \rightarrow 0$. We also note from Theorem~\ref{Theorem: Lower bounds} that there exists a test statistic such that the corresponding local permutation test fails to control the type I error rate when $nM^{-4} \rightarrow \infty$. To demonstrate both results, we use the test statistic in~(\ref{Eq: unweighted statistic}) by varying the number of bins $M \in \{n, \ceil{n^{1/2}}, \ceil{n^{1/4}}, \ceil{n^{1/10}}\}$. The result is given in the right panel of Figure~\ref{Figure: Discrete Case}. As we can see from the result, the type I error is well controlled when $M$ is chosen such that $nM^{-4} \rightarrow 0$. On the other hand, the error tends to increase when $nM^{-4} \rightarrow \infty$, which coincides with our theory.

\subsection{Experiment 3} \label{Section: Experiment 3}
In our third experiment, we illustrate type I and II error control of the double-binning test in Theorem~\ref{Theorem: Type II error of the double-binned test} by setting $M = b = \ceil{n^{2/5}}$. The permutation $p$-value of the double-binning procedure is approximated similarly as $\widehat{p}_{\mathrm{perm}}$ in (\ref{Eq: approximated p-value}) but by drawing $B$ cyclic permutations from $\boldsymbol{\Pi}_{\mathrm{cyclic}}$ without replacement. As before, we choose $B=100$ for our third simulation as well. To demonstrate the performance, we let $Z$ have a uniform distribution over the interval $[0,1]$ and $X,Y$ be Bernoulli random variables with the following conditional probability mass functions: 
\begin{align} \label{Eq: conditional pmf}
	p_{X|Z}(X=1|z) = p_{Y|Z}(Y=1|z) = e^{\sin(\theta z)} /4.
\end{align}
The considered distribution depends on the parameter $\theta$, which controls the smoothness of the conditional probability mass function. In particular, the conditional marginals~(\ref{Eq: conditional pmf}) become more wiggly as $\theta$ increases, which makes it more difficult to control the type I error under the null. 

\begin{figure}[t!]
	\begin{center}		
		\begin{minipage}[b]{0.46\textwidth}
			\includegraphics[width=\textwidth]{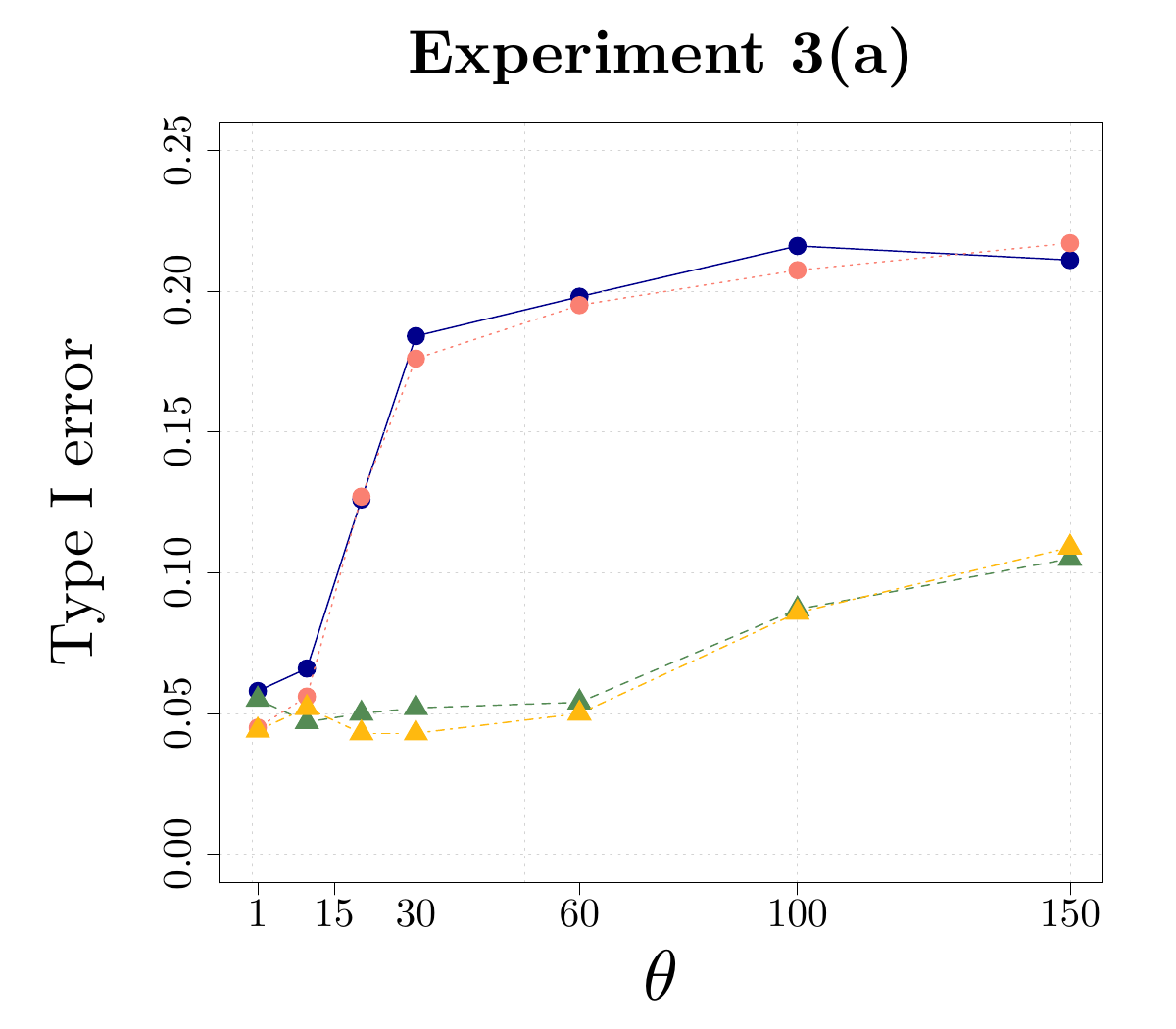}
		\end{minipage} 
		\hskip 2em
		\begin{minipage}[b]{0.46\textwidth}
			\includegraphics[width=\textwidth]{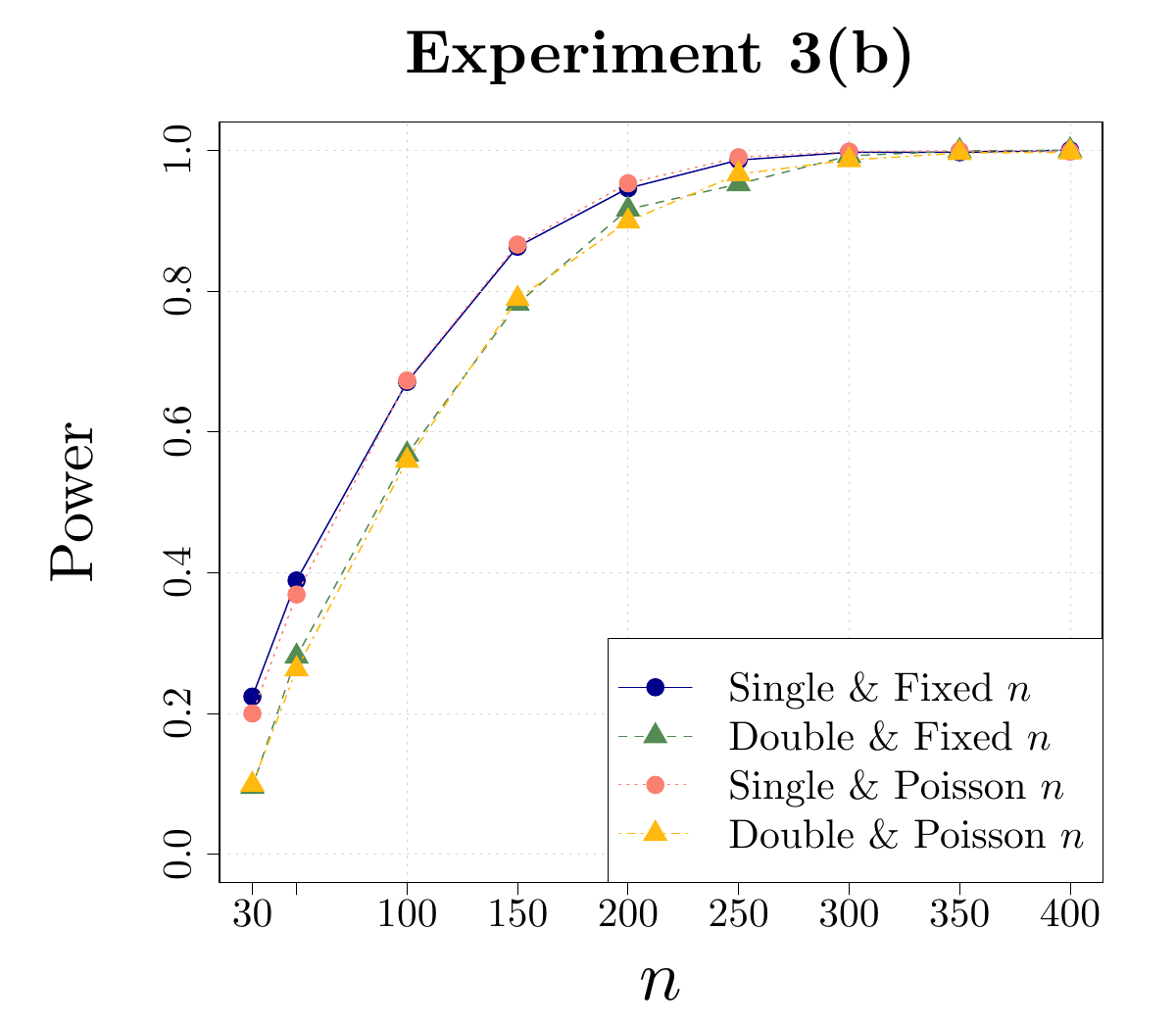}
		\end{minipage}
		\caption{\small Illustration of the type I error and power of the single- and double-binning tests using $T_{\mathrm{CI}}$. The left-panel considers the null distribution with conditional marginals~(\ref{Eq: conditional pmf}) and shows that the double-binning test has better type I error control than the single-binning test over different values of smoothness parameter $\theta$. The right-panel is concerned with the power as $n$ varies, demonstrating that the single-binning performs better than the double-binning test in terms of power. We also include simulation results under Poissonization, which illustrates that the Possionization version performs similarly as the fixed-$n$ counterpart. A more detailed explanation can be found in Section~\ref{Section: Experiment 3}.} \label{Figure: Double-binning}
	\end{center}
\end{figure}

\begin{enumerate}[(a).]
	\item \textbf{Type I error.} To illustrate type I error control, we consider the null distribution with the conditional marginals~(\ref{Eq: conditional pmf}). We then draw $n=100$ samples from the null distribution, and compute the test statistic $T_{\mathrm{CI}}$ as well as the $p$-value. The finite-sample type I error is approximated by Monte Carlo simulations and the result is provided in the left panel of Figure~\ref{Figure: Double-binning}. As a reference point, we also consider the single-binning test based on the same test statistic with $M = \ceil{n^{2/5}}$ and its type I error rate is also provided in the left panel of Figure~\ref{Figure: Double-binning}. From the result, we see that the type I error of the single-binning test increases with $\theta$ much faster than that of the double-binning test. This empirical result supports Proposition~\ref{Proposition: Type I error of the double-binned test} that claims that the double-binning test is valid over a larger class of null distributions than the corresponding single-binning test. 
	\item \textbf{Power.} To illustrate the power performance, we consider an alternative distribution with the same marginals~(\ref{Eq: conditional pmf}) with $\theta = 1$. In particular, by writing $f(z) = e^{\sin(z)} /4$, we set 
	\begin{align*}
		&p_{XY|Z}(X=1,Y=1|z) = f(z)^2 + f(z) / 5, \\[.5em]
		& p_{XY|Z}(X=0,Y=0|z) = \{1 - f(z)\}^2 + f(z) / 5  \quad \text{and} \\[.5em]
		& p_{XY|Z}(X=1,Y=0|z) = p_{XY|Z}(X=0,Y=1|z) = 4f(z)/5 - f(z)^2.
	\end{align*}
	One can check that the above joint distribution is a valid alternative distribution where the conditional joint distribution differs from the product of the conditional marginal distributions. With $n$ draws from the above distribution, we compute the same test statistic and $p$-value as before and approximate the power of the test by changing $n \in \{30,50,100,150,\ldots,350,400\}$. The right-panel of Figure~\ref{Figure: Double-binning} collects the power approximates for both single-binning and double-binning tests. Overall, the power of the single-binning test is higher than that of the double-binning test, but the difference seems marginal, especially when the power is close to one. This may be viewed as empirical evidence of Theorem~\ref{Theorem: Type II error of the double-binned test}, which shows that both tests have the same power up to a constant factor in certain regimes. 
\end{enumerate}
In both the panels of Figure~\ref{Figure: Double-binning}, we also present the type I error and power of the corresponding tests under Poisson sampling. Specifically, for $N \sim \text{Poisson}(n)$, we draw i.i.d.~samples $(X^N,Y^N,Z^N)$ from the joint distribution of $(X,Y,Z)$ and compute each permutation test using $(X^N,Y^N,Z^N)$. As we can see, the results under Poisson sampling are not significantly different from the previous results with the fixed sample size~$n$, 
and we anticipate that our theoretical results will continue to hold with a fixed sample-size.

\section{Discussion} \label{Section: Discussion}
In this paper, we investigated several statistical properties of the local permutation method for CI testing. We started by presenting a new hardness result of CI testing, which, along with the recent work of \cite{shah2020hardness}, motivates us to consider reasonable assumptions under which CI testing becomes possible. Under certain smoothness assumptions, we provided upper bounds for the type I error of the local permutation test and further showed that these bounds are tight in some cases. Turning to the power, we demonstrated that the local permutation test can retain minimax power, while rigorously controlling the type I error, under certain circumstances. In particular, we showed that the local permutation tests using the same test statistics in \cite{canonne2018testing,neykov2020minimax} have the same power guarantee. However, compared to the previous tests, the type I error of the local permutation test is guaranteed over a smaller set of null distributions in the continuous case of $Z$. To this end, we introduced and analyzed a double-binning strategy, which mitigates this drawback.

\vskip .5em

\noindent \textbf{Future directions.} We close by discussing several interesting directions for future work.
\begin{itemize}
	\item \textbf{Adaptive binning strategy.} Throughout this paper, we have been mainly concerned with equal-sized bins. This strategy, as we saw earlier, can lead to optimal CI tests from a minimax perspective. However, when there exists a local structure on the distribution of $Z$, many of bins would be empty. In this case, it would be more desirable to consider an adaptive binning scheme that uses different sizes of bins over different regions of $Z$. This strategy, potentially data-dependent, requires a more delicate analysis for both type I and type II error control, which we leave to future work. As mentioned before, it would also be interesting to see whether it is possible to develop an adaptive test to the unknown smoothness parameters without sacrificing power much. 
	\item \textbf{Different metrics.} In this work, we have focused on the smoothness conditions based on the generalized Hellinger distance and R{\'e}nyi divergence. It may be possible to obtain different and potentially sharper validity conditions by considering other metrics. In addition, one can impose a smoothness assumption on higher order derivatives of a conditional distribution and see whether its improves the validity result. It is also worth investigating the minimax power of the local permutation test in different metrics other than the TV distance.
	\item \textbf{Other test statistics.} While our validity result can be applied to any binning-based statistic, the power analysis was specifically based on U-statistics with discrete-type kernels. We believe that it is also possible to obtain similar minimax power results using other test statistics. In particular, exploring the power of RKHS-based test statistics~\citep[e.g.][]{fukumizu2008kernel} is a promising direction for future work. This can be potentially explored by building on the recent work of \cite{meynaoui2019adaptive,kim2020minimax} that investigate the minimax power of unconditional independence tests based on Gaussian kernels.	
	\item \textbf{Depoissonization.} As in the previous work~\citep{canonne2018testing,balakrishnan2019hypothesis,neykov2020minimax}, it was crucial to use Poissonization technique for our power analysis. Since a Poisson random variable is tightly concentrated around its mean, it sounds plausible that the same power guarantee can be achieved by the local permutation test without Poissonization as empirically demonstrated in Section~\ref{Section: Experiment 3}. However, a formal proof is not available at the current stage, which we leave as future work. 	
	\item \textbf{Multivariate $Z$.} Our results on type I error control show that the validity of a local permutation test crucially relies on the maximum diameter of bins, which is well-defined even when $Z$ is a multivariate random vector. Our results on minimax power, on the other hand, only deal with the univariate case of $Z$. Indeed, tight minimax separation rates for CI testing are only known for the case when the dimension of $Z$ is either one or two~\citep{neykov2020minimax}. As the purpose of our work is to demonstrate that local permutation tests can achieve the same optimal power as their theoretical counterparts, we focus only on the univariate $Z$ case. It therefore remains for future work to establish minimax separation rates for a general multivariate case and see whether local permutation tests can achieve these rates.
	\item \textbf{Theoretical vs.~practical calibration.} While the double-binning strategy allows us to consider less smooth null distributions than the corresponding single-binning test, there still remain settings where minimax-optimal tests are available \citep{neykov2020minimax} but where we do not have practical methods to calibrate these optimal tests. Understanding these gaps and more clearly establishing the fundamental limits for \emph{practical} conditional independence tests are an important direction for future work.
\end{itemize}

\section*{Acknowledgements}
This work was partially supported by funding from the NSF grants DMS-1713003, DMS-2113684 and CIF-1763734, as well as Amazon AI and a Google Research Scholar Award to SB. IK was partially supported by EPSRC grant EP/N031938/1 and Yonsei University Research Fund of 2021-22-0332.

\bibliographystyle{apalike}
\bibliography{reference}

\appendix

\allowdisplaybreaks

\section{Proofs of the results in the main text} \label{Section: proofs}

In this section, we collect the proofs of the results in the main text. Some auxiliary lemmas required for the proofs can be found in Appendix~\ref{Section: Auxiliary lemmas}. We start by introducing additional notation that simplifies our presentation.

\paragraph{Additional notation.} For sequences $a_n$ and $b_n$, we write $a_n \lesssim b_n$ or $b_n \gtrsim a_n$ if there exists an absolute constant $C>0$ such that $a_n \leq C b_n$ for all $n$. In addition, $a_n = o(b_n)$ means that $a_n/b_n \rightarrow 0$ as $n \rightarrow \infty$. As convention, $\|\cdot \|_p$ represents the $L_p$ norm. $C,C_1,C_2,\ldots$ denote some constants whose value may differ in different places.

\subsection{Proof of Theorem~\ref{Theorem: Hardness of CI testing}} \label{Section: Proof of Theorem: Hardness of CI testing}
As mentioned in the main text, the proof of this result is highly motivated by the impossibility result of distribution-free conditional predictive inference (see e.g.,~Lemma A.1 of \citealt{foygel2019limits} and Lemma 1 of \citealt{barber2020distribution}). We also mention again that negative results for testing argued via sampling with replacement can be traced back to Example 1 of \cite{gretton2012kernel}. Following the construction of Lemma 1 in \cite{barber2020distribution}, let $\mathcal{L} = \{(X_1,Y_1,Z_1),\ldots,(X_J,Y_J,Z_J)\} \overset{\mathrm{i.i.d.}}{\sim} P_{X,Y,Z} \in \mathcal{P}_1$ be ghost samples for some $J \geq n(n-1)$. Given $\mathcal{L}$, draw $\{(X_1',Y_1',Z_1'),\ldots,(X_n',Y_n',Z_n')\}$ from $\mathcal{L}$ \emph{without replacement}. We denote this conditional sampling distribution by $P_{\mathrm{without}}$. Then by marginalizing over $\mathcal{L}$, we simply have the identity
\begin{align} \label{Eq: marginalizing over L}
	\mE_{P_{X,Y,Z}^n}[\phi] = \mE_{\mathcal{L}} [\mE_{P_{\mathrm{without}}} [\phi|\mathcal{L}]].
\end{align} 
Next consider drawing $\{(X_1',Y_1',Z_1'),\ldots,(X_n',Y_n',Z_n')\}$ from $\mathcal{L}$ \emph{with replacement} and denote this conditional sampling distribution by $P_{\mathrm{with}}$. By \cite{stam1978distance}, the total variation distance between $P_{\mathrm{with}}$ and $P_{\mathrm{without}}$ can be bounded by $n(n-1)/J$. Therefore, using the identity~(\ref{Eq: marginalizing over L}), we have
\begin{align*}
	\mE_{P_{X,Y,Z}^n}[\phi]  ~\leq~  \mE_{\mathcal{L}} [\mE_{P_{\mathrm{with}}} [\phi|\mathcal{L}]] + \frac{n(n-1)}{J}.
\end{align*}
Let $\mathcal{A}$ be the event that $Z_1,\ldots,Z_J~\text{are distinct}$. Under the event $\mathcal{A}$, we have $X \independent Y |Z$ for $(X,Y,Z) \sim P_{\mathrm{with}}$ since the conditional distribution of $X,Y$ given any $Z$ takes a single value~(i.e.~conditionally degenerate). Therefore we should have 
\begin{align*}
	\mE_{P_{\mathrm{with}}} [\phi|\mathcal{L},\mathcal{A}] \leq \alpha \quad \text{and} \quad \mE_{P_{\mathrm{with}}} [\phi|\mathcal{L},\mathcal{A}^c] \leq 1
\end{align*}
as $\sup_{P \in \mathcal{P}_{0,\mathrm{disc}}}\mE_P[\phi] \leq \alpha$ from our condition and $\phi \leq 1$. Combining the results together with the law of total expectation, we have for any $P_{X,Y,Z} \in \mathcal{P}_1$
\begin{align*}
	\mE_{P_{X,Y,Z}^n}[\phi] ~ \leq ~ & \mE_{\mathcal{L}}[	\mE_{P_{\mathrm{with}}} [\phi|\mathcal{L},\mathcal{A}] \times \mP(\mathcal{A}|\mathcal{L})] +  \mE_{\mathcal{L}}[	\mE_{P_{\mathrm{with}}} [\phi|\mathcal{L},\mathcal{A}^c] \times \mP(\mathcal{A}^c|\mathcal{L})] + \frac{n(n-1)}{J} \\[.5em]
	\leq ~ & \alpha \times \mP(\mathcal{A}) + \mP(\mathcal{A}^c) + \frac{n(n-1)}{J}.
\end{align*}
Therefore the desired result follows.

\begin{remark}  \leavevmode \normalfont
\begin{itemize} \label{Remark: Results on impossibility}
\item From Equation~(2.1) of \cite{stam1978distance}, the result still holds if we replace $\frac{n(n-1)}{J}$ in the bound~(\ref{Eq: hardness result}) with $2 - 2\frac{J!}{(J-n)!J^n}$ where one can prove that $2 - 2\frac{J!}{(J-n)!J^n} \leq  \frac{n(n-1)}{J}$ by induction.
\item Suppose that $Z$ has a uniform distribution supported on a set of size $M_n:=M > J_n:=J$. Then $\rho_{J,P}=\mP(\mathcal{A})$ is the probability that there is no collision, which can be computed as
\begin{align*}
	\rho_{J,P}= \frac{M}{M} \times \frac{M-1}{M} \times \cdots \times \frac{M - (J-1)}{M} = \frac{M!}{M^J(M-J)!}.
\end{align*}
This probability is lower and upper bounded by
\begin{align*}
	1 - \frac{J(J-1)}{2M}\leq \rho_{J,P} \leq \exp\left({-\frac{J(J-1)}{2M}}\right).
\end{align*}
See Equation~(2.4) of \cite{stam1978distance} for details. Therefore, in this case, the test $\phi$ has asymptotically no power when
\begin{align*}
	\frac{J_n^2}{M_n} \rightarrow 0 \quad \text{and} \quad \frac{n^2}{J_n} \rightarrow 0.
\end{align*}
Now let $a_n$ be a positive sequence that goes to zero arbitrarily slowly as $n \rightarrow \infty$. Then by taking $J_n^2 = M_n \times a_n$, the above conditions become equivalent to
\begin{align*}
	\frac{n^4}{M_n \times a_n} \rightarrow 0.
\end{align*}
This essentially means that any valid test has asymptotically no power in this setup if $M_n$ increases much faster than $n^4$.
\item We also note that the same proof carries over even if we assume that $(X,Y,Z)$ is supported on a compact set.
\end{itemize}
\end{remark}

\vskip 2em

\subsection{Proof of Lemma~\ref{Lemma: Error bound in terms of the Hellinger distance}}
Continuing from Lemma~\ref{Lemma: Error bound in terms of the TV distance}, we have the list of inequalities:
\begin{align*}
	 \mE_{P_{X,Y,Z}^n}[\phi_{\mathrm{perm},n}] ~\leq~ & \alpha + \mathcal{D}_{\mathrm{TV}}\big(Q_{X,Y,\widetilde{Z}}^n,\widetilde{Q}_{X,Y,\widetilde{Z}}^n\big)	\overset{\text{(i)}}{\leq} ~ \alpha + \sqrt{2} \mathcal{D}_{\mathrm{H}}\big(Q_{X,Y,\widetilde{Z}}^n,\widetilde{Q}_{X,Y,\widetilde{Z}}^n\big) \\[.5em]
	\overset{\text{(ii)}}{\leq} ~ & \alpha + \sqrt{2n} \mathcal{D}_{\mathrm{H}}\big(Q_{X,Y,\widetilde{Z}},\widetilde{Q}_{X,Y,\widetilde{Z}}\big)  \\[.5em]
	\overset{\text{(iii)}}{=}~&  \alpha +  \Bigg\{2n \sum_{m=1}^M q_{\widetilde{Z}}(m)  \times \mathcal{D}_{{H}}^2\big(Q_{X,Y|\widetilde{Z}=m},\widetilde{Q}_{X,Y|\widetilde{Z}=m}\big) \Bigg\}^{1/2},
\end{align*}
where step~(i) uses the well-known inequality $\mathcal{D}_{\mathrm{TV}}(P,Q) \leq \sqrt{2} \mathcal{D}_{\mathrm{H}}(P,Q)$ and step~(ii) use subadditivity of the squared Hellinger distance for product measures and step~(iii) can be verified by the definition of the Hellinger distance. This completes the proof of Lemma~\ref{Lemma: Error bound in terms of the Hellinger distance}.

\subsection{Proof of Theorem~\ref{Theorem: Validity of the permutation test under Hellinger smooth}} \label{Section: Proof of Theorem: Validity of the permutation test under Hellinger smooth}
We start with the case of $\gamma \in [1,2]$. In this case, the first inequality~(\ref{Eq: Hellinger inequality 1}) of Lemma~\ref{Lemma: Generalized Hellinger distance} shows that 
\begin{align*}
	& \mathcal{D}_{\mathrm{H}}^2\big(P_{X|Z=Z_m},P_{X|Z=Z_m^{'}}\big) \leq \mathcal{D}_{\gamma,\mathrm{H}}^\gamma\big(P_{X|Z=Z_m},P_{X|Z=Z_m^{'}}\big) \quad \text{and} \\[.5em]
	& \mathcal{D}_{\mathrm{H}}^2\big(P_{Y|Z=Z_m},P_{Y|Z=Z_m^{''}}\big) \leq \mathcal{D}_{\gamma,\mathrm{H}}^\gamma\big(P_{Y|Z=Z_m},P_{Y|Z=Z_m^{''}}\big). 
\end{align*}
Therefore, from Lemma~\ref{Lemma: Bound on the Hellinger distance}, it holds that
\begin{align*}
	& \mathcal{D}_{\mathrm{H}}^2 \big(Q_{X,Y|\widetilde{Z}=m},\widetilde{Q}_{X,Y|\widetilde{Z}=m}\big) \\[.5em] 
	\leq ~ & 6 \mE_{Z_m,Z_m^{'},Z_m^{''}}\big[ \mathcal{D}_{\mathrm{H}}^2\big(P_{X|Z=Z_m},P_{X|Z=Z_m^{'}}\big) \times \mathcal{D}_{\mathrm{H}}^2\big(P_{Y|Z=Z_m},P_{Y|Z=Z_m^{''}}\big)\big] \\[.5em]
	\leq ~ & 6 \mE_{Z_m,Z_m^{'},Z_m^{''}}\big[ \mathcal{D}_{\gamma,\mathrm{H}}^\gamma \big(P_{X|Z=Z_m},P_{X|Z=Z_m^{'}}\big) \times \mathcal{D}_{\gamma,\mathrm{H}}^\gamma \big(P_{Y|Z=Z_m},P_{Y|Z=Z_m^{''}}\big)\big] \\[.5em]
	\overset{\text{(i)}}{\leq} ~ & 6 L^{2\gamma} \sup_{z,z' \in B_m}\delta^{2\gamma}(z,z')  ~	\overset{\text{(ii)}}{\leq} ~ 6 L^{2\gamma} h^{2\gamma} \quad \text{for any $m \in [M]$,}
\end{align*}
where step~(i) uses $\gamma$-Hellinger Lipschitzness and step~(ii) follows by the definition of $h$ in (\ref{Eq: maximum diameter}). Applying this result to the type I error bound in terms of the Hellinger distance in Lemma~\ref{Lemma: Error bound in terms of the Hellinger distance} yields the result. 

Next we turn to the case of $\gamma > 2$. In this case, Corollary~\ref{Corollary: Hellinger bound} proves
\begin{align*}
	& \mathcal{D}_{\mathrm{H}}^2\big(P_{X|Z=Z_m},P_{X|Z=Z_m^{'}}\big) \leq C_\gamma \mathcal{D}_{\gamma,\mathrm{H}}^2 \big(P_{X|Z=Z_m},P_{X|Z=Z_m^{'}}\big) \quad \text{and} \\[.5em]
	& \mathcal{D}_{\mathrm{H}}^2\big(P_{Y|Z=Z_m},P_{Y|Z=Z_m^{''}}\big) \leq C_\gamma' \mathcal{D}_{\gamma,\mathrm{H}}^2 \big(P_{Y|Z=Z_m},P_{Y|Z=Z_m^{''}}\big),
\end{align*}
where $C_\gamma$ and $C_\gamma'$ are constants that only depend on $\gamma$. Then following the same steps as before, it holds that 
\begin{align*}
	\mathcal{D}_{\mathrm{H}}^2 \big(Q_{X,Y|\widetilde{Z}=m},\widetilde{Q}_{X,Y|\widetilde{Z}=m}\big)
	\leq ~ C_\gamma^{''} L^{2} h^{2} \quad \text{for any $m \in [M]$,}
\end{align*}
where $C_{\gamma}^{''}$ is another constant that only depends on $\gamma$. Combining this inequality with Lemma~\ref{Lemma: Error bound in terms of the Hellinger distance} yields the result. This completes the proof of Theorem~\ref{Theorem: Validity of the permutation test under Hellinger smooth}.

\subsection{Proof of Theorem~\ref{Theorem: Validity of the permutation test under Renyi smooth}}
The proof of Theorem~\ref{Theorem: Validity of the permutation test under Renyi smooth} follows exactly the same lines of the proof of Theorem~\ref{Theorem: Validity of the permutation test under Hellinger smooth} by replacing Corollary~\ref{Corollary: Hellinger bound} with Lemma~\ref{Lemma: Renyi divergence}. For this reason, we omit the details.

\vskip 2em

\subsection{Proof of Theorem~\ref{Theorem: Lower bounds}}

We proceed by considering the case of $\gamma$-Hellinger Lipschitzness (Appendix~\ref{Section: Hellinger case 1} for $2 \leq \gamma$ and Appendix~\ref{Section: Hellinger case 2} for $1 \leq \gamma < 2$) and $\gamma$-R{\'e}nyi Lipschitzness (Appendix~\ref{Section: Renyi}) in order. We start with a high-level sketch of the proof.

\begin{itemize}
	\item As discussed in Remark~\ref{Remark: Lower bounds}, we would like to show that the permutation test $\phi_{\mathrm{perm},n}^\dagger$ has power greater than $1-\beta$ for the hypothesis~(\ref{Eq: modified hypothesis}). This in turn implies our desired claim that $\phi_{\mathrm{perm},n}^\dagger$ has the type I error rate greater than $1-\beta$ for $H_0:X \independent Y | Z$. In this regard, we need to consider a powerful test that can reliably distinguish between $Q_{X,Y,\widetilde{Z}}$ and $\widetilde{Q}_{X,Y,\widetilde{Z}}$ to prove the result. To this end, we analyze the test statistic~(\ref{Eq: unweighted statistic}): 
	\begin{align*}
		T_{\mathrm{CI}} = \sum_{m \in [M]} \mathds{1}(\sigma_m \geq 4) \sigma_m U(\bW_m).
	\end{align*}
	\item However, it is quite complicated to study $T_{\mathrm{CI}}$ directly as it is a sum of dependent variables. Our strategy to overcome this complication is to construct a distributional setting where $T_{\mathrm{CI}} \approx \sigma_1 U(\bW_1)$ with high probability. More formally, we let the marginal density of $Z$ be
	\begin{align} \label{Eq: density of z}
		p_Z(z) = \begin{cases*}
			M - \epsilon, & \text{if $z \in [0, M^{-1}]$,} \\
			\frac{\epsilon}{M-1}, & \text{if $z \in (M^{-1}, 1]$,} \\
			0, & \text{otherwise,}
		\end{cases*}
	\end{align}
	for some $\epsilon:=\epsilon_n >0$, which may vary depending on $n$. Intuitively, when $\epsilon$ is sufficiently small (this will become precise in the proof), we only observe samples from the first bin with high probability and thus our test statistic will be dominated by the first term of the test statistic, i.e.
	\begin{align*}
		\sum_{m \in [M]} \mathds{1}(\sigma_m \geq 4) \sigma_m U(\bW_m)  \approx \sigma_1 U(\bW_1).
	\end{align*}
	If this is the case, then we are essentially testing for unconditional independence between $X$ and $Y$ within the first bin, and we can leverage the results of \cite{kim2020minimax} to investigate the power of the permutation test based on $\sigma_1 U(\bW_1)$. 
	\item In particular, \cite{kim2020minimax} present sufficient conditions using the first two moments of a test statistic under which the permutation test based on $U(\bW_1)$ has significant power. At a high-level, these sufficient conditions guarantee that the test statistic is greater than the critical value of the permutation test with high probability, which yields high power. 
	\item Throughout this proof, we only deal with the type II error (or $1-$power) of $\phi_{\mathrm{perm},N}$ under Poissonization since that of $\phi_{\mathrm{perm},n}^\dagger$ is upper bounded by
	\begin{equation}
	\begin{aligned} \label{Eq: bound on type II error}
		\mE[1 - \phi_{\mathrm{perm},n}^\dagger] ~\leq~ & \mE[1 - \phi_{\mathrm{perm},N}] + \mP(N > n) \\[.5em]
		\leq ~ &  \mE[1 - \phi_{\mathrm{perm},N}] + e^{-n/8},
	\end{aligned}
	\end{equation}
	where the second inequality uses an exponential Poisson tail bound \citep[e.g.][]{canonne2019note}. Therefore, once we establish that the type II error of $\phi_{\mathrm{perm},N}$ is small, it follows that the type II error of $\phi_{\mathrm{perm},n}^\dagger$ is also small by taking $n$ large enough. 
\end{itemize}

The next section formalizes it under $\gamma$-Hellinger Lipschitzness for $2 \leq \gamma$.

\subsubsection{$\gamma$-Hellinger Lipschitzness for $2 \leq \gamma$} \label{Section: Hellinger case 1}
Due to the inequality~(\ref{Eq: bound on type II error}), it suffices to show that the type II error of $\phi_{\mathrm{perm},N}$ for the hypothesis~(\ref{Eq: modified hypothesis}) is sufficiently small under the given conditions. To proceed, let us denote the entire sample by $\mathcal{W}_N=\{(X_i,Y_i,Z_i)\}_{i=1}^N$ and denote the expectation and the variance operator under the permutation law by $\mE_{\boldsymbol{\pi}}[\cdot | \boldsymbol{\sigma}, \mathcal{W}_N]$ and $\text{Var}_{\boldsymbol{\pi}}[\cdot | \boldsymbol{\sigma}, \mathcal{W}_N]$ where we recall $\boldsymbol{\sigma} = \{\sigma_1,\ldots,\sigma_M\}$. First of all, we make an observation that 
\begin{align} \label{Eq: mean zero property}
	\mE_{\boldsymbol{\pi}}[T_{\mathrm{CI}}^{\boldsymbol{\pi}} | \boldsymbol{\sigma}, \mathcal{W}_N] =  \sum_{m \in [M]} \mathds{1}(\sigma_m \geq 4) \sigma_m 	\mE_{\boldsymbol{\pi}}[U(\bW_m^{\pi_m})| \boldsymbol{\sigma}, \mathcal{W}_N] = 0,
\end{align}
where we use the fact that $\mE_{\boldsymbol{\pi}}[U(\bW_m^{\pi_m})| \boldsymbol{\sigma}, \mathcal{W}_N]=0$ when $\sigma_m \geq 4$~\citep[see Appendix I of][]{kim2020minimax}. Thus a modification of Lemma~3.1 of \cite{kim2020minimax} shows that if
\begin{align} \label{Eq: intermediate goal}
	\mE_{P_{X,Y,Z}^N,N}[T_{\mathrm{CI}}] \geq \sqrt{\frac{2\tV_{P_{X,Y,Z}^N,N}[T_{\mathrm{CI}}]}{\beta}} + \sqrt{\frac{2\mE_{P_{X,Y,Z}^N,N}\{\text{Var}_{\boldsymbol{\pi}}[T_{\mathrm{CI}}^{\boldsymbol{\pi}} | \boldsymbol{\sigma}, \mathcal{W}_N]\}}{\alpha \beta}}
\end{align} 
holds then the power of $\phi_{\mathrm{perm},N}$ is lower bounded by $1-\beta$. Therefore our goal is to prove the inequality~(\ref{Eq: intermediate goal}) under the given conditions.

In what follows, we suppress the dependence on $P_{X,Y,Z}^N,N$ in the expectation and the variance for notational convenience, and often write $U(\bW_m)$ as $U_m$ if there is no confusion. The inequality~(\ref{Eq: intermediate goal}) is verified by showing that there exists $P_{X,Y,Z} \in \mathcal{P}_{0,\mathrm{H},\gamma,\delta}(L)$ such that 
\begin{subequations}
	\begin{align} \label{Eq: step 1}
			&\mE[T_{\mathrm{CI}}] ~ \gtrsim ~ nh_n^4,  \\[.5em] \label{Eq: step 2}
			&\tV[T_{\mathrm{CI}}] ~\lesssim~ \mE[T_{\mathrm{CI}}] + 1, \\[.5em] \label{Eq: step 3}
			&\mE\{\text{Var}_{\boldsymbol{\pi}}[T_{\mathrm{CI}}^{\boldsymbol{\pi}} | \boldsymbol{\sigma}, \mathcal{W}_N]\} ~\lesssim~ 1.
	\end{align}
\end{subequations}
When these inequalities are fulfilled for large $n$, the condition~(\ref{Eq: intermediate goal}) is automatically satisfied since it is assumed that $nh_n^4 \rightarrow \infty$. Therefore our task boils down to confirming (\ref{Eq: step 1}), (\ref{Eq: step 2}) and (\ref{Eq: step 3}). 

\vskip 1em

\noindent \textbf{Proof of the inequality~(\ref{Eq: step 1}).} Let us define $p_i = \mP(Z \in B_i)$ for $i=1,\ldots,M$. Since the expectation of $U_m$ is non-negative, the expected value of $T_{\mathrm{CI}}$ can be lower bounded by
\begin{align} \label{Eq: expectation lower bound}
	\mE[T_{\mathrm{CI}}] ~\geq~ & \mE[\mathds{1}(\sigma_1 \geq 4) \sigma_1 U_1] \\[.5em]  \nonumber
	\overset{\text{(i)}}{\geq} ~ & (1 - 5e^{-1}/2) \min \{n p_1, (np_1)^4\} \|Q_{XY|\widetilde{Z}=1} -  Q_{X\cdot|\widetilde{Z}=1} Q_{\cdot Y|\widetilde{Z}=1} \|_2^2,
\end{align}
where (i) makes use of Lemma 3.1 of \cite{canonne2018testing} and the fact that $\sigma_1 \sim \text{Pois}(np_1)$. We also use the observation that $\mE[U_1 |\sigma_1] \mathds{1}(\sigma_1 \geq 4) =  \|Q_{XY|\widetilde{Z} = 1} -  Q_{X\cdot|\widetilde{Z}=1} Q_{\cdot Y|\widetilde{Z}=1} \|_2^2 \mathds{1}(\sigma_1 \geq 4)$. By our choice of the density of $Z$ in (\ref{Eq: density of z}), we have $p_1 = 1 - M^{-1}_n \epsilon_n$ and we set $\epsilon_n \leq 1/2$. In this case, $np_1 \geq n/2$ and thus
\begin{align} \label{Eq: expectation lower bound 2}
	\mE[T_{\mathrm{CI}}] ~\gtrsim~ n  \|Q_{XY|\widetilde{Z}=1} -  Q_{X\cdot|\widetilde{Z}=1} Q_{\cdot Y|\widetilde{Z}=1} \|_2^2.
\end{align}
Therefore it is sufficient to construct $P_{X|Z}$ and $P_{Y|Z}$ that satisfy (i)~$P_{X|Z}P_{Y|Z}P_Z \in \mathcal{P}_{0,\mathrm{H},\gamma,\delta}(L)$ for $\gamma \geq 2$ and (ii)~the squared $L_2$ distance between the corresponding $Q_{XY|\widetilde{Z}=1}$ and $Q_{X\cdot|\widetilde{Z}=1} Q_{\cdot Y|\widetilde{Z}=1}$ is lower bounded by
\begin{align} \label{Eq: lower bound of l2}
	 \|Q_{XY|\widetilde{Z}=1} -  Q_{X\cdot|\widetilde{Z}=1} Q_{\cdot Y|\widetilde{Z}=1} \|_2^2 \gtrsim h_n^4.
\end{align}
We construct such example below. For notational convenience, we often suppress the dependence of $n$ in $h_n$ and $M_n$, and also write $h = M^{-1}$ as we work on the equi-partition. 

\vskip 1em 

\noindent \textbf{Example of the conditional distribution in $\mathcal{P}_{0,H,\gamma,\delta}(L)$.} Suppose $X$ and $Y$ are Bernoulli random variables. Then the Lipschitz condition in Definition~\ref{Definition: Hellinger Lipschitzness} can be equivalently written as
\begin{align*}
	& \big| p_{X|Z}^{1/\gamma}(x=1|z) - p_{X|Z}^{1/\gamma}(x=1|z') \big|^\gamma +  \big|p_{X|Z}^{1/\gamma}(x=0|z) - p_{X|Z}^{1/\gamma}(x=0|z') \big|^\gamma \\[.5em]
	\leq ~ & 2L^\gamma |z - z'|^\gamma. 
\end{align*} 
For any $M \geq 1$, the following probability mass functions are well-defined:
\begin{equation}
	\begin{aligned} \label{Eq: Example of the distribution}
		p_{X|Z}(x=1|z) = \begin{cases}
			\frac{1}{2} - \frac{1}{2M} + z, \quad & \text{for $z \in [0,1/M]$,} \\[.5em]
			\frac{1}{2} + \frac{1}{2M}, \quad & \text{for $z \in (1/M,1]$,} \\[.5em]
			0, \quad & \text{otherwise,}
		\end{cases} \\[.5em]
		p_{X|Z}(x=0|z) = \begin{cases}
			\frac{1}{2} + \frac{1}{2M} - z, \quad & \text{for $z \in [0,1/M]$,} \\[.5em]
			\frac{1}{2} - \frac{1}{2M}, \quad &  \text{for $z \in (1/M,1]$,}  \\[.5em]
			0, \quad & \text{otherwise.}
		\end{cases}
	\end{aligned}
\end{equation}
\begin{figure}[t!]
	\begin{center}		
		\includegraphics[width=22em]{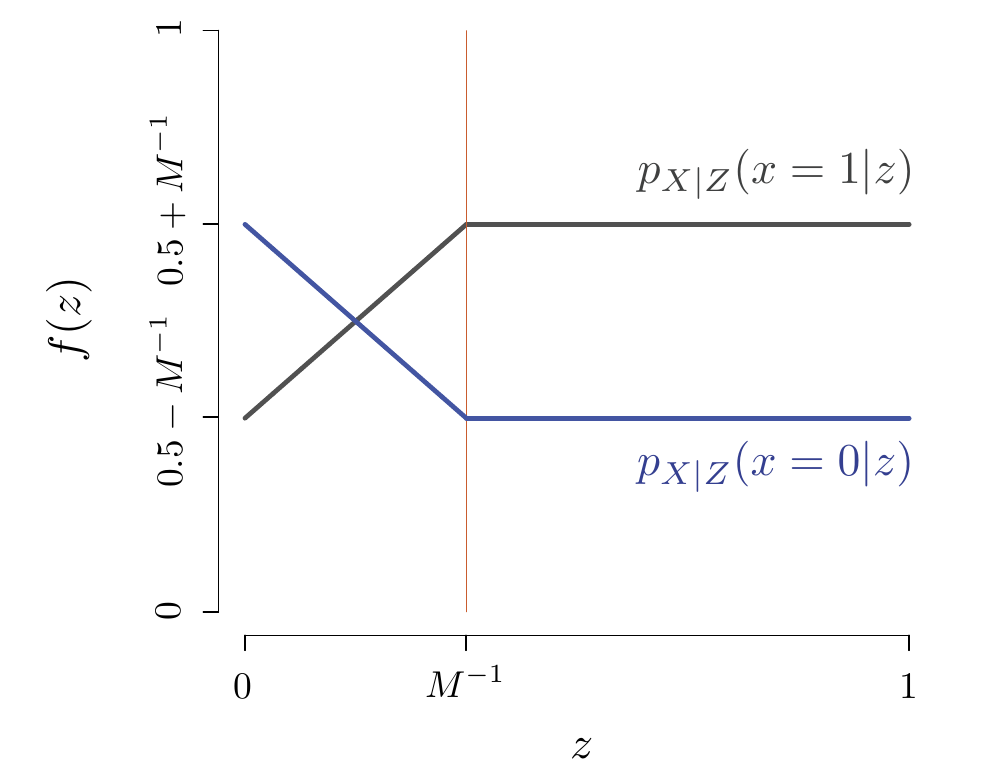} 
		\caption{\small Illustration of $p_{X|Z}(x=1|z)$ and $p_{X|Z}(x=0|z)$ described in (\ref{Eq: Example of the distribution}) as a function of $z$.} \label{Figure: density}
	\end{center}
\end{figure}
See Figure~\ref{Figure: density} for an illustration. In addition, by taking $M$ large enough, we see that both $p_{X|Z}(x=1|z)$ and $p_{X|Z}(x=0|z)$ are bounded below and also above by, let's say, $[1/4,3/4]$ for $z \in [0,1]$. Thus, using Taylor's theorem, one can verify that
\begin{align*}
	& \big| p_{X|Z}^{1/\gamma}(x=1|z) - p_{X|Z}^{1/\gamma}(x=1|z') \big|^\gamma +  \big|p_{X|Z}^{1/\gamma}(x=0|z) - p_{X|Z}^{1/\gamma}(x=0|z') \big|^\gamma \\[.5em]
	\leq~ & C_\gamma |z - z'|^\gamma,
\end{align*}
for some constant $C_\gamma > 0$. Thus it is seen that the distribution with the probability mass function~(\ref{Eq: Example of the distribution}) satisfies Definition~\ref{Definition: Hellinger Lipschitzness} with $L^\gamma = C_\gamma /2$. Similarly we let $p_{Y|Z}(y=1|z) = p_{X|Z}(x=1|z)$ and thus $p_{Y|Z}(y|z)$ is also Lipschitz in $\gamma$-Hellinger distance.

Next we calculate $q_{XY|\widetilde{Z}}(x,y|m=1)$ and $q_{X\cdot|\widetilde{Z}}(x|m=1) q_{\cdot Y|\widetilde{Z}}(y|m=1)$. To start with $q_{XY|\widetilde{Z}}(x,y|m=1)$, recall that
\begin{align*}
	d\widetilde{P}_{Z| Z \in B_1}(z) = \frac{dP_Z(z)}{\mP(Z \in [0,h])}
\end{align*}
and also recall that $Z$ has the density in (\ref{Eq: density of z}). Then since $h = M^{-1}$
\begin{align*}
	q_{XY|\widetilde{Z}}(x=1,y=1|m=1) = ~ & \int_0^{h}  p_{X|Z}(x=1|z) p_{Y|Z}(y=1|z)d\widetilde{P}_{Z| Z \in B_1}(z)  \\[.5em]
	= ~ & h^{-1} \int_0^{h} p_{X|Z}(x=1|z) p_{Y|Z}(y=1|z) dz \\[.5em] 
	=~ & h^{-1} \int_0^{h} \left( \frac{1}{2} - \frac{h}{2} + z \right)^2 dz \\[.5em]
	= ~ & \frac{h^2}{12} + \frac{1}{4}.
\end{align*}
Similarly,
\begin{align*}
	q_{X\cdot|\widetilde{Z}}(x=1|m=1)= q_{\cdot Y|\widetilde{Z}}(y=1|m=1)  = \int_0^{h}  p_{X|Z}(x=1|z) d\widetilde{P}_{Z| Z \in B_1}(z) = \frac{1}{2}.
\end{align*}
For the binary case, the squared $L_2$ distance simply becomes
\begin{align*}
 & \|Q_{XY|\widetilde{Z}=1} -  Q_{X\cdot|\widetilde{Z}=1} Q_{\cdot Y|\widetilde{Z}=1} \|_2^2 \\[.5em]
 = ~ & 4 |q_{XY|\widetilde{Z}}(x=1,y=1|m=1) - q_{X\cdot|\widetilde{Z}}(x=1|m=1) q_{\cdot Y|\widetilde{Z}}(y=1|m=1) |^2 \\[.5em]
 = ~ & \frac{h^4}{144}.
\end{align*}
Therefore the condition~(\ref{Eq: lower bound of l2}) is fulfilled and the inequality~(\ref{Eq: step 1}) follows. 

\vskip 1em

\noindent \textbf{Proof of the inequality~(\ref{Eq: step 2}).} By modifying the proof of Theorem 5.2 in \cite{neykov2020minimax}, the variance of $T_{\mathrm{CI}}$ is upper bounded by
\begin{align} \nonumber
	\text{Var}[T_{\mathrm{CI}}] ~\lesssim ~ &  \mE[T_{\mathrm{CI}}] +\mE \sum_{m \in [M]} \mathds{1}(\sigma_m \geq 4)  \\[.5em] \nonumber
	\lesssim ~ &  \mE[T_{\mathrm{CI}}] + \sum_{m \in [M]} \mP(\sigma_m \geq 1)  \\[.5em]
	\lesssim ~ &   \mE[T_{\mathrm{CI}}] + 1 - e^{-n p_1} + \sum_{m =2}^M \bigl(1 - e^{-n p_m}\bigr), \label{Eq: temp}
\end{align}
where we use the fact that $\sigma_m \sim \text{Pois}(np_m)$, and $C$ is some positive constant. Now by choosing $\epsilon_n>0$ in the density function~(\ref{Eq: density of z}) such that 
\begin{align} \label{Eq: choice of epsilon 2}
	\epsilon_n ~\leq~  \frac{M(M-1)}{n} \log \left( \frac{M-1}{M-2} \right),
\end{align}
the summation in (\ref{Eq: temp}) can be controlled as $\sum_{m =2}^M 
(1 - e^{-n p_m})  \leq 1$. In more detail, for any $m \geq 2$, the binned probability $p_m$ is the same as $\frac{\epsilon_n}{M(M-1)}$. Thus we have a series of equivalent conditions as
\begin{align*}
\sum_{m =2}^M \bigl( 1 - e^{-n p_m} \bigr)  \leq 1 ~\Longleftrightarrow~ & (M-1) (1 - e^{-\frac{n\epsilon_n}{M(M-1)}}) \leq 1 \\[.5em]
\Longleftrightarrow ~& \frac{M-2}{M-1} \leq e^{-\frac{n\epsilon_n}{M(M-1)}} \\[.5em]
\Longleftrightarrow ~& \epsilon_n \leq  \frac{M(M-1)}{n} \log \left( \frac{M-1}{M-2} \right).
\end{align*}
With the choice of such $\epsilon_n$ in the density of $Z$~(\ref{Eq: density of z}) (and also let $\epsilon_n \leq 1/2$ from the previous step), the variance is upper bounded by 
\begin{align*}
	\text{Var}[T_{\mathrm{CI}}] \lesssim \mE[T_{\mathrm{CI}}] + 1.
\end{align*}
Therefore the inequality~(\ref{Eq: step 2}) follows.

\vskip 1em

\noindent \textbf{Proof of the inequality~(\ref{Eq: step 3}).} Let us denote $U(\bW_m^{\pi_m})$ by $U_m^{\pi_m}$ for simplicity. Then using the independence between $\pi_1,\ldots,\pi_M$, the conditional variance of $T_{\mathrm{CI}}^{\boldsymbol{\pi}}$ is 
\begin{align*}
	\text{Var}_{\boldsymbol{\pi}} \big[T_{\mathrm{CI}}^{\boldsymbol{\pi}} | \boldsymbol{\sigma}, \mathcal{W}_N \big] = \sum_{m \in [M]} \mathds{1}(\sigma_m \geq 4) \sigma_m^2 \text{Var}_{\boldsymbol{\pi}}[ U_m^{\pi_m} | \boldsymbol{\sigma}, \mathcal{W}_N]. 
\end{align*}
Based on the result in Appendix J of \cite{kim2020minimax},
\begin{align*}
	\mE\{ \text{Var}_{\boldsymbol{\pi}}[ U_m^{\pi_m} | \boldsymbol{\sigma}, \mathcal{W}_N] | \boldsymbol{\sigma}\}  ~\lesssim ~ & \sigma_m^{-2} \max \big\{ \|Q_{XY|\widetilde{Z}=m}\|_2^2, \|Q_{X\cdot|\widetilde{Z}=m} Q_{\cdot Y|\widetilde{Z}=m} \|_2^2 \big\} \\[.5em]
	\lesssim ~ &  \sigma_m^{-2},
\end{align*}
where the second inequality holds since $\|P\|_2^2 \leq \|P\|_1^2 = 1$ for a discrete probability distribution $P$. Therefore, with the choice of $\epsilon_n$ as in (\ref{Eq: choice of epsilon 2}), 
\begin{align*}
	\mE\{\text{Var}_{\boldsymbol{\pi}}[T_{\mathrm{CI}}^{\boldsymbol{\pi}} | \boldsymbol{\sigma}, \mathcal{W}_N]\} ~\lesssim ~ \sum_{m \in [M]} \mP(\sigma_m \geq 4) ~ \lesssim ~  1,
\end{align*}
where the last inequality holds since $\sigma_m \sim \text{Pois}(np_m)$ and
\begin{align*}
\sum_{m \in [M]} \mP(\sigma_m \geq 4)  \leq \sum_{m \in [M]} \mP(\sigma_m \geq 1) \leq 1 + \sum_{m =2}^M \bigl( 1 - e^{-n p_m} \bigr) \leq 2.
\end{align*}
This verifies the inequality~(\ref{Eq: step 3}). 

\vskip 1em
 
Indeed, our proof in this subsection goes through for the case where $1 \leq \gamma < 2$ as well. However, the condition $\sqrt{n}h_n^2 \rightarrow \infty$ is not tight for this case, and we need a more careful analysis detailed in the next subsection.

\vskip 1em 

\subsubsection{$\gamma$-Hellinger Lipschitzness for $1 \leq \gamma < 2$}  \label{Section: Hellinger case 2}
The underlying idea of the proof of this result is the same as in the previous section. That is, we want to show that the inequality~(\ref{Eq: intermediate goal}) holds when $\sqrt{n} h_n^\gamma \rightarrow \infty$ for some $P_{X,Y,Z} \in \mathcal{P}_{0,\mathrm{H},\gamma,\delta}(L)$ with $\gamma \in [1,2]$. However, we need to obtain a sharper bound for the variance of the test statistic in order to achieve the desired rate, which requires some nontrivial effort.

Suppose that $X$ and $Y$ are binary random variables taking a value between $\{0,1\}$. We then set 
\begin{equation}
	\begin{aligned} \label{Eq: conditional densities under gamma hellinger}
		& p_{X|Z}(x=1|z) = 
		\begin{cases}
			z^\gamma, \quad & \text{for $z \in [0,2^{-1/\gamma}]$,} \\[.5em]
			2^{-1}, \quad & \text{for $z \in [2^{-1/\gamma},1]$,} \\[.5em]
			0, \quad & \text{otherwise.}
		\end{cases} \\[.5em]
		& p_{X|Z}(x=0|z) = 
		\begin{cases}
			1-z^\gamma, \quad & \text{for $z \in [0,2^{-1/\gamma}]$,} \\[.5em]
			2^{-1}, \quad & \text{for $z \in [2^{-1/\gamma},1]$,} \\[.5em]
			0, \quad & \text{otherwise.}
	\end{cases}
	\end{aligned}
\end{equation}
Similarly, we let $p_{X|Z}(x=1|z) = p_{Y|Z}(y=1|z)$ and $p_{X|Z}(x=0|z) = p_{Y|Z}(y=0|z)$. It can be checked that $P_{X|Z}$ and $P_{Y|Z}$ with the above probability mass functions are $\gamma$-Hellinger smooth for $\gamma \geq 1$ using Taylor's theorem. Here we note that the Lipschitz constant depends on $\gamma$ and is assumed to be fixed. Also note that this example is more sophisticated than the one in (\ref{Eq: Example of the distribution}). For instance, when $\gamma = 1$, the distribution~(\ref{Eq: conditional densities under gamma hellinger}) is TV smooth but not smooth in the Hellinger distance (since $\sqrt{x}$ is not a Lipschitz function). In addition to the condition on the conditional distributions~(\ref{Eq: conditional densities under gamma hellinger}), we let $Z$ have the marginal density as in (\ref{Eq: density of z}) with a sufficiently small $0 < \epsilon_n \leq 1/2$ (specified later).

Under this setting, we will show that the mean, the variance and the expected conditional variance of the test statistic $T_{\mathrm{CI}}$ are bounded by
\begin{subequations}
\begin{align} \label{Eq: step 1a}
	& \mE[T_{\mathrm{CI}}] ~\gtrsim~ n h_n^{4\gamma}, \\[.5em] \label{Eq: step 2a}
	& \mathrm{Var}[T_{\mathrm{CI}}] ~\lesssim~ \mE[T_{\mathrm{CI}}]  h_n^{2\gamma} + h_n^{4\gamma}, \\[.5em] \label{Eq: step 2c}
	& \mE \big\{ \mathrm{Var}_{\boldsymbol{\pi}} \big[T_{\mathrm{CI}}^{\boldsymbol{\pi}} | \boldsymbol{\sigma}, \mathcal{W}_N] \big] \big\} ~\lesssim ~  h_n^{4\gamma}.
\end{align}  
\end{subequations}
If this is the case, then the inequality~(\ref{Eq: intermediate goal}) holds as $\sqrt{n} h_n^\gamma \rightarrow \infty$, which completes the proof. In what follows, we prove these inequalities in order. 

\vskip 1em

\noindent \textbf{Proof of the inequality~(\ref{Eq: step 1a}).} As we define the density of $Z$ to be (\ref{Eq: density of z}) with $\epsilon_n \leq 1/2$, the inequality~(\ref{Eq: expectation lower bound 2}) follows in this case as well, that is
\begin{align} \label{Eq: lower bond for the expectation}
		\mE[T_{\mathrm{CI}}] ~\gtrsim~ n  \|Q_{XY|\widetilde{Z}=1} -  Q_{X\cdot|\widetilde{Z}=1} Q_{\cdot Y|\widetilde{Z}=1} \|_2^2.
\end{align}
Therefore it is enough to show that $\|Q_{XY|\widetilde{Z}=1} -  Q_{X\cdot|\widetilde{Z}=1} Q_{\cdot Y|\widetilde{Z}=1} \|_2^2 \gtrsim h^{4\gamma}$. Then the inequality~(\ref{Eq: step 1a}) follows. In fact, when $h \leq 1/2$, 
\begin{align*}
	q_{XY|\widetilde{Z}}(x=1,y=1|m=1) ~= ~ & \int_0^{h}  p_{X|Z}(x=1|z) p_{Y|Z}(y=1|z)d\widetilde{P}_{Z| Z \in B_1}(z)  \\[.5em]
	= ~ & h^{-1} \int_0^{h} p_{X|Z}(x=1|z) p_{Y|Z}(y=1|z) dz \\[.5em] 
	=~ & h^{-1} \int_0^{h} z^{2 \gamma} dz = \frac{h^{2\gamma}}{2\gamma + 1}.
\end{align*}
Similarly,
\begin{align*}
	q_{X\cdot|\widetilde{Z}}(x=1|m=1)= q_{\cdot Y|\widetilde{Z}}(y=1|m=1)  = \int_0^{h}  p_{X|Z}(x=1|z) d\widetilde{P}_{Z| Z \in B_1}(z) = \frac{h^\gamma}{\gamma+1}.
\end{align*}
When $X$ and $Y$ are binary, the squared $L_2$ norm simply becomes 
\begin{align}  \nonumber 
	 & \|Q_{XY|\widetilde{Z}=1} -  Q_{X\cdot|\widetilde{Z}=1} Q_{\cdot Y|\widetilde{Z}=1} \|_2^2 \\[.5em] \label{Eq: alternative expression of l2 norm}
	 = ~ & 4 \big\{ q_{XY|\widetilde{Z}}(x=1,y=1|m=1) - q_{X\cdot|\widetilde{Z}}(x=1|m=1) q_{\cdot Y|\widetilde{Z}}(y=1|m=1) \big\}^2 \\[.5em] \nonumber 
	 = ~ & 4h^{4\gamma}  \left( \frac{1}{2\gamma + 1} - \frac{1}{(\gamma + 1)^2} \right)^2, 
\end{align}
which concludes the inequality~(\ref{Eq: step 1a}). 

\vskip 1em

\noindent \textbf{Proof of the inequality~(\ref{Eq: step 2a}).} Bounding the variance requires a bit more work. First the law of total variance yields
\begin{align} \label{Eq: variance decomposition}
	\mathrm{Var}[T_{\mathrm{CI}}] =  \mE\{\mathrm{Var}[T_{\mathrm{CI}}| \boldsymbol{\sigma} ]\} +  \mathrm{Var}\{\mE[T_{\mathrm{CI}}| \boldsymbol{\sigma} ]\}.
\end{align}
To start with the first term, we use the independence between different bins conditional on $\boldsymbol{\sigma} $ and see
\begin{align*}
	\mE\{\mathrm{Var}[T_{\mathrm{CI}}| \boldsymbol{\sigma} ]\} ~=~ & \sum_{m \in [M]} \mE\{\mathds{1}(\sigma_m \geq 4) \sigma_m^2 \mathrm{Var}[U_m| \boldsymbol{\sigma}]\} \\[.5em]
	\leq ~ & \mE\{\mathds{1}(\sigma_1 \geq 4) \sigma_1^2 \mathrm{Var}[U_1| \sigma_1]\} + \sum_{m=2}^M \mE[\mathds{1}(\sigma_m \geq 4) \sigma_m^2],
\end{align*}
where the second inequality uses $|U_m| \leq 1$. Since we have freedom to choose $\epsilon$ in the density of $Z$, we can make the second term smaller than 
\begin{align} \label{Eq: condition on eps}
	\sum_{m=2}^M \mE[\mathds{1}(\sigma_m \geq 4) \sigma_m^2] ~ \lesssim ~h^{4\gamma}.
\end{align}
For example, since $\sigma_m \sim \text{Pois}(np_m)$ where $p_m= \frac{\epsilon}{M(M-1)}$ for $m=2,\ldots,M$, applying Cauchy--Schwarz inequality yields 
\begin{align*}
	\sum_{m=2}^M \mE[\mathds{1}(\sigma_m \geq 4) \sigma_m^2] ~ \lesssim ~ & M  \sqrt{\mP(\sigma_2 \geq 4) \mE[\sigma_2^4]} \\[.5em]
	\lesssim ~ &  M \sqrt{1 - e^{-np_2}},
\end{align*}
where the second inequality follows since $\mP(\sigma_2 \geq 4) \leq \mP(\sigma_2 \geq 1) = 1 - e^{-np_2}$ and $\mE[\sigma_2^4] = np_2(1 + 7np_2 + 6n^2p_2^2 + n^3 p_2^3) \leq 15$ by assuming $np_2 \leq 1$, which is satisfied when $\epsilon \lesssim 1/(nh^2)$. Therefore, the inequality~(\ref{Eq: condition on eps}) holds, for instance, when
\begin{align*}
	\epsilon ~\lesssim~ \frac{1}{nh^2} \log \left(  \frac{1}{1 - h^{4\gamma+1}}\right).
\end{align*} 

Next we analyze $ \mE\{\mathds{1}(\sigma_1 \geq 4) \sigma_1^2 \mathrm{Var}[U_1| \sigma_1]\}$. Intuitively the conditional variance of $U_1|\sigma_1$ becomes smaller for large $M$ (equivalently small $h$) under our distributional setting~(\ref{Eq: conditional densities under gamma hellinger}). In particular, observe that $q_{XY|\widetilde{Z}}(x=0,y=0|m=1) \rightarrow 1$ and thus $(X,Y)|Z \in B_1$ becomes degenerate at $X=0,Y=0$ as $M$ increases. However the bound for $\mathrm{Var}[U_1 | \sigma_1]$ given in Lemma 5.1 of \cite{neykov2020minimax} \citep[and also see the proof of Proposition 5.3 in][]{kim2020minimax}:
\begin{align*}
	\mathrm{Var}[U_1 | \sigma_1] ~ \lesssim ~&  \frac{\mE[U_1|\sigma_1] \max \big\{\|Q_{XY|\widetilde{Z}=1}\|_2, \| Q_{X\cdot|\widetilde{Z}=1} Q_{\cdot Y|\widetilde{Z}=1}\|_2 \big\}}{\sigma_1} \\[.5em]
	+ ~ & \frac{\max\big\{\|Q_{XY|\widetilde{Z}=1}\|_2^2, \| Q_{X\cdot|\widetilde{Z}=1} Q_{\cdot Y|\widetilde{Z}=1} \|_2^2 \big\}}{\sigma_1^2} \quad \text{for $\sigma_1 \geq 4$,}
\end{align*}
does not capture this intuition. It turns out that this variance bound is not tight for the binary case. In fact, we can obtain a sharper bound as
\begin{equation}
\begin{aligned} \label{Eq: sharper bound}
	\mathrm{Var}[U_1 | \sigma_1] ~ \lesssim ~&  \frac{\mE[U_1|\sigma_1]\max \big\{ q_{XY|\widetilde{Z}}(1,1|1), q_{X\cdot|\widetilde{Z}}(1|1) q_{\cdot Y|\widetilde{Z}}(1|1) \big\}}{\sigma_1} \\[.5em]
	+ ~ & \frac{\max \big\{ q_{XY|\widetilde{Z}}^2(1,1|1), q_{X\cdot|\widetilde{Z}}^2(1|1) q_{\cdot Y|\widetilde{Z}}^2(1|1) \big\}}{\sigma_1^2} \quad \text{for $\sigma_1 \geq 4$,}
\end{aligned}
\end{equation}
where $q_{XY|\widetilde{Z}}(1,1|1) = q_{XY|\widetilde{Z}}(x=1,y=1|m=1)$ and the other quantities are similarly defined. We will prove this inequality at the end of this section. Therefore, using the previous calculations of $q_{XY|\widetilde{Z}}$, $q_{X\cdot|\widetilde{Z}}$ and $q_{\cdot Y|\widetilde{Z}}$, the (conditional) variance is bounded by
\begin{align} \label{Eq: variance upper bound}
	\mathrm{Var}[U_1 | \sigma_1] ~\lesssim~ \frac{ \mE[U_1|\sigma_1] h^{2 \gamma} }{\sigma_1} +  \frac{h^{4 \gamma}}{\sigma_1^2} \quad \text{for $\sigma_1 \geq 4$}.
\end{align}
In summary, we have
\begin{equation}
\begin{aligned} \label{Eq: upper bound for the first term}
	\mE\{\mathrm{Var}[T_{\mathrm{CI}}| \boldsymbol{\sigma} ]\}  ~ \lesssim ~ & \mE\{\mathds{1}(\sigma_1 \geq 4) \sigma_1 \mE[U_1|\sigma_1]\} h^{2 \gamma}  + h^{4\gamma} \\[.5em]
	\lesssim ~ & \mE[T_{\mathrm{CI}}] h^{2 \gamma}  + h^{4\gamma}. 
\end{aligned}
\end{equation}

Next we focus on the second term $\mathrm{Var}\{\mE[T_{\mathrm{CI}}| \boldsymbol{\sigma} ]\}$ in (\ref{Eq: variance decomposition}) and see 
\begin{align*}
	\mathrm{Var}\{ \mE[T_{\mathrm{CI}}| \boldsymbol{\sigma}]\} ~=~ & \sum_{m \in [M]} \mathrm{Var}[\sigma_m \mathds{1}(\sigma_m \geq 4)] \times \|Q_{XY|\widetilde{Z}=m} -  Q_{X\cdot|\widetilde{Z}=m} Q_{\cdot Y|\widetilde{Z}=m} \|_2^4 \\[.5em]
	\overset{\text{(i)}}{\lesssim} ~ & \sum_{m \in [M]} \mE[\sigma_m \mathds{1}(\sigma_m \geq 4)] \times \|Q_{XY|\widetilde{Z}=m} -  Q_{X\cdot|\widetilde{Z}=m} Q_{\cdot Y|\widetilde{Z}=m} \|_2^4 \\[.5em]
	\overset{\text{(ii)}}{\lesssim} ~ & np_1 \times \big|q_{X,Y|\widetilde{Z}}(1,1|1) - q_{X\cdot|\widetilde{Z}}(1|1) q_{\cdot Y|\widetilde{Z}}(1|1)\big|^4 + h^{4\gamma} \\[.5em]
	\overset{\text{(iii)}}{\lesssim} ~ & \mE[T_{\mathrm{CI}}] h^{2\gamma} + h^{4\gamma},
\end{align*}
where step~(i) uses Claim 2.1 of \cite{canonne2018testing}, step~(ii) can be deduced by the inequality~(\ref{Eq: condition on eps}) and the identity~(\ref{Eq: alternative expression of l2 norm}). To give a more detail on step~(ii), consider a decomposition
\begin{align*}
& \sum_{m \in [M]} \mE[\sigma_m \mathds{1}(\sigma_m \geq 4)] \times \|Q_{XY|\widetilde{Z}=m} -  Q_{X\cdot|\widetilde{Z}=m} Q_{\cdot Y|\widetilde{Z}=m} \|_2^4 \\[.5em]
= ~ & \underbrace{\mE[\sigma_1 \mathds{1}(\sigma_1 \geq 4)] \times \|Q_{XY|\widetilde{Z}=1} -  Q_{X\cdot|\widetilde{Z}=1} Q_{\cdot Y|\widetilde{Z}=1} \|_2^4}_{\mathrm{(I)}} \\[.5em]
&+ \underbrace{\sum_{m=2}^M \mE[\sigma_m \mathds{1}(\sigma_m \geq 4)] \times \|Q_{XY|\widetilde{Z}=m} -  Q_{X\cdot|\widetilde{Z}=m} Q_{\cdot Y|\widetilde{Z}=m} \|_2^4}_{\mathrm{(II)}}.
\end{align*}
For the first term~(I), we make use of the identity~(\ref{Eq: alternative expression of l2 norm}) and $\mE[\sigma_1] = np_1$, which yields $\mathrm{(I)} \lesssim np_1 \times \big|q_{X,Y|\widetilde{Z}}(1,1|1) - q_{X\cdot|\widetilde{Z}}(1|1) q_{\cdot Y|\widetilde{Z}}(1|1)\big|^4$. For the second term, note that the $L_2$ norm of $Q_{XY|\widetilde{Z}=m} -  Q_{X\cdot|\widetilde{Z}=m} Q_{\cdot Y|\widetilde{Z}=m}$ is bounded by two for any $m$. Therefore, $\mathrm{(II)} \lesssim \sum_{m=2}^M \mE[\sigma_m^2 \mathds{1}(\sigma_m \geq 4)] \lesssim h^{4\gamma}$ by the inequality~(\ref{Eq: condition on eps}). This explains step~(ii). Step~(iii) follows by using the lower bound for $\mE[T_{\mathrm{CI}}]$ in (\ref{Eq: lower bond for the expectation}) and the observation that $|q_{XY|\widetilde{Z}}(1,1|1)-q_{X\cdot|\widetilde{Z}}(1|1) q_{\cdot Y|\widetilde{Z}}(1|1)|^2 \lesssim h^{4\gamma} \leq h^{2\gamma}$. Combining this with (\ref{Eq: upper bound for the first term}) yields the claim~(\ref{Eq: step 2a}). The last step is to verify the sharper bound~(\ref{Eq: sharper bound}).

\vskip 1em

\noindent \textbf{Verification of the bound~(\ref{Eq: sharper bound}).} Recall that the kernel of $U_1$ is 
\begin{align*}
	h_{i_1,i_2,i_3,i_4}= \frac{1}{4!} \sum_{\pi \in \boldsymbol{\Pi}_4} \sum_{x \in [\ell_1], y \in [\ell_2]} \psi_{\pi(1) \pi(2)}(x,y) \psi_{\pi(3)\pi(4)}(x,y),
\end{align*}
where $\boldsymbol{\Pi}_4$ is the set of all permutations of $\{i_1,i_2,i_3,i_4\}$ and 
\begin{align*}
	\psi_{ij}(x,y) = \mathds{1}(X_{i,1}=x,Y_{i,1}=y) -  \mathds{1}(X_{i,1}=x) \mathds{1}(Y_{j,1}=y).
\end{align*}
Importantly, for the binary case, we observe that 
\begin{align} \label{Eq: another kernel}
	h_{i_1,i_2,i_3,i_4} = g_{i_1,i_2,i_3,i_4}:= \frac{4}{4!} \sum_{\pi \in \boldsymbol{\Pi}_4}  \psi_{\pi(1)\pi(2)}(1,1) \psi_{\pi(3)\pi(4)}(1,1). 
\end{align}
This means that $U_1$ is the same as another U-statistic defined with the kernel $g_{i_1,i_2,i_3,i_4}$. Having this observation, we work with the U-statistic associated with $g_{i_1,i_2,i_3,i_4}$. We then follow basically the same lines of the proofs in \cite{neykov2020minimax} and Appendix J of \cite{kim2020minimax} and see that the bound~(\ref{Eq: sharper bound}) is satisfied (more specifically, follow the same lines of the proof in Appendix J of \cite{kim2020minimax} by setting $d_1=d_2=1$).

\vskip 1em

\noindent \textbf{Proof of the inequality~(\ref{Eq: step 2c}).} Similar to the proof of (\ref{Eq: step 3}), we can obtain 
\begin{align*}
	\mathrm{Var}_{\boldsymbol{\pi}} \big[T_{\mathrm{CI}}^{\boldsymbol{\pi}} | \boldsymbol{\sigma}, \mathcal{W}_N \big] ~=~ & \sum_{m \in [M]} \mathds{1}(\sigma_m \geq 4) \sigma_m^2 \mathrm{Var}_{\boldsymbol{\pi}}[ U_m^{\pi_m} | \boldsymbol{\sigma}, \mathcal{W}_N] \\
	\leq ~ & \mathds{1}(\sigma_1 \geq 4) \sigma_1^2 \mathrm{Var}_{\boldsymbol{\pi}}[ U_1^{\pi_1} | \boldsymbol{\sigma}, \mathcal{W}_N] +  \sum_{m = 2}^M \mathds{1}(\sigma_m \geq 4) \sigma_m^2,
\end{align*}
where the second inequality follows by $|U_m^{\pi_m}|\leq 1$. By taking the expectation, the second term can be smaller than $h^{4\gamma}$ by the condition~(\ref{Eq: condition on eps}). On the other hand, recall that for the binary case, $U_1$ is equivalent to the U-statistic based on the kernel~$g_{i_1,i_2,i_3,i_4}$ in (\ref{Eq: another kernel}). Thus, based on the results in Appendix I and Appendix J of \cite{kim2020minimax} (more specifically by letting $d_1=d_2=1$ therein),
\begin{align*}
	\mE\{ \mathrm{Var}_{\boldsymbol{\pi}}[ U_1^{\pi_1} | \boldsymbol{\sigma}, \mathcal{W}_N] | \boldsymbol{\sigma}\}  ~\lesssim ~ & \sigma_1^{-2} \max \big\{ q_{XY|\widetilde{Z}}^2(1,1|1), q_{X\cdot|\widetilde{Z}}^2(1|1) q_{\cdot Y|\widetilde{Z}}^2(1|1) \big\} \\[.5em]
	\lesssim ~ & \sigma_1^{-2} h^{4\gamma}. 
\end{align*}
Consequently, we obtain the desired bound $\mE \big\{ \mathrm{Var}_{\boldsymbol{\pi}} \big[T_{\mathrm{CI}}^{\boldsymbol{\pi}} | \boldsymbol{\sigma}, \mathcal{W}_N] \big] \big\} ~\lesssim ~  h_n^{4\gamma}$.

\subsubsection{$\gamma$-R{\'e}nyi Lipschitzness for $0 < \gamma$}  \label{Section: Renyi}
The proof of this part is essentially the same as that in the case of $\gamma$-Hellinger Lipschitzness for $2 \leq \gamma$ (Appendix~\ref{Section: Hellinger case 1}). In particular, we consider the distribution setting in (\ref{Eq: Example of the distribution}) and show that the considered distribution is R{\'e}nyi smooth for all $\gamma > 0$. Once this is established, the other parts of the proof follow exactly the same lines of the proof as before in Appendix~\ref{Section: Hellinger case 1}.

Note that, for a sufficiently large $M$, $p_{X|Z}(x=1|z)$ and $p_{X|Z}(x=0|z)$ are bounded below and above by constants for any $z \in [0,1]$, say,
\begin{align} \label{Eq: bound on pmf}
	p_{X|Z}(x=1|z) \in [1/4,3/4]\quad \text{and} \quad p_{X|Z}(x=0|z) \in [1/4,3/4].
\end{align}
Therefore, when $0 < \gamma \leq 2$, there exist some constants $C_1,C_2>0$ such that 
\begin{align*}
		\mathcal{D}_{\gamma\text{,\emph{R}}}(P_{X|Z=z} \| P_{X|Z=z'}) ~ \overset{\text{(i)}}{\leq} ~ & \mathcal{D}_{\chi^2} (P_{X|Z=z} \| P_{X|Z=z'}) \\[.5em]
	\overset{\text{(ii)}}{\leq} ~ & C_1 \| P_{X|Z=z} - P_{X|Z=z'}\|_2^2 \\[.5em]
	\overset{\text{(iii)}}{\leq} ~ & C_2 |z - z'|^2,
\end{align*}
where step~(i) uses the fact that R{\'e}nyi divergence is nondecreasing in $\gamma$~\citep[e.g.~Theorem 3 of][]{van2014renyi} and the inequality~(\ref{Eq: Renyi inequality}). Additionally, step~(ii) holds since we set $P_{X|Z=z'}$ is lower bounded by some constant and step~(iii) follows directly by the construction of $P_{X|Z=z}$ in (\ref{Eq: Example of the distribution}). Similarly we let $P_{Y|Z=z} = P_{X|Z=z}$ and thus both $P_{X|Z=z}$ and $P_{Y|Z=z}$ are R{\'e}nyi smooth when $0 < \gamma \leq 2$.

We now consider the case where $\gamma > 2$. Again, we make use of the distribution given in (\ref{Eq: Example of the distribution}) and assume that $p_{X|Z}(x=1|z)$ and $p_{X|Z}(x=0|z)$ are bounded below and above by constants for any $z \in [0,1]$ as in (\ref{Eq: bound on pmf}). Since $X$ takes a binary value, R{\'e}nyi divergence of $P_{X|Z=z}$ from $P_{X|Z=z'}$ simply becomes 
\begin{align*}
	& \mathcal{D}_{\gamma\text{,\emph{R}}}(P_{X|Z=z} \| P_{X|Z=z'}) \\[.5em]
	 = ~ & \frac{1}{\gamma - 1} \log \biggl\{ \biggl(  \frac{p_{X|Z}(x=1|z)}{p_{X|Z}(x=1|z')}\biggr)^\gamma p_{X|Z}(x=1|z') +  \biggl(  \frac{p_{X|Z}(x=0|z)}{p_{X|Z}(x=0|z')}\biggr)^\gamma p_{X|Z}(x=0|z') \biggr\}.
\end{align*}
When both $z,z' \in (1/M,1]$, R{\'e}nyi divergence becomes zero and hence there is nothing to prove. Next assume that both $z,z' \in [0,1/M]$. In this case, R{\'e}nyi divergence can be written more explicitly as
\begin{align*}
	& \mathcal{D}_{\gamma\text{,\emph{R}}}(P_{X|Z=z} \| P_{X|Z=z'})  \\[.5em]
  = ~ & \frac{1}{\gamma - 1}  \log \biggl\{ \biggl( \frac{\frac{1}{2} - \frac{1}{2M} + z}{ \frac{1}{2} - \frac{1}{2M} + z' }\biggr)^\gamma \biggl(\frac{1}{2} - \frac{1}{2M} + z'\biggr) +  \biggl( \frac{\frac{1}{2} + \frac{1}{2M} - z}{ \frac{1}{2} + \frac{1}{2M} - z' }\biggr)^\gamma \biggl(\frac{1}{2} + \frac{1}{2M} - z'\biggr) \biggr\} \\[.5em]
  := ~ & f(z).
\end{align*}
By viewing the above expression as a function of $z$, Taylor's theorem yields that there exists $\xi \in (0,1/M)$ such that 
\begin{align*}
	f(z) = f(z') + f'(z')(z-z') + f''(\xi)(z-z')^2.
\end{align*}
The first term $f(z')$ is zero. We also note that the first derivative satisfies $f'(z')=0$ since
\begin{align*}
	f'(z) = \frac{1}{\gamma - 1} \times \frac{\gamma \biggl( \frac{\frac{1}{2} - \frac{1}{2M} + z}{ \frac{1}{2} - \frac{1}{2M} + z' }\biggr)^{\gamma-1} - \gamma \biggl( \frac{\frac{1}{2} + \frac{1}{2M} - z}{ \frac{1}{2} + \frac{1}{2M} - z' }\biggr)^{\gamma-1}}{ \biggl( \frac{\frac{1}{2} - \frac{1}{2M} + z}{ \frac{1}{2} - \frac{1}{2M} + z' }\biggr)^\gamma \biggl(\frac{1}{2} - \frac{1}{2M} + z'\biggr) +  \biggl( \frac{\frac{1}{2} + \frac{1}{2M} - z}{ \frac{1}{2} + \frac{1}{2M} - z' }\biggr)^\gamma \biggl(\frac{1}{2} + \frac{1}{2M} - z'\biggr)}.
\end{align*}
Furthermore, it can be seen that $\sup_{z \in [0,1/M]}|f''(z)| \leq C_{\gamma,M}$ where $C_{\gamma,M}$ is a positive constant depending only on $\gamma$ and $M$. This shows that $ \mathcal{D}_{\gamma\text{,\emph{R}}}(P_{X|Z=z} \| P_{X|Z=z'}) \leq C_{\gamma,M}(z - z')^2$ for all $z,z' \in [0,1/M]$. Suppose now that $z \in [0,1/M]$ and $z' \in(1/M,1]$. In this case, since $p_{X|Z}(x=1|z')= p_{X|Z}(x=z|1/M)$ and $p_{X|Z}(x=0|z') = p_{X|Z}(x=0|1/M)$ for all $z' \in (1/M,1]$, the previous result implies that $ \mathcal{D}_{\gamma\text{,\emph{R}}}(P_{X|Z=z} \| P_{X|Z=z'}) \leq C_{\gamma,M}(z - 1/M)^2 \leq C_{\gamma,M}(z - z')^2$. The final case where $z \in (1/M,1]$ and $z' \in [0,1/M]$ can be similarly handled, which concludes that $ \mathcal{D}_{\gamma\text{,\emph{R}}}(P_{X|Z=z} \| P_{X|Z=z'}) \leq C_{\gamma,M}(z - z')^2$ for all $z,z' \in [0,1]$. Hence, by letting $P_{Y|Z=z} = P_{X|Z=z}$, both $P_{X|Z=z}$ and $P_{Y|Z=z}$ are R{\'e}nyi smooth for $\gamma \geq 2$.

In summary, the distribution in (\ref{Eq: Example of the distribution}) is R{\'e}nyi smooth for $\gamma > 0$ with a different Lipschitz constant and therefore the result $\sup_{P_{X,Y,Z} \in \mathcal{P}_{0, \text{\emph{R}},\gamma,\delta}(L)}  \mE_{P_{X,Y,Z}^N,N}[\phi_{\text{{perm}},n}^\dagger] \geq  1 - \beta$ follows.

\vskip 2em

\subsection{Proof of Theorem~\ref{Theorem: Type II error under discrete setting}} \label{Section: Proof of Theorem: Type II error under discrete setting}
We start by proving the type II error bound for $\phi_{\mathrm{perm},2}$ and then turn to $\phi_{\mathrm{perm},1}$. 

\vskip 1em 

\noindent \textbf{$\bullet$ Type II error of $\phi_{\mathrm{perm},2}$.} We note that by Lemma~\ref{Lemma: quantile}, the type II error of $\phi_{\mathrm{perm},2}$ can be expressed as 
\begin{align*}
	\mE_{P_{X,Y,Z}^N,N}[1 - \phi_{\mathrm{perm},2}]  = \mP_{P_{X,Y,Z}^N,N}[T_{\mathrm{CI},W} \leq q_{1-\alpha}],
\end{align*}
where $q_{1-\alpha}$ is the $1-\alpha$ quantile of the empirical distribution of $T_{\mathrm{CI},W}^{\boldsymbol{\pi}_1},\ldots,T_{\mathrm{CI},W}^{\boldsymbol{\pi}_K}$. Therefore, it is enough to work with the test that rejects the null when $T_{\mathrm{CI},W} > q_{1-\alpha}$. For simplicity, we use the notation $\boldsymbol{\mathcal{W}}_{\mathrm{weight}} := \{\bW_{X,m}, \bW_{Y,m}\}_{m\in[M]}$ to denote samples used for weights. Given this notation, we first claim that 
\begin{align} \nonumber
	& \mP_{P_{X,Y,Z}^N,N}[T_{\mathrm{CI},W} \leq q_{1-\alpha}] \\[.5em]
	\leq~ & \mP_{P_{X,Y,Z}^N,N}\Big\{T_{\mathrm{CI},W} < \zeta\sqrt{\min\{n,M\} + \mE_{P_{X,Y,Z}^N,N}[T_{\mathrm{CI},W}| \boldsymbol{\mathcal{W}}_{\mathrm{weight}}, \boldsymbol{\sigma}]}\Big\} + \frac{1}{200}, \label{Eq: first goal}
\end{align}
where $\zeta$ is a sufficiently large constant as in the critical value of the test in \cite{canonne2018testing}. We then show that under the given condition, 
\begin{align} \label{Eq: second goal}
	\mP_{P_{X,Y,Z}^N,N}\Big\{T_{\mathrm{CI},W} < \zeta\sqrt{\min\{n,M\} + \mE_{P_{X,Y,Z}^N,N}[T_{\mathrm{CI},W}| \boldsymbol{\mathcal{W}}_{\mathrm{weight}}, \boldsymbol{\sigma}]}\Big\} \leq \frac{1}{200},
\end{align}
which establishes the desired result. To ease the notation, we suppress dependence on $P_{X,Y,Z}^N$ and $N$ in the expectation and probability function, throughout this proof. 

\vskip 1em

\noindent \textbf{1. Proof of the first claim~(\ref{Eq: first goal}).} To show the first claim~(\ref{Eq: first goal}), let us define an event 
\begin{align*}
	\mathcal{A} := \{ q_{1-\alpha} \leq \zeta\sqrt{\min(n,M) + \mE_{P_{X,Y,Z}^N,N}[T_{\mathrm{CI},W}| \boldsymbol{\mathcal{W}}_{\mathrm{weight}}, \boldsymbol{\sigma}]} \},
\end{align*}
and show that
\begin{align} \label{Eq: goal}
	\mP(\mathcal{A}) \geq \frac{199}{200}.
\end{align}
If this is the case, then the first claim~(\ref{Eq: first goal}) follows by the union bound. To begin, let us simply denote the permutation distribution of $T_{\mathrm{CI},W}^{\boldsymbol{\pi}}$ by
\begin{align*}
	\mathbb{P}_{\boldsymbol{\pi}}\big[T_{\mathrm{CI},W}^{\boldsymbol{\pi}} \leq t \big] = \frac{1}{K} \sum_{\boldsymbol{\pi}_i \in \boldsymbol{\Pi}} \mathds{1}\big[ T_{\mathrm{CI},W}^{\boldsymbol{\pi}_i} \leq t \big].
\end{align*}
More generally, we let $\mathbb{P}_{\boldsymbol{\pi}}[\cdot]$ be the probability measure in terms of $\boldsymbol{\pi}$ uniformly distributed over $\boldsymbol{\Pi}$ conditional on everything else. We denote the expectation and the variance operator under the permutation law by $\mE_{\boldsymbol{\pi}}[\cdot]$ and $\mathrm{Var}_{\boldsymbol{\pi}}[\cdot]$ for simplicity. As in the case of (\ref{Eq: mean zero property}) for $T_{\mathrm{CI}}^{\boldsymbol{\pi}}$, we make an observation that $T_{\mathrm{CI},W}^{\boldsymbol{\pi}}$ is centered under the permutation law, i.e.
\begin{align*}
	\mE_{\boldsymbol{\pi}}[T_{\mathrm{CI},W}^{\boldsymbol{\pi}}] =  \sum_{m \in [M]} \mathds{1}(\sigma_m \geq 4) \sigma_m  \omega_m\mE_{\boldsymbol{\pi}}[U_W(\bW_m^{\pi_m})] = 0,
\end{align*}
which can be shown using the result in Appendix I of \cite{kim2020minimax}. Therefore, for any $t>0$, Chebyshev's inequality yields
\begin{align*}
	\mathbb{P}_{\boldsymbol{\pi}}\big\{T_{\mathrm{CI},W}^{\boldsymbol{\pi}} \geq t \big\} ~\leq~ \frac{1}{t^2} \mE_{\boldsymbol{\pi}}\big[(T_{\mathrm{CI},W}^{\boldsymbol{\pi}})^2\big].
\end{align*}
Then, by the definition of the $1-\alpha$ quantile, it holds that $q_{1-\alpha} \leq \sqrt{\alpha^{-1}\mE_{\boldsymbol{\pi}}\big[(T_{\mathrm{CI},W}^{\boldsymbol{\pi}})^2\big]}$ with probability one. Notice that $\mE_{\boldsymbol{\pi}}\big[(T_{\mathrm{CI},W}^{\boldsymbol{\pi}})^2\big]$ is a non-negative random variable, changing its value depending on $\{\bW_{X,m}, \bW_{Y,m}, \bW_{XY,m}\}_{m=1}^M$. Then Markov's inequality verifies that the following inequality 
\begin{align*}
	\mE_{\boldsymbol{\pi}}\big[(T_{\mathrm{CI},W}^{\boldsymbol{\pi}})^2\big] <  400 \mE \big\{  \mE_{\boldsymbol{\pi}}\big[(T_{\mathrm{CI},W}^{\boldsymbol{\pi}})^2\big] | \boldsymbol{\mathcal{W}}_{\mathrm{weight}}, \boldsymbol{\sigma} \big\}
\end{align*}
holds with probability at least $1-1/400$. Based on the results of \cite{canonne2018testing}, \cite{neykov2020minimax} and \cite{kim2020minimax}, we further prove below that 
\begin{align} \label{Eq: expectation bound}
	400\alpha^{-1} \mE \big\{  \mE_{\boldsymbol{\pi}}\big[(T_{\mathrm{CI},W}^{\boldsymbol{\pi}})^2\big] | \boldsymbol{\mathcal{W}}_{\mathrm{weight}}, \boldsymbol{\sigma} \big\} \leq \zeta^{1/2} \big( \min \{n,M\} + \mE[T_{\mathrm{CI},W}| \boldsymbol{\mathcal{W}}_{\mathrm{weight}}, \boldsymbol{\sigma}]\big),
\end{align}
with probability at least $1 - 1/400$ by taking $\zeta$ sufficiently large (depending on $\alpha$, but recall that $\alpha$ is assumed to be a fixed constant). This verifies the probability bound~(\ref{Eq: goal}) for $\mathcal{A}$, completing the proof of the first claim~(\ref{Eq: first goal}).

\vskip 1em

\noindent \textbf{Detail of the inequality~(\ref{Eq: expectation bound}).} First of all, since $\pi_1,\ldots,\pi_M$ are independent, it can be verified that
\begin{align} \label{Eq: conditional variance}
	\mE_{\boldsymbol{\pi}}\big[(T_{\mathrm{CI},W}^{\boldsymbol{\pi}})^2\big]  = \sum_{m \in [M]} \mathds{1}(\sigma_m \geq 4) \sigma_m^2 \omega_m^2 \mE_{\boldsymbol{\pi}}\big[U_W^2(\bW_m^{\pi_m})\big].
\end{align}
In order to analyze the expectation of $\mE_{\boldsymbol{\pi}}\big[U_W^2(\bW_m^{\pi_m})\big]$, we need some preliminary results. To this end, for $m=1,\ldots,M$, let 
\begin{align*}
	& g_{X,m}(x_1,x_2) := \sum_{x \in [\ell_1]} \frac{\mathds{1}(x_1 = x) \mathds{1}(x_2 = x)}{1+a_{x,m}} \quad \text{and} \\[.5em]
	& g_{Y,m}(y_1,y_2) := \sum_{y \in [\ell_2]} \frac{\mathds{1}(y_1 = y) \mathds{1}(y_2 = y)}{1+a_{y,m}},
\end{align*}
and for $p:[\ell_1] \times [\ell_2] \mapsto \mathbb{R}$, denote 
\begin{align*}
	\|p\|_{2,\boldsymbol{a}_m}^2 :=  \sum_{x \in [\ell_1]} \sum_{y \in [\ell_2]} \frac{p^2(x,y)}{(1+a_{x,m})(1+a_{y,m})}.
\end{align*}
Since $a_{x,m}, a_{y,m} \geq 0$, the above bivariate functions are non-negative and bounded by one. Therefore
\begin{align*}
	& \mE[g_{X,m}^2(X_{1,m},X_{2,m})g_{Y,m}^2(Y_{1,m},Y_{2,m})|\bW_{X,m}, \bW_{Y,m} ] \\[.5em]
	\leq ~ & \mE[g_{X,m}(X_{1,m},X_{2,m})g_{Y,m}(Y_{1,m},Y_{2,m})|\bW_{X,m}, \bW_{Y,m} ] ~=~ \|P_{XY|Z=m}\|_{2,\boldsymbol{a}_m}^2.
\end{align*}
Similar calculations show
\begin{align*}
	\mE[g_{X,m}^2(X_{1,m},X_{2,m})g_{Y,m}^2(Y_{3,m},Y_{4,m})|\bW_{X,m}, \bW_{Y,m} ] ~ \leq ~  \|P_{X|Z=m}P_{Y|Z=m}\|_{2,\boldsymbol{a}_m}^2
\end{align*}
and
\begin{align*}
	& \mE[g_{X,m}^2(X_{1,m},X_{2,m})g_{Y,m}^2(Y_{1,m},Y_{3,m})|\bW_{X,m}, \bW_{Y,m} ] \\[.5em]
	\leq ~ &  \frac{1}{2}\|P_{XY|Z=m}\|_{2,\boldsymbol{a}_m}^2 + \frac{1}{2} \|P_{X|Z=m}P_{Y|Z=m}\|_{2,\boldsymbol{a}_m}^2,
\end{align*}
where the last inequality uses the inequality $xy \leq \frac{1}{2}x^2 + \frac{1}{2}y^2$. In addition, using the triangle inequality, we have
\begin{align*}
	\|P_{XY|Z=m}\|_{2,\boldsymbol{a}_m}^2 ~\leq~ 2\|P_{X,Y|m} - P_{X|Z=m}P_{Y|Z=m}\|_{2,\boldsymbol{a}_m}^2 + 2	\|P_{X|Z=m}P_{Y|Z=m}\|_{2,\boldsymbol{a}_m}^2.
\end{align*}
Since we consider half-permutation method~(Remark~\ref{Remark: full- vs half-permutation}) for $T_{\mathrm{CI},W}$, we can apply the bound for the expected variance of a permuted U-statistic in \cite{kim2020minimax} to $U_W^2(\bW_m^{\pi_m})$. More specifically, by the argument made in Appendix I of \cite{kim2020minimax} along with the above results, we have
\begin{align*}
	& \mathds{1}(\sigma_m \geq 4)  \mE\big[ \mE_{\boldsymbol{\pi}}\big[U_W^2(\bW_m^{\pi_m})\big] | \bW_{X,m}, \bW_{Y,m} ,\sigma_m \big] \\[.5em]
	\lesssim ~ &  \mathds{1}(\sigma_m \geq 4) \sigma_m^{-2} \max \big\{ \|P_{XY|Z=m}\|_{2,\boldsymbol{a}_m}^2, \|P_{X|Z=m}P_{Y|Z=m}\|_{2,\boldsymbol{a}_m}^2\big\}   \\[.5em]
	\lesssim ~  &  \mathds{1}(\sigma_m \geq 4) \sigma_m^{-2} \big\{ \|P_{X,Y|m} - P_{X|Z=m}P_{Y|Z=m}\|_{2,\boldsymbol{a}_m}^2 +  \|P_{X|Z=m}P_{Y|Z=m}\|_{2,\boldsymbol{a}_m}^2 \big\}.
\end{align*}
This result together with~(\ref{Eq: conditional variance}) leads to the following bound
\begin{align*}
	\mE \Big[ \mE_{\boldsymbol{\pi}}\big[(T_{\mathrm{CI},W}^{\boldsymbol{\pi}})^2\big] | \boldsymbol{\mathcal{W}}_{\mathrm{weight}}, \boldsymbol{\sigma}\big] \Big] ~ \lesssim ~  \big\{ (\mathrm{I}) + (\mathrm{II}) \big\},
\end{align*}
where
\begin{align*}
	& (\mathrm{I}) ~:=~ \sum_{m \in [M]} \mathds{1}(\sigma_m \geq 4) \omega_m^2 \|P_{X|Z=m}P_{Y|Z=m}\|_{2,\boldsymbol{a}_m}^2, \\[.5em]
	& (\mathrm{II}) ~:=~ \sum_{m \in [M]} \mathds{1}(\sigma_m \geq 4) \omega_m^2 \|P_{X,Y|m} - P_{X|Z=m}P_{Y|Z=m}\|_{2,\boldsymbol{a}_m}^2.
\end{align*}
To simplify the first term, the result of \cite{canonne2018testing} and \cite{neykov2020minimax} yields
\begin{align*}
	\mE[\|P_{X|Z=m}P_{Y|Z=m}\|_{2,\boldsymbol{a}_m}^2]  ~\leq~ \frac{1}{(1 + t_{1,m})(1+t_{2,m})},
\end{align*}
which further shows that  
\begin{align*}
	\mE\big[(\mathrm{I})| \boldsymbol{\sigma} \big] ~ \leq ~ \sum_{m \in [M]} \mathds{1}(\sigma_m \geq 4) \frac{\omega_m^2}{(1 + t_{1,m})(1+t_{2,m})} ~\lesssim ~ \min\{N,M\}. 
\end{align*}
Recall that $N \sim \text{Pois}(n)$, which highly concentrates around its mean. Thus, by Markov's inequality as well as the union bound, one can have $(\mathrm{I}) \lesssim  \min\{n,M\}$ with probability at least $1 - 1/400$. For the second term, since $\omega_m \leq \sigma_m$, we have $(\mathrm{II}) \leq \mE\big[T_{\mathrm{CI},W}| \boldsymbol{\mathcal{W}}_{\mathrm{weight}}, \boldsymbol{\sigma} \big]$. Therefore, with probability at least $1 - 1/400$, it holds that 
\begin{align*}
	(\mathrm{I}) + (\mathrm{II})  ~ \lesssim ~ \min\{n,M\} + \mE\big[T_{\mathrm{CI},W}| \boldsymbol{\mathcal{W}}_{\mathrm{weight}}, \boldsymbol{\sigma} \big],
\end{align*} 
which proves the inequality~(\ref{Eq: expectation bound}).

\vskip 1em

\noindent \textbf{2. Proof of the second claim~(\ref{Eq: second goal}).} Next we prove the second claim~(\ref{Eq: second goal}) under the condition on $\epsilon$ imposed in \cite{canonne2018testing}. In particular, for $\ell_1 \geq \ell_2$ and $\inf_{Q \in \mathcal{P}_{0,[M]}}\mathcal{D}_{\mathrm{TV}}(P_{X,Y,Z}, Q) \geq \varepsilon$, equation~(1) of \cite{canonne2018testing} assumes that $\varepsilon$  satisfies
\begin{equation}
\begin{aligned} \label{Eq: condition on varepsilon}
	n ~ \geq ~ \zeta' \times \max \Biggl( & \min \Biggl\{ \frac{M^{7/8}\ell_1^{1/4}\ell_2^{1/4} }{\varepsilon}, \frac{M^{6/7} \ell_1^{2/7}\ell_2^{2/7} }{\varepsilon^{8/7}}  \Biggr\}, \\[.5em]
	& \frac{M^{3/4}\ell_1^{1/2}\ell_2^{1/2} }{\varepsilon},  \frac{M^{2/3}\ell_1^{2/3}\ell_2^{1/3} }{\varepsilon^{4/3}},\frac{M^{1/2}\ell_1^{1/2}\ell_2^{1/2} }{\varepsilon^{2}}  \Biggr),
\end{aligned}
\end{equation}
where $\zeta'$ is a sufficiently large positive constant. A sufficient condition for this inequality in the form of $\varepsilon \gtrsim f(M,\ell_1,\ell_2,n)$ is 
\begin{align*}
	\varepsilon ~\geq ~ \zeta'' \times \max \Biggl(  \frac{M^{7/8}\ell_1^{1/4}\ell_2^{1/4}}{n}, & \frac{M^{3/4}\ell_1^{1/2}\ell_2^{1/2}}{n}, \frac{M^{1/2} \ell_1^{1/2} \ell_2^{1/4}}{n^{3/4}}, \\[.5em]
	& \frac{M^{1/4}\ell_1^{1/4}\ell_2^{1/4}}{n^{1/2}},  \frac{M^{3/4} \ell_1^{1/4} \ell_2^{1/4}}{n^{7/8}}  \Biggr),
\end{align*}
where $\zeta''$ is some large positive constant. Let us define an event $\mathcal{E}$ such that $\mE\big[T_{\mathrm{CI},W}| \boldsymbol{\mathcal{W}}_{\mathrm{weight}}, \boldsymbol{\sigma} \big] \gtrsim  \sqrt{\zeta' \min\{n,M\}}$ and 
\begin{align*}
	& \mathrm{Var}\big[T_{\mathrm{CI},W}| \boldsymbol{\mathcal{W}}_{\mathrm{weight}}, \boldsymbol{\sigma} \big] \\[.5em] 
	\lesssim ~ & \min\{n,M\} + \sqrt{\min\{n,M\}} \mE\big[T_{\mathrm{CI},W}| \boldsymbol{\mathcal{W}}_{\mathrm{weight}}, \boldsymbol{\sigma} \big] +  \mE\big[T_{\mathrm{CI},W}| \boldsymbol{\mathcal{W}}_{\mathrm{weight}}, \boldsymbol{\sigma} \big]^{3/2}.
\end{align*}
If $\inf_{Q \in \mathcal{P}_{0,[M]}}\mathcal{D}_{\mathrm{TV}}(P_{X,Y,Z}, Q) \geq \varepsilon$, Lemma 5.4 of \cite{canonne2018testing} guarantees that $\mathbb{P}(\mathcal{E}) \geq 399/400$. Under this event, by choosing $\zeta'$ sufficiently large depending on $\zeta$, we see that 
\begin{align*}
	\zeta\sqrt{\min\{n,M\} + \mE[T_{\mathrm{CI},W}| \boldsymbol{\mathcal{W}}_{\mathrm{weight}}, \boldsymbol{\sigma}]} ~\leq~ \frac{1}{2} \mE[T_{\mathrm{CI},W}| \boldsymbol{\mathcal{W}}_{\mathrm{weight}}, \boldsymbol{\sigma}].
\end{align*}
Thus 
\begin{align*}
	& \mP \Big[T_{\mathrm{CI},W} < \zeta\sqrt{\min\{n,M\} + \mE[T_{\mathrm{CI},W}| \boldsymbol{\mathcal{W}}_{\mathrm{weight}}, \boldsymbol{\sigma}]} \Big] \\[.5em]
	\leq ~ & \mP \Big[ T_{\mathrm{CI},W} \leq  \frac{1}{2} \mE[T_{\mathrm{CI},W}| \boldsymbol{\mathcal{W}}_{\mathrm{weight}}, \boldsymbol{\sigma}] | \mathcal{E}  \Big] + \frac{1}{400} \\[.5em]
	\leq ~ & \frac{1}{200},
\end{align*}
where the last inequality uses Chebyshev's inequality, as detailed in the proof of Lemma 5.6 in \cite{canonne2018testing}.

\vskip 1em
 
\noindent \textbf{$\bullet$ Type II error of $\phi_{\mathrm{perm},1}$.} The type II error bound for $\phi_{\mathrm{perm},1}$ follows similarly as $\phi_{\mathrm{perm},2}$. Indeed it is simpler to prove the result of $\phi_{\mathrm{perm},1}$ than $\phi_{\mathrm{perm},2}$ as $T_{\mathrm{CI}}$ does not involve randomness from weights. We note from Lemma~\ref{Lemma: quantile} again that the type II error of $\phi_{\mathrm{perm},1}$ is equivalent to $\mP_{P_{X,Y,Z}^N,N}[T_{\mathrm{CI}} \leq q_{1-\alpha}]$ where $q_{1-\alpha}$ is the $1-\alpha$ quantile of the empirical distribution of $T_{\mathrm{CI}}^{\boldsymbol{\pi}_1},\ldots,T_{\mathrm{CI}}^{\boldsymbol{\pi}_K}$ (with some abuse of notation). We then prove the results in two steps as before. In the first step, we show that
\begin{align} \label{Eq: first claim}
	\mP_{P_{X,Y,Z}^N,N}[T_{\mathrm{CI}} \leq q_{1-\alpha}] 	~\leq~  \mP_{P_{X,Y,Z}^N,N}\big[T_{\mathrm{CI}} < \zeta\sqrt{\min\{n,M\}}\big] + \frac{1}{200}.
\end{align}
In the second step, we verify that 
\begin{align}  \label{Eq: second claim}
	\mP_{P_{X,Y,Z}^N,N}\big[T_{\mathrm{CI}} < \zeta\sqrt{\min\{n,M\}}\big]  \leq \frac{1}{200},
\end{align}
under the condition on $\varepsilon$ imposed by \cite{canonne2018testing}. 

\vskip 1em

\noindent \textbf{1. Proof of the first claim~(\ref{Eq: first claim}).} To start with the first claim~(\ref{Eq: first claim}), by Chebyshev's inequality along with Markov's inequality, we have 
\begin{align*}
	q_{1-\alpha} \leq \sqrt{200\alpha^{-1} \mE\big\{ \mE_{\boldsymbol{\pi}}\big[(T_{\mathrm{CI}}^{\boldsymbol{\pi}})^2\big]\big\}},
\end{align*}
with probability at least $1 - 1/200$. Again, by the independence between $\pi_1,\ldots,\pi_M$ and following the analysis in the proof of Theorem~\ref{Theorem: Lower bounds} --- especially the proof of the inequality~(\ref{Eq: step 3}), the conditional expectation becomes 
\begin{equation}
	\begin{aligned} \label{Eq: expected variance of permuted U-stat}
		\mE\big\{ \mE_{\boldsymbol{\pi}}\big[(T_{\mathrm{CI}}^{\boldsymbol{\pi}})^2\big] | \boldsymbol{\sigma} \big\} ~=~ & \sum_{m \in [M]} \mathds{1}(\sigma_m \geq 4) \sigma_m^2 \mE\big[U_W^2(\bW_m^{\pi_m})| \boldsymbol{\sigma}\big] \\[.5em]
		~\overset{\mathrm{(i)}}{=}~& \sum_{m \in [M]} \mathds{1}(\sigma_m \geq 4) \sigma_m^2 \mE\big[\mathrm{Var}_{\boldsymbol{\pi}} \{U_W(\bW_m^{\pi_m})| \bW_m, \boldsymbol{\sigma}\}| \boldsymbol{\sigma} \big] \\[.5em]
		\overset{\mathrm{(ii)}}{\lesssim} ~ &  \sum_{m \in [M]} \mathds{1}(\sigma_m \geq 1) \overset{\mathrm{(iii)}}{\lesssim}  \min\{N,M\},
	\end{aligned}
\end{equation}
where step~(i) uses the fact that $\mE_{\boldsymbol{\pi}} \{U_W(\bW_m^{\pi_m})| \bW_m, \boldsymbol{\sigma}\} = 0$, step~(ii) follows based on the analysis in  the proof of the inequality~(\ref{Eq: step 3}) and step~(iii) follows by combining the two inequalities $\sum_{m \in [M]} \mathds{1}(\sigma_m \geq 1) \leq \sum_{m \in [M]} 1 = M$ and $\sum_{m \in [M]} \mathds{1}(\sigma_m \geq 1) \leq \sum_{m \in [M]} \sigma_m = N$. By taking the expectation over $\boldsymbol{\sigma}$, we then have $\mE\big\{ \mE_{\boldsymbol{\pi}}\big[(T_{\mathrm{CI}}^{\boldsymbol{\pi}})^2\big]\big\} \lesssim \min\{n,M\}$. Therefore, the inequality~(\ref{Eq: first claim}) holds by taking $\zeta$ sufficiently large and the union bound. 

\vskip 1em

\noindent \textbf{2. Proof of the second claim~(\ref{Eq: second claim}).} The proof of this result directly follows by Section 3.1 of \cite{canonne2018testing}. In particular, for fixed $\ell_1$ and $\ell_2$, and $\inf_{Q \in \mathcal{P}_{0,[M]}}\mathcal{D}_{\mathrm{TV}}(P_{X,Y,Z}, Q) \geq \varepsilon$, \cite{canonne2018testing} assume that $\varepsilon$ satisfies
\begin{align}  \label{Eq: condition on varepsilon 2}
	n ~ \geq ~\zeta' \times \max \Biggl( \frac{M^{1/2}}{\varepsilon^2}, \min \Biggl( \frac{M^{7/8}}{\varepsilon}, \frac{M^{6/7}}{\varepsilon^{8/7}} \Biggr) \Biggr),
\end{align}
where $\zeta'$ is a sufficiently large positive constant. A sufficient condition for this inequality in the form of $\varepsilon \gtrsim f(M,n)$ is 
\begin{align*}
	\varepsilon ~ \geq ~ \zeta'' \max \Biggl( \frac{M^{1/4}}{n^{1/2}}, \frac{M^{7/8}}{n}, \frac{M^{3/4}}{n^{7/8}} \Biggr),
\end{align*}
for some large $\zeta''>0$ as mentioned in the main text. Under this condition, Section 3.1 of \cite{canonne2018testing} shows that the second claim~(\ref{Eq: second claim}) holds by taking $\zeta$ sufficiently large. Finally, combining the two inequalities~(\ref{Eq: first claim}) and (\ref{Eq: second claim}) yields the desired result. This completes the proof of Theorem~\ref{Theorem: Type II error under discrete setting}.

\vskip 2em

\subsection{Proof of Theorem~\ref{Theorem: Type II error under continuous setting I}}  \label{Section: Proof of Theorem: Type II error under continuous setting I}
We first study the type II error bound of $\phi_{\mathrm{perm},1}^\dagger$ and then turn to the type II error bound of $\phi_{\mathrm{perm},2}^\dagger$. 

\vskip 1em

\noindent \textbf{$\bullet$ Type II error of $\phi_{\mathrm{perm},1}^\dagger$.} The proof of this part is similar to that of Theorem~\ref{Theorem: Type II error under discrete setting}, especially the result of $\phi_{\mathrm{perm},1}$. The only difference is that we are now dealing with the discretized probability $Q_{X,Y,\widetilde{Z}}$ defined in Section~\ref{Section: Validity}, rather than the original probability $P_{X,Y,Z}$. Since the upper bound for the expected variance $\mE\big\{ \mE_{\boldsymbol{\pi}}\big[(T_{\mathrm{CI}}^{\boldsymbol{\pi}})^2\big] | \boldsymbol{\sigma} \big\}$ in (\ref{Eq: expected variance of permuted U-stat}) does not depend on the underlying discrete distribution, the same proof carries through and we can obtain the bound~(\ref{Eq: first claim}) for $Q_{X,Y,\widetilde{Z}}$ as well. In particular, since we set $\min\{n,M\} \leq M = \ceil{n^{2/5}}$, we have
\begin{align*}
	\mE_{P_{X,Y,Z}^N,N}[1 - \phi_{\mathrm{perm},1}^\dagger]  ~ \leq ~ \mP_{Q_{X,Y,\widetilde{Z}}^N,N}\big[T_{\mathrm{CI}} < \zeta n^{1/5}\big] + \frac{1}{200} + e^{-n/8},
\end{align*}
where we note that $\mP_{Q_{X,Y,\widetilde{Z}}^N,N}\big[T_{\mathrm{CI}} < \zeta n^{1/5}\big] = \mP_{P_{X,Y,Z}^N,N}\big[T_{\mathrm{CI}} < \zeta n^{1/5}\big]$ since $T_{\mathrm{CI}}$ is defined only through the binned data. We also note that $e^{-n/8}$ comes from the truncation probability of $N$ as in~(\ref{Eq: bound on type II error}). This probability is essentially the type II error of the test of \cite{neykov2020minimax} by adjusting the value of $\zeta$, and thus we have
\begin{align*}
		\sup_{P_{X,Y,Z} \in \mathcal{P}_{1,[0,1],\mathrm{TV}}(L): \inf_{Q \in \mathcal{P}_{0,[0,1]}} \mathcal{D}_{\mathrm{TV}}(P_{X,Y,Z},Q) \geq \varepsilon} \mP_{P_{X,Y,Z}^N,N} \big[T_{\mathrm{CI}} < \zeta n^{1/5}\big] ~\leq~ \frac{1}{100}.
\end{align*}
See the proof of Theorem 5.2 in \cite{neykov2020minimax} for more details. Therefore, we conclude that $\phi_{\mathrm{perm},1}^\dagger$ has the same power guarantee as $\phi_{\mathrm{NBW},1}$.

\vskip 1em

\noindent \textbf{$\bullet$ Type II error of $\phi_{\mathrm{perm},2}^\dagger$.} Given that $\phi_{\mathrm{perm},2}^\dagger$ uses the same test statistic $T_{\mathrm{CI},W}$ as in Theorem~\ref{Theorem: Type II error under discrete setting}, the proof of this result is similar to that of Theorem~\ref{Theorem: Type II error under discrete setting}, especially the result of $\phi_{\mathrm{perm},2}$. The only difference is, again, that we now work with the discretized probability distribution $Q_{X,Y,\widetilde{Z}}$ defined in Section~\ref{Section: Validity}. Following exactly the same lines in the proof of Theorem~\ref{Theorem: Type II error under discrete setting}, there exists an absolute constant $\zeta' > 0$ such that 
\begin{align*}
	&\mE_{P_{X,Y,Z}^N,N}[1 - \phi_{\mathrm{perm},2}^\dagger]   \\[.5em]
	\leq~ & \mP_{Q_{X,Y,\widetilde{Z}}^N,N}\Big\{T_{\mathrm{CI},W} < \zeta'\sqrt{\min\{n,M\} + \mE_{Q_{X,Y,\widetilde{Z}}^N,N}[T_{\mathrm{CI},W} | \boldsymbol{\mathcal{W}}_{\mathrm{weight}}, \boldsymbol{\sigma}]}\Big\} + \frac{1}{200} + e^{-n/8},
\end{align*}
where the second line also holds by replacing $Q_{X,Y,\widetilde{Z}}^N$ with $P_{X,Y,Z}^N$ in the probability and the expectation since $T_{\mathrm{CI},W}$ is defined only through the binned data. 

To further bound the given probability, we recall the results in \cite{neykov2020minimax}. By ignoring the dependence on $Q_{X,Y,\widetilde{Z}}^N$, the following two bounds hold with probability at least $399/400$ under the given conditions for the guarantee~(\ref{Eq: type II error 2}):
\begin{align*}
	& \mathrm{Var}\big[T_{\mathrm{CI},W} | \boldsymbol{\mathcal{W}}_{\mathrm{weight}}, \boldsymbol{\sigma} \big] \\[.5em]
	\lesssim ~ &  M + (\sqrt{M} +1) \mE[T_{\mathrm{CI},W} | \boldsymbol{\mathcal{W}}_{\mathrm{weight}}, \boldsymbol{\sigma}]  + [T_{\mathrm{CI},W} | \boldsymbol{\mathcal{W}}_{\mathrm{weight}}, \boldsymbol{\sigma}]^{3/2}
\end{align*}
and 
\begin{align*}
	\mE[T_{\mathrm{CI},W} | \boldsymbol{\mathcal{W}}_{\mathrm{weight}}, \boldsymbol{\sigma}] \gtrsim  \sqrt{\zeta M},
\end{align*}
where $\zeta > 0$ is a sufficiently large constant given in $\phi_{\mathrm{NBW},2}$, which differs from $\zeta'$ above. See Lemma C.10 of \cite{neykov2020minimax} for details on these bounds. Let $\mathcal{E}$ be an event that both inequalities hold simultaneously. Under this event, by taking $\zeta$ sufficiently large depending on $\zeta'$, we have
\begin{align*}
	\zeta'\sqrt{\min\{n,M\} + \mE[T_{\mathrm{CI},W} | \boldsymbol{\mathcal{W}}_{\mathrm{weight}}, \boldsymbol{\sigma}]} ~\leq~ \frac{1}{2}\mE[T_{\mathrm{CI},W} | \boldsymbol{\mathcal{W}}_{\mathrm{weight}}, \boldsymbol{\sigma}].
\end{align*}
Then by Chebyshev's inequality and Lemma C.13 of \cite{neykov2020minimax},
\begin{align*}
	& \mP \Big[T_{\mathrm{CI},W} < \zeta'\sqrt{\min\{n,M\} + \mE[T_{\mathrm{CI},W}| \boldsymbol{\mathcal{W}}_{\mathrm{weight}}, \boldsymbol{\sigma}]} \Big] \\[.5em]
	\leq ~ & \mP \Big[ T_{\mathrm{CI},W} \leq  \frac{1}{2} \mE[T_{\mathrm{CI},W}| \boldsymbol{\mathcal{W}}_{\mathrm{weight}}, \boldsymbol{\sigma}] | \mathcal{E}  \Big] + \frac{1}{400} \\[.5em]
	\leq ~ & \frac{1}{200}.
\end{align*}
This completes the proof of the second part. 

\vskip 2em

\subsection{Proof of Theorem~\ref{Theorem: Type II error under continuous setting II}}
The proof of this result is essentially the same as that of the second part of Theorem~\ref{Theorem: Type II error under continuous setting I}, which builds on the proof of Theorem~\ref{Theorem: Type II error under discrete setting}. The only difference is that we are now in a situation where $X,Y,Z$ are all discretized. Under this binning scheme, the mean and the variance of $T_{\mathrm{CI},W}$ have been studied in \cite{neykov2020minimax}. In particular, Lemma G.9 of \cite{neykov2020minimax} yields that the following two events hold under the conditions of the theorem with probability at least $399/400$:
\begin{align*}
	& \mathrm{Var}\big[T_{\mathrm{CI},W} | \boldsymbol{\mathcal{W}}_{\mathrm{weight}}, \boldsymbol{\sigma} \big] \\[.5em]
	\lesssim ~ &  \big\{ M + (\sqrt{M} +1) \mE[T_{\mathrm{CI},W} | \boldsymbol{\mathcal{W}}_{\mathrm{weight}}, \boldsymbol{\sigma}]  + [T_{\mathrm{CI},W} | \boldsymbol{\mathcal{W}}_{\mathrm{weight}}, \boldsymbol{\sigma}]^{3/2} \big\}
\end{align*}
and 
\begin{align*}
	\mE[T_{\mathrm{CI},W} | \boldsymbol{\mathcal{W}}_{\mathrm{weight}}, \boldsymbol{\sigma}] \gtrsim \sqrt{\zeta M},
\end{align*}
where $\zeta > 0$ is a sufficiently large constant given in $\phi_{\mathrm{NBW},2}$. Building upon this result, we can proceed by following the same lines of the proof of Theorem~\ref{Theorem: Type II error under continuous setting I} to complete the proof. We omit the details.

\vskip 2em

\subsection{Proof of Proposition~\ref{Proposition: Type I error of the double-binned test}}
The proof of this result follows by modifying the proof of the validity of usual permutation tests based on exchangeability~\citep[e.g.][]{lehmann2006testing}. Here we provide a detail for completeness. Let $\mathcal{X}:= \{(X_i,Y_i,\widetilde{Z}_i)\}_{i=1}^n \sim \widetilde{Q}_{X,Y,\widetilde{Z}}^n$ be a set of $n$ i.i.d.~draws from the CI projection~$\widetilde{Q}_{X,Y,\widetilde{Z}}^n$. For each $\boldsymbol{\pi}, \boldsymbol{\pi}' \in \boldsymbol{\Pi}_{\mathrm{cyclic}}$, write $\mathcal{X}^{\boldsymbol{\pi}} = \{(X_i,Y_{\boldsymbol{\pi}(i)},\widetilde{Z}_i)\}_{i=1}^n$ and $\mathcal{X}^{\boldsymbol{\pi} \circ \boldsymbol{\pi}'} = \{(X_i,Y_{\boldsymbol{\pi} \circ \boldsymbol{\pi}'(i)},\widetilde{Z}_i)\}_{i=1}^n$ where $\boldsymbol{\pi} \circ \boldsymbol{\pi}'$ denotes the composition of two permutation maps $\boldsymbol{\pi}$ and $\boldsymbol{\pi}'$. Throughout this proof, we use the notation $T(\mathcal{X}^{\boldsymbol{\pi}}) = T_{\mathrm{CI}}^{\boldsymbol{\pi}}$ and $T(\mathcal{X}) = T_{\mathrm{CI}}$. Let $T^{(1)}(\mathcal{X}) \leq \ldots \leq T^{(K_\ast)}(\mathcal{X})$ be the ordered statistics of $T(\mathcal{X}^{\boldsymbol{\pi}_1}),\ldots,T(\mathcal{X}^{\boldsymbol{\pi}_{K_\ast}})$ where $\{\boldsymbol{\pi}_1,\ldots,\boldsymbol{\pi}_{K_\ast}\} = \boldsymbol{\Pi}_{\mathrm{cyclic}}$. Note that the set $\{T(\mathcal{X}^{\boldsymbol{\pi}_1}),\ldots,T(\mathcal{X}^{\boldsymbol{\pi}_{K_\ast}}) \}$ has the same components as $\{ T(\mathcal{X}^{\boldsymbol{\pi}_1 \circ \boldsymbol{\pi} }),\ldots,T(\mathcal{X}^{\boldsymbol{\pi}_{K_\ast} \circ \boldsymbol{\pi}})\}$ for any $\boldsymbol{\pi} \in \boldsymbol{\Pi}_{\mathrm{cyclic}}$. Therefore, $T^{(i)}(\mathcal{X}) = T^{(i)}(\mathcal{X}^{\boldsymbol{\pi}})$ for any $i=1,\ldots,K_\ast$ and $\boldsymbol{\pi} \in \boldsymbol{\Pi}_{\mathrm{cyclic}}$. 

Having the above preliminaries in place, suppose we reject the null when $T(\mathcal{X}) > T^{(k)}(\mathcal{X})$ where $k = \ceil{(1-\alpha)K_\ast}$. Then, by letting $\boldsymbol{\pi}_0$ have the uniform distribution on $\boldsymbol{\Pi}_{\mathrm{cyclic}}$, we have
\begin{align*}
	\mP_{\widetilde{Q}_{X,Y,\widetilde{Z}}^n}\big\{ T(\mathcal{X}) > T^{(k)}(\mathcal{X}) \big\} ~\overset{(\mathrm{i})}{=}~ & \mP_{\widetilde{Q}_{X,Y,\widetilde{Z}}^n}\big\{ T(\mathcal{X}^{\boldsymbol{\pi}_0}) > T^{(k)}(\mathcal{X}^{\boldsymbol{\pi}_0}) \big\}  \\[.5em]
	\overset{(\mathrm{ii})}{=}~ & \mP_{\widetilde{Q}_{X,Y,\widetilde{Z}}^n}\big\{ T(\mathcal{X}^{\boldsymbol{\pi}_0}) > T^{(k)}(\mathcal{X}) \big\},
\end{align*}
where step~(i) uses the fact that $\mathcal{X}$ and $\mathcal{X}^{\boldsymbol{\pi}}$ have the same distribution for all $\boldsymbol{\pi} \in \boldsymbol{\Pi}_{\mathrm{cyclic}}$ under the law $\widetilde{Q}_{X,Y,\widetilde{Z}}^n$, and step~(ii) holds due to $T^{(k)}(\mathcal{X}) = T^{(k)}(\mathcal{X}^{\boldsymbol{\pi}})$. Now since $\boldsymbol{\pi}_0$ is uniformly distributed on $\boldsymbol{\Pi}_{\mathrm{cyclic}}$, we have
\begin{align*}
	\mP_{\widetilde{Q}_{X,Y,\widetilde{Z}}^n}\big\{ T(\mathcal{X}^{\boldsymbol{\pi}_0}) > T^{(k)}(\mathcal{X}) \big\} = \mE_{\widetilde{Q}_{X,Y,\widetilde{Z}}^n}\Bigg[ \frac{1}{K_\ast} \sum_{\boldsymbol{\pi}_i \in \boldsymbol{\Pi}_{\mathrm{cyclic}}} \mathds{1}\big\{ T(\mathcal{X}^{\boldsymbol{\pi}_{i}}) > T^{(k)}(\mathcal{X}) \big\}  \Bigg] \leq \alpha,
\end{align*}
where the last inequality holds since
\begin{align*}
	 \frac{1}{K_\ast} \sum_{\boldsymbol{\pi}_i \in \boldsymbol{\Pi}_{\mathrm{cyclic}}} \mathds{1}\big\{ T(\mathcal{X}^{\boldsymbol{\pi}_{i}}) > T^{(k)}(\mathcal{X}) \big\}  \leq \alpha.
\end{align*}
Finally, by Lemma~\ref{Lemma: quantile}, the permutation test from Algorithm~\ref{alg:double-binning permutation procedure} can be equivalently written as $\phi_{\mathrm{perm},n} = \mathds{1}\{T(\mathcal{X}) > T^{(k)}(\mathcal{X})\}$. Therefore, we have 
\begin{align*}
	\mE_{\widetilde{Q}_{X,Y,\widetilde{Z}}^n} [\phi_{\mathrm{perm},n}] \leq \alpha.
\end{align*}
Having this inequality, the other steps are exactly the same as the proofs of Theorem~\ref{Theorem: Validity of the permutation test under Hellinger smooth} and Theorem~\ref{Theorem: Validity of the permutation test under Renyi smooth}. This completes the proof of Proposition~\ref{Proposition: Type I error of the double-binned test}.

\subsection{Proof of Theorem~\ref{Theorem: Type II error of the double-binned test}}
The main idea of the proof of this result is similar to that of Theorem~\ref{Theorem: Type II error under continuous setting I}. However, it requires a much more involved analysis for the mean and the variance due to a sophisticated dependence structure among local permutations over smaller bins. For the single-binning method, there is a symmetric structure of the U-statistic under the (global) permutation, which allows us to show that the mean of the permuted U-statistic is zero and its variance is bounded by some tractable quantities~\citep[e.g.~Theorem 5.1 of][]{kim2020minimax}. However, this is no longer the case for the double-binning method, which leads to a technical challenge. 

 Throughout this proof, we only deal with the type II error (or $1-$power) of $\phi_{\mathrm{perm},N}$ since the result of $\phi_{\mathrm{perm},N}^\dagger$ follows directly by the inequality~(\ref{Eq: bound on type II error}). For a sufficiently large $\zeta>0$ (depending on $L$), \cite{neykov2020minimax} prove that the test
$\phi_{\mathrm{NBW},1}' = \mathds{1}(T_{\mathrm{CI}} \geq \zeta n^{1/5})$ has a small type II error in the sense of (\ref{Eq: type II error 1}) (indeed, without $e^{-n/8}$ term) when $\varepsilon \geq cn^{-2/5}$ for a sufficiently large $c$ depending on ($\zeta, L, \ell_1, \ell_2$). Moreover, since $\mathcal{P}'_{1,[0,1],\mathrm{TV}}(L) \subseteq \mathcal{P}_{1,[0,1],\mathrm{TV}}(L)$, it holds that
\begin{align}
	\sup_{P_{X,Y,Z} \in \mathcal{P}'_{1,[0,1],\mathrm{TV}}(L): \inf_{Q \in \mathcal{P}_{0,[0,1]}} \mathcal{D}_{\mathrm{TV}}(P_{X,Y,Z},Q) \geq \varepsilon} \mP_{P_{X,Y,Z}^N,N} [T_{\mathrm{CI}} \geq \zeta n^{1/5}] \leq \frac{1}{200}.
\end{align}
Now let $q_{1-\alpha}$ be the $1-\alpha$ quantile of the empirical distribution of $T_{\mathrm{CI}}^{\boldsymbol{\pi}_1},\ldots,T_{\mathrm{CI}}^{\boldsymbol{\pi}_{K_\ast}}$ where $\boldsymbol{\pi}_i \in \boldsymbol{\Pi}_{\mathrm{cyclic}}$. Our main goal is to show that $q_{1-\alpha} \leq \zeta' n^{1/5}$ for a sufficiently large $\zeta'$ with high probability (say $1/200$) under the given conditions. Since the type II error of $\phi_{\mathrm{perm},N}$ can be written as $\mP(T_{\mathrm{CI}} \leq q_{1-\alpha})$ in view of Lemma~\ref{Lemma: quantile}, this bound for the quantile together with the union bound implies that the type II error of $\phi_{\mathrm{perm},N}$ is smaller than $1/100$ by taking $\zeta$ sufficiently large. 

We break down the proof into several pieces for readability. 
\begin{itemize}
	\item \textbf{Step 1. Bounding the quantile.} In the first step, we prove that the $1-\alpha$ quantile is upper bounded by 
	\begin{align*}
		q_{1-\alpha} ~\lesssim~ \mE[T_{\mathrm{CI}}^{\boldsymbol{\pi}}]  + \sqrt{\mathrm{Var}[T_{\mathrm{CI}}^{\boldsymbol{\pi}}]},
	\end{align*}
	with high probability where $\boldsymbol{\pi} \sim F_{\boldsymbol{\pi}}$ defined in (\ref{Eq: distribution function}) below. Here we note that the expectation and the variance are with respect to both randomness from $\boldsymbol{\pi}$ and samples. 
	\item \textbf{Step 2. Bounding the expectation.} In the second step, we upper bound the expectation of $T_{\mathrm{CI}}^{\boldsymbol{\pi}}$ where $\boldsymbol{\pi} \sim F_{\boldsymbol{\pi}}$ in (\ref{Eq: distribution function}) and prove that
	\begin{align*}
		\mE[T_{\mathrm{CI}}^{\boldsymbol{\pi}}] ~=~ O\Bigg( \frac{n}{M^2} + 1 + \frac{M^2}{n} \Bigg),
	\end{align*}
	\item \textbf{Step 3. Bounding the variance.} In the third step, we upper bound the variance of $T_{\mathrm{CI}}^{\boldsymbol{\pi}}$ where $\boldsymbol{\pi} \sim F_{\boldsymbol{\pi}}$ in (\ref{Eq: distribution function}) and prove that
	\begin{align*}
		\mathrm{Var}[T_{\mathrm{CI}}^{\boldsymbol{\pi}}] ~=~ O\left( \frac{n^2}{M^4} + \frac{n}{M^2} + \frac{M^2}{n} + M \right).
	\end{align*}
	\item \textbf{Step 4. Completing the proof.} Lastly, we combine the previous results and prove the desired type II error bound for $\phi_{\mathrm{perm},N}$. 
\end{itemize}

\vskip 1em

\noindent \textbf{Step 1. Bounding the quantile.} We start by bounding the $1-\alpha$ quantile
\begin{align*}
	q_{1-\alpha} = \inf \bigg\{x \in \mathbb{R}: 1-\alpha \leq  \frac{1}{K_\ast} \sum_{\boldsymbol{\pi}_i \in \boldsymbol{\Pi}_{\mathrm{cyclic}}} \mathds{1} \big(T_{\mathrm{CI}}^{\boldsymbol{\pi}_i} \leq x \big) \bigg\}. 
\end{align*}
Consider a subset of $\boldsymbol{\Pi}_{\mathrm{cyclic}}$ where we remove the identity permutation within each small bin. Let us denote the resulting subset by $ \boldsymbol{\Pi}_{\mathrm{cyclic}}^\mathrm{-id}$ whose cardinality is $K_\ast^{-\mathrm{id}} := \prod_{(m,k) \in [M] \times [b]} \max\{\sigma_{m,k}-1,1\}$. By letting $\alpha_\ast$ such that $1 - \alpha_\ast =  \min\big\{(1-\alpha) K_\ast / K_\ast^\mathrm{-id},1 \big\}$, define 
\begin{align*}
	q_{1-\alpha_\ast}^{\mathrm{-id}} ~=~ & \inf \bigg\{x \in \mathbb{R}: \min\bigg\{(1-\alpha), \frac{K_\ast^{\mathrm{-id}}}{K_\ast}\bigg\} \leq  \frac{1}{K_\ast} \sum_{\boldsymbol{\pi}_i \in \boldsymbol{\Pi}_{\mathrm{cyclic}}^\mathrm{-id}} \mathds{1} \big( T_{\mathrm{CI}}^{\boldsymbol{\pi}_i} \leq x \big) \bigg\}  \\[.5em]
	= ~ & \inf \bigg\{x \in \mathbb{R}: \min\bigg\{ (1-\alpha) \frac{K_\ast}{K_\ast^\mathrm{-id}}, 1 \bigg\} \leq  \frac{1}{K_\ast^\mathrm{-id}} \sum_{\boldsymbol{\pi}_i \in \boldsymbol{\Pi}_{\mathrm{cyclic}}^\mathrm{-id}} \mathds{1} \big( T_{\mathrm{CI}}^{\boldsymbol{\pi}_i} \leq x \big) \bigg\}.
\end{align*}
Since $\sum_{\boldsymbol{\pi}_i \in \boldsymbol{\Pi}_{\mathrm{cyclic}}} \mathds{1} \big(T_{\mathrm{CI}}^{\boldsymbol{\pi}_i} \leq x \big)$ is a sum of non-negative numbers, we have $q_{1-\alpha} \leq q_{1-\alpha_\ast}^{\mathrm{-id}}$. Therefore
\begin{align*}
	\mP(T_{\mathrm{CI}} \leq q_{1-\alpha}) ~\leq~ \mP(T_{\mathrm{CI}} \leq q_{1-\alpha_\ast}^{\mathrm{-id}}).
\end{align*} 
Let us define an event $\mathcal{A}$ as
\begin{align} \label{Eq: event A}
	\mathcal{A} = \bigg\{\frac{K_\ast}{K_\ast^\mathrm{-id}} \leq \frac{1 - \alpha/2}{1-\alpha}  \bigg\},
\end{align}
and let $q_{1-\alpha/2}^{\mathrm{-id}}$ be the $1-\alpha/2$ quantile of the empirical distribution 
\begin{align} \label{Eq: distribution function}
	F_{\boldsymbol{\pi}}(x):=\frac{1}{K_\ast^\mathrm{-id}} \sum_{\boldsymbol{\pi}_i \in \boldsymbol{\Pi}_{\mathrm{cyclic}}^\mathrm{-id}} \mathds{1} \big( T_{\mathrm{CI}}^{\boldsymbol{\pi}_i} \leq x \big).
\end{align}
Then, since $q_{1-\alpha_\ast}^{\mathrm{-id}} \leq q_{1-\alpha/2}^{\mathrm{-id}}$ under $\mathcal{A}$, we have 
\begin{align*}
	\mP(T_{\mathrm{CI}} \leq q_{1-\alpha_\ast}^{\mathrm{-id}}) ~\leq~ \mP(T_{\mathrm{CI}} \leq q_{1-\alpha/2}^{\mathrm{-id}}) + \mP(\mathcal{A}^c).
\end{align*}
In what follows, we use the notation $\mE_{\boldsymbol{\pi}}[\cdot]$ and $\mathrm{Var}_{\boldsymbol{\pi}}[\cdot]$ to denote the expectation and variance operators with respect to $F_{\boldsymbol{\pi}}$ and let $\mP_{\boldsymbol{\pi}}[\cdot] = \mE_{\boldsymbol{\pi}}[\mathds{1}(\cdot)]$. Given this notation, we further upper bound $q_{1-\alpha/2}^{\mathrm{-id}}$. To this end, we follow the proof of Lemma 3.1 of \cite{kim2020minimax}. First, for $t>0$, by Chebyshev's inequality,
\begin{align*}
	\mP_{\boldsymbol{\pi}}(T_{\mathrm{CI}}^{\boldsymbol{\pi}} - \mE_{\boldsymbol{\pi}}[T_{\mathrm{CI}}^{\boldsymbol{\pi}}] \geq t)  \leq~  \frac{1}{t^2} \mathrm{Var}_{\boldsymbol{\pi}}[T_{\mathrm{CI}}^{\boldsymbol{\pi}}]. 
\end{align*}
By setting the right-hand side to be $\alpha/2$, we see that $q_{1-\alpha/2}^{\mathrm{-id}}$ is deterministically bounded by 
\begin{align*}
	q_{1-\alpha/2}^{\mathrm{-id}} \leq \mE_{\boldsymbol{\pi}}[T_{\mathrm{CI}}^{\boldsymbol{\pi}}] + \sqrt{\frac{2}{\alpha}\mathrm{Var}_{\boldsymbol{\pi}}[T_{\mathrm{CI}}^{\boldsymbol{\pi}}]}.
\end{align*}
Next applying Chebyshev's inequality and Markov's inequality yields
\begin{align*}
	& \mP\bigg( \mE_{\boldsymbol{\pi}}[T_{\mathrm{CI}}^{\boldsymbol{\pi}}] \geq \mE[T_{\mathrm{CI}}^{\boldsymbol{\pi}}] + \sqrt{\frac{2}{\beta} \mathrm{Var}\big\{ \mE_{\boldsymbol{\pi}}[T_{\mathrm{CI}}^{\boldsymbol{\pi}}]\big\}} \bigg) \leq \frac{\beta}{2} \quad \text{and} \\[.5em]
	& \mP \bigg( \frac{2}{\alpha}\mathrm{Var}_{\boldsymbol{\pi}}[T_{\mathrm{CI}}^{\boldsymbol{\pi}}] \geq \frac{4}{\alpha \beta}\mE\big\{\mathrm{Var}_{\boldsymbol{\pi}}[T_{\mathrm{CI}}^{\boldsymbol{\pi}}]\big\}  \bigg) \leq \frac{\beta}{2}.
\end{align*}
From the above observation together with the union bound, it holds with probability at least $1-\beta$ that 
\begin{align*}
	q_{1-\alpha/2}^{\mathrm{-id}} ~\leq~ & \mE[T_{\mathrm{CI}}^{\boldsymbol{\pi}}]  + \sqrt{\frac{2}{\beta} \mathrm{Var}\big\{ \mE_{\boldsymbol{\pi}}[T_{\mathrm{CI}}^{\boldsymbol{\pi}}]\big\}} + \sqrt{\frac{4}{\alpha \beta}\mE\big\{\mathrm{Var}_{\boldsymbol{\pi}}[T_{\mathrm{CI}}^{\boldsymbol{\pi}}]\big\}} \\[.5em]
	\leq ~ &   \mE[T_{\mathrm{CI}}^{\boldsymbol{\pi}}]  + \sqrt{\frac{32}{\alpha \beta}\mathrm{Var}[T_{\mathrm{CI}}^{\boldsymbol{\pi}}]},
\end{align*}
where the last inequality uses $\mathrm{Var}\big\{ \mE_{\boldsymbol{\pi}}[T_{\mathrm{CI}}^{\boldsymbol{\pi}}]\big\} + \mE\big\{\mathrm{Var}_{\boldsymbol{\pi}}[T_{\mathrm{CI}}^{\boldsymbol{\pi}}]\big\} \leq \mathrm{Var}[T_{\mathrm{CI}}^{\boldsymbol{\pi}}]$. In summary, by taking $\alpha = 1/100$ and $\beta = 1/200$, we have obtained that
\begin{align} \label{Eq: summary}
	\mP(T_{\mathrm{CI}} \leq q_{1-\alpha}) ~\leq~ \mP\Big(T_{\mathrm{CI}} \lesssim \mE[T_{\mathrm{CI}}^{\boldsymbol{\pi}}] + \sqrt{\mathrm{Var}[T_{\mathrm{CI}}^{\boldsymbol{\pi}}]}\Big) + \frac{1}{200} + \mP(\mathcal{A}^c).
\end{align}

\vskip 1em

\noindent \textbf{Step 2. Upper bounding the expectation.} Recall that 
\begin{align*}
	T_{\mathrm{CI}}^{\boldsymbol{\pi}} ~=~&  \sum_{m \in [M]} \mathds{1}(\sigma_m \geq 4) \sigma_m U(\bW_m^{\boldsymbol{\pi}}),
\end{align*}
where $\boldsymbol{\pi}$ is uniformly distributed on $\boldsymbol{\Pi}_{\mathrm{cyclic}}^\mathrm{-id}$. When $\sigma_m < 4$, $\mathds{1}(\sigma_m \geq 4) \sigma_m U(\bW_m^{\boldsymbol{\pi}})$ becomes zero, which does not contribute to $T_{\mathrm{CI}}^{\boldsymbol{\pi}}$. Hence, we assume $\sigma_m \geq 4$ and 	show that 
\begin{align} \label{Eq: bounding the expectation}
	\mE[U(\bW_m^{\boldsymbol{\pi}})|\sigma_m] = O\left( \frac{1}{M^2} + \frac{1}{M\sigma_m} + \frac{1}{\sigma_m^2} \right) \quad \text{for each $m \in [M]$.}
\end{align}
Recall $U(\bW_m^{\boldsymbol{\pi}})$ is the U-statistic with the kernel that takes the form
\begin{align*}
	\varphi_{i,j,k,\ell}^{\boldsymbol{\pi}}:=\sum_{x \in [\ell_1], y \in [\ell_2]} & \big\{ \mathds{1}(X_{i,m}=x,Y_{{\pi}_{i,m}}=y) -  \mathds{1}(X_{i,m}=x) \mathds{1}(Y_{{\pi}_{j,m}}=y) \big\} \times \\[.5em]
	&	\big\{ \mathds{1}(X_{k,m}=x,Y_{{\pi}_{k,m}}=y) -  \mathds{1}(X_{k,m}=x) \mathds{1}(Y_{{\pi}_{\ell,m}}=y) \big\},
\end{align*}
where $\{i,j,k,\ell\}$ are distinct integers between $1$ and $\sigma_m$, and $\{\pi_{i,m}, \pi_{j,m}, \pi_{k,m}, \pi_{\ell,m}\}$ are components of $\boldsymbol{\pi}$ that are used to permute the corresponding components of $\{Y_{i,m},Y_{j,m},Y_{k,m},Y_{\ell,m}\}$. Due to a complicated local structure of $\boldsymbol{\pi}$, it is not trivial to compute the expectation of this kernel. This is in contrast to the single-binning case where this expectation of the kernel is zero under the local permutation law. We instead upper bound the expectation using the following observation. Suppose that there is no collision, meaning that neither of $Y_{i,m}, Y_{k,m}$ corresponds to $Y_{\pi_{i,m}}, Y_{\pi_{j,m}},Y_{\pi_{k,m}}, Y_{\pi_{\ell,m}}$. In this case, $\{X_{i,m},X_{j,m}, Y_{\pi_{i,m}}, Y_{\pi_{j,m}},Y_{\pi_{k,m}}, Y_{\pi_{\ell,m}}\}$ are mutually independent, and the expectation of the kernel is bounded above by
\begin{equation}
\begin{aligned} \label{Eq: bounding the varphi}
		\big|\mE[\varphi_{i,j,k,\ell}^{\boldsymbol{\pi}}|\sigma_m]\big| ~\leq~ \Bigg| \max_{k_1,k_2,k_3,k_4 \in [b]} \Bigg[\sum_{x \in [\ell_1], y \in [\ell_2]} & \{p_{m,k_1}(x)q_{m,k_1}(y) - p_{m,k_1}(x)q_{m,k_2}(y) \} \times \\[.5em]
		& \{p_{m,k_3}(x)q_{m,k_3}(y) - p_{m,k_3}(x)q_{m,k_4}(y) \}\Bigg]\Bigg|,
\end{aligned}
\end{equation}
where $p_{m,k}(x):= \mP(X=x|Z \in B_{m,k})$ and $q_{m,k}(y) :=\mP(Y=y|Z \in B_{m,k})$. Furthermore, since we assume $P_{X,Y,Z} \in \mathcal{P}'_{1,[0,1],\mathrm{TV}(L)}$,
\begin{align*}
	 & \mathcal{D}_{\mathrm{TV}}(q_{m,k_1}, q_{m,k_2}) \\[.5em]
	  = ~ &  \frac{1}{2} \sum_{y \in [\ell_2]} \bigg| \int_{B_{m,k_1}}p_{Y|Z}(y|z)  d\widetilde{P}_{Z|Z \in B_{m,k_1}}(z) - \int_{B_{m,k_2}}p_{Y|Z}(y|z')  d\widetilde{P}_{Z|Z \in B_{m,k_2}}(z')\bigg| \\[.5em]
	 \leq ~ & \frac{1}{2} \int_{B_{m,k_1}} \int_{B_{m,k_2}}\sum_{y \in [\ell_2]} \big| p_{Y|Z}(y|z)  - p_{Y|Z}(y|z')  \big| d\widetilde{P}_{Z|Z \in B_{m,k_1}}(z)  d\widetilde{P}_{Z|Z \in B_{m,k_2}}(z') \\[.5em]
	 \leq ~ & \frac{1}{2} \max_{z \in B_{m,k_1},z' \in B_{m,k_2}}\delta(z,z') \leq \frac{C_1}{M},
\end{align*}
for some constant $C_1>0$. Therefore $\big|\mE[\varphi_{i,j,k,\ell}^{\boldsymbol{\pi}}]\big|$ is further bounded by
\begin{align*}
	& \big|\mE[\varphi_{i,j,k,\ell}^{\boldsymbol{\pi}}|\sigma_m]\big| \\[.5em]
	\leq ~ & \Bigg| \max_{k_1,k_2,k_3,k_4 \in [b]} \Bigg[ \sum_{x \in [\ell_1]} p_{m,k_1}(x)p_{m,k_3}(x) \sum_{y \in [\ell_2]} |q_{m,k_1}(y) - q_{m,k_2}(y)| \sum_{y \in [\ell_2]} |q_{m,k_3}(y) - q_{m,k_4}(y)|  \Bigg] \Bigg| \\[.5em]
	\leq ~ & 4 \max_{k_1,k_2,k_3,k_4 \in [b]} \mathcal{D}_{\mathrm{TV}}(q_{m,k_1}, q_{m,k_2}) \times \mathcal{D}_{\mathrm{TV}}(q_{m,k_3}, q_{m,k_4}) \\[.5em]
	\leq ~ & \frac{C_2}{M^2}.
\end{align*}
Suppose that there is only one collision, i.e.~there exists only one dependent pair of $(X,Y)$ among $\{X_{i,m},X_{j,m}, Y_{\pi_{i,m}}, Y_{\pi_{j,m}},Y_{\pi_{k,m}}, Y_{\pi_{\ell,m}}\}$. Since we deal with a cyclic permutation excluding the identity map, there are only a few possibilities of one collision. In particular, a cyclic permutation, which is not the identity map, is a derangement~(i.e.~no point will be at its original place). This guarantees that $X_{i,m}$ and $Y_{\pi_{i,m}}$ are independent for each $i$ and $m$, which simplifies our analysis.

\vskip 1em

\noindent \textbf{Case 1)} Suppose that $X_{i,m}$ and $Y_{\pi_{j,m}}$ collide. In this case, $\{X_{k,m}, Y_{\pi_{k,m}}, Y_{\pi_{\ell,m}}\}$ are still independent, which leads to the bound of the expectation 
\begin{align*}
	\big|\mE[\varphi_{i,j,k,\ell}^{\boldsymbol{\pi}} |\sigma_m]\big| ~:=~ & \bigg| \sum_{x \in [\ell_1], y \in [\ell_2]} \mE\big\{ \mathds{1}(X_{i,m}=x,Y_{{\pi}_{i,m}}=y) -  \mathds{1}(X_{i,m}=x) \mathds{1}(Y_{{\pi}_{j,m}}=y) \big\} \times \\[.5em]
	&	~~~~~~~~~~~~~~ \mE\big\{ \mathds{1}(X_{k,m}=x,Y_{{\pi}_{k,m}}=y) -  \mathds{1}(X_{k,m}=x) \mathds{1}(Y_{{\pi}_{\ell,m}}=y) \big\} \bigg| \\[.5em]
	\leq ~ &  2 \sum_{x \in [\ell_1], y \in [\ell_2]} \big| \mE\big\{ \mathds{1}(X_{k,m}=x,Y_{{\pi}_{k,m}}=y) -  \mathds{1}(X_{k,m}=x) \mathds{1}(Y_{{\pi}_{\ell,m}}=y) \big\} \big| \\[.5em]
	\leq~ & 2  \sup_{k_1,k_2 \in [b]} \sum_{x \in [\ell_1], y \in [\ell_2]} p_{m,k_1}(x) \times \big| q_{m,k_1}(y) - q_{m,k_2}(y)\big| \\[.5em]
	\leq ~ & \frac{C_3}{M},
\end{align*}
where the last inequality follows similarly as the previous analysis. 

\vskip 1em

\noindent \textbf{Case 2)} Next, suppose that $X_{i,m}$ and $Y_{\pi_{k,m}}$ collide. In this case, $Y_{\pi_{i,m}}$ and $Y_{\pi_{j,m}}$ are independent of the others. Using this property, the expectation is bounded by 
\begin{align*}
	\big|\mE[\varphi_{i,j,k,\ell}^{\boldsymbol{\pi}}|\sigma_m]\big| ~:=~ & \bigg| \sum_{x \in [\ell_1], y \in [\ell_2]} \mE \big\{\mathds{1}(Y_{{\pi}_{i,m}}=y) -  \mathds{1}(Y_{{\pi}_{j,m}}=y) \big\} \mE[\mathds{1}(X_{k,m}=x)] \times \\[.5em]
	&	~~~~~~~~~~~~~~ \mE\big[\mathds{1}(X_{i,m}=x)\big\{ \mathds{1}(Y_{{\pi}_{k,m}}=y) -   \mathds{1}(Y_{{\pi}_{\ell,m}}=y) \big\} \big] \bigg| \\[.5em]
	\leq~ & 2  \sup_{k_1,k_2 \in [b]} \sum_{x \in [\ell_1], y \in [\ell_2]} p_{m,k_1}(x) \times \big| q_{m,k_1}(y) - q_{m,k_2}(y)\big| \\[.5em]
	\leq ~ & \frac{C_3}{M}.
\end{align*}
The other cases with one collision can be bounded similarly by $O(M^{-1})$.

Now consider another representation of the U-statistic:
\begin{align*}
	U(\bW_m^{\boldsymbol{\pi}}) ~=~ \frac{(\sigma_m-4)!}{\sigma_m!} \sum_{\{i,j,k,\ell\} \in \boldsymbol{i}_4^{\sigma_m}} \varphi_{i,j,k,\ell}^{\boldsymbol{\pi}},
\end{align*}
where $\boldsymbol{i}_4^{\sigma_m}$ denotes the set of all $4$-tuples drawn without replacement from $[\sigma_m]$. Next we decompose the summation into three cases depending on whether there exists a collision among 
\begin{align} \label{Eq: 4-tuples}
	\{(X_{i,m},Y_{\pi_{i,m}}),(X_{j,m},Y_{\pi_{j,m}}),(X_{k,m},Y_{\pi_{k,m}}),(X_{\ell,m},Y_{\pi_{\ell,m}})\}.
\end{align}
In particular,
\begin{align*}
	\sum_{\{i,j,k,\ell\} \in \boldsymbol{i}_4^{\sigma_m}}\varphi_{i,j,k,\ell}^{\boldsymbol{\pi}}  ~=~ &  \sum_{\{i,j,k,\ell\} \in \boldsymbol{i}_4^{\sigma_m} \cap A_{\mathrm{collision}}^0} \varphi_{i,j,k,\ell}^{\boldsymbol{\pi}} + \sum_{\{i,j,k,\ell\} \in \boldsymbol{i}_4^{\sigma_m} \cap A_{\mathrm{collision}}^1}\varphi_{i,j,k,\ell}^{\boldsymbol{\pi}} \\[.5em]
	& + \sum_{\{i,j,k,\ell\} \in \boldsymbol{i}_4^{\sigma_m} \cap A_{\mathrm{collision}}^{1+}}\varphi_{i,j,k,\ell}^{\boldsymbol{\pi}},
\end{align*}
where 1) $\boldsymbol{i}_4^{\sigma_m} \cap A_{\mathrm{collision}}^0$ is the subset of $\boldsymbol{i}_4^{\sigma_m}$ such that there exists no collision, 2) $\boldsymbol{i}_4^{\sigma_m} \cap A_{\mathrm{collision}}^1$ is the subset of $\boldsymbol{i}_4^{\sigma_m}$ such that there exists only one collision, and 3) $\boldsymbol{i}_4^{\sigma_m} \cap A_{\mathrm{collision}}^{1+}$ is the subset of $\boldsymbol{i}_4^{\sigma_m}$ such that there exists more than one collision. 

Now we upper bound each summation. First of all, by the previous analysis,
\begin{align*}
	\frac{(\sigma_m-4)!}{\sigma_m!}  \sum_{\{i,j,k,\ell\} \in \boldsymbol{i}_4^{\sigma_m} \cap A^0_{\mathrm{collision}}} \mE[\varphi_{i,j,k,\ell}^{\boldsymbol{\pi}}|\sigma_m]  ~=~ O\left( \frac{1}{M^2} \right).
\end{align*}
Next, since $\boldsymbol{\pi} \in \boldsymbol{\Pi}_{\mathrm{cyclic}}^\mathrm{-id}$ is a derangement by construction, if we pick four observations~(\ref{Eq: 4-tuples}) with no collision (i.e.~$\{i,j,k,\ell\}$ and $\{\pi_i,\pi_j,\pi_k,\pi_{\ell}\}$ have no shared component) out of $\sigma_m$ observations, there are at least $\sigma_m(\sigma_m-2)(\sigma_m-4)(\sigma_m-6)$ ways. Thus in order to have a collision, we have at most $O(\sigma_m^3)$ choices. Moreover, as analyzed before, the expectation of the kernel is bounded by $O(M^{-1})$ if there is one collision. Therefore
\begin{align*}
	\frac{(\sigma_m-4)!}{\sigma_m!}  \sum_{\{i,j,k,\ell\} \in \boldsymbol{i}_4^{\sigma_m} \cap A_{\mathrm{collision}}^1} \mE[\varphi_{i,j,k,\ell}^{\boldsymbol{\pi}}|\sigma_m]  ~=~  O\left( \frac{1}{M \sigma_m} \right).
\end{align*}
Lastly, there are at most $\asymp \sigma_m^2$ cases to have more than one collision. In this case, using the fact that the kernel $\varphi_{i,j,k,\ell}^{\boldsymbol{\pi}}$ is bounded between $-1$ and $1$, we have
\begin{align*}
	\frac{(\sigma_m-4)!}{\sigma_m!}  \sum_{\{i,j,k,\ell\} \in \boldsymbol{i}_4^{\sigma_m} \cap A_{\mathrm{collision}}^{1+}} \mE[\varphi_{i,j,k,\ell}^{\boldsymbol{\pi}}|\sigma_m]  ~=~  O\left( \frac{1}{\sigma_m^2} \right).
\end{align*}
These three bounds imply the claim~(\ref{Eq: bounding the expectation}). Therefore
\begin{align*}
	\mE[T_{\mathrm{CI}}^{\boldsymbol{\pi}}] ~=~  O\Bigg( \frac{n}{M^2} + 1 + \sum_{m \in [M]} \mE\big[\mathds{1}(\sigma_m \geq 4) \sigma_m^{-1}\big] \Bigg).
\end{align*}
By recalling $\sigma_m \sim \mathrm{Pois}(np_m)$, we have the following inequalities
\begin{align} \label{Eq: poisson reciprocal}
	\mE\bigg[\mathds{1}(\sigma_m \geq 4) \frac{1}{\sigma_m}\bigg] ~\leq~ \frac{5}{4}\mE\bigg[\mathds{1}(\sigma_m \geq 4) \frac{1}{\sigma_m + 1}\bigg] ~\leq~ \frac{5}{4} \mE\bigg[ \frac{1}{\sigma_m + 1}\bigg].
\end{align}
Furthermore,
\begin{align*}
	\mE\bigg[ \frac{1}{\sigma_m + 1}\bigg] ~\leq~ \frac{1 - e^{-np_m}}{np_m}  ~\leq~ \frac{M}{c_{\mathrm{low}}n},
\end{align*}
where the last inequality uses the condition $p_m \geq c_{\mathrm{low}}M^{-1}$ for all $m \in [M]$. Therefore, we have 
\begin{align*}
	\mE[T_{\mathrm{CI}}^{\boldsymbol{\pi}}] ~=~  O\Bigg( \frac{n}{M^2} + 1 + \frac{M^2}{n} \Bigg).
\end{align*}

\vskip 1em

\noindent \textbf{Step 3. Upper bounding the variance.}
Next we upper bound the variance of the U-statistic~$U(\bW_m^{\boldsymbol{\pi}})$. To compute the variance of $U(\bW_m^{\boldsymbol{\pi}})$, first set up some notation. Similarly as in Section~\ref{Section: Tightness}, let
\begin{align*}
	\psi_{ij}^{m,\boldsymbol{\pi}}(x,y) = \mathds{1}(X_{i,m}=x,Y_{\pi_{i,m}}=y) -  \mathds{1}(X_{i,m}=x) \mathds{1}(Y_{\pi_{j,m}}=y),
\end{align*}
and define a symmetrized kernel function as 
\begin{align*} 
	h_{i_1,i_2,i_3,i_4}^{m,\boldsymbol{\pi}}= \frac{1}{4!} \sum_{\pi \in \boldsymbol{\Pi}_4} \sum_{x \in [\ell_1], y \in [\ell_2]} \psi^{m,\boldsymbol{\pi}}_{\pi(1) \pi(2)}(x,y) \psi^{m,\boldsymbol{\pi}}_{\pi(3)\pi(4)}(x,y)
\end{align*}
where $\boldsymbol{\Pi}_4$ is the set of all permutations of $\{i_1,i_2,i_3,i_4\}$. 
Then the U-statistic can be written as 
\begin{align*}
	U(\bW_m^{\boldsymbol{\pi}}) ~=~ \frac{1}{\binom{\sigma_m}{4}} \sum_{i_1<i_2<i_3<i_4: (i_1,i_2,i_3,i_4) \in[\sigma_m]} h_{i_1,i_2,i_3,i_4}^{m,\boldsymbol{\pi}},
\end{align*}
and it variance has the form as 
\begin{align*}
	\mathrm{Var}[U(\bW_m^{\boldsymbol{\pi}})|\sigma_m] ~=~ \frac{1}{\binom{\sigma_m}{4}^2} \sum_{\substack{i_1<i_2<i_3<i_4 \\: (i_1,i_2,i_3,i_4) \in[\sigma_m]}} \sum_{\substack{j_1<j_2<j_3<j_4 \\: (j_1,j_2,j_3,j_4) \in[\sigma_m]}} \mathrm{Cov}(h_{i_1,i_2,i_3,i_4}^{m,\boldsymbol{\pi}}, h_{j_1,j_2,j_3,j_4}^{m,\boldsymbol{\pi}}).
\end{align*}
Here and hereafter, we often ignore the conditional event $\sigma_m$ when we write the expectation and the (co)variance operator for simplicity. 

Now we decompose the summation into a few different parts depending on the number of common indices between $\{i_1,i_2,i_3,i_4\}$ and $\{j_1,j_2,j_3,j_4\}$. In particular, let $A_{\mathrm{com}}^{k}$ be the set of indices $\{i_1,i_2,i_3,i_4,j_1,j_2,j_3,j_4\} \in [\sigma_m]$ such that $i_1<i_2<i_3<i_4$, $j_1<j_2<j_3<j_4$ and there are $k$ elements in common between $\{i_1,i_2,i_3,i_4\}$ and $\{j_1,j_2,j_3,j_4\}$. Using this notation, the variance can be written as 
\begin{align*}
	& \mathrm{Var}[U(\bW_m^{\boldsymbol{\pi}})|\sigma_m] \\[.5em]
	 ~=~ & \frac{1}{\binom{\sigma_m}{4}^2} \sum_{A_{\mathrm{com}}^{0}} \mathrm{Cov}(h_{i_1,i_2,i_3,i_4}^{m,\boldsymbol{\pi}}, h_{j_1,j_2,j_3,j_4}^{m,\boldsymbol{\pi}}) + \frac{1}{\binom{\sigma_m}{4}^2} \sum_{A_{\mathrm{com}}^{1}} \mathrm{Cov}(h_{i_1,i_2,i_3,i_4}^{m,\boldsymbol{\pi}}, h_{j_1,j_2,j_3,j_4}^{m,\boldsymbol{\pi}}) \\[.5em]
	 + ~ &  \frac{1}{\binom{\sigma_m}{4}^2} \sum_{A_{\mathrm{com}}^{1+}} \mathrm{Cov}(h_{i_1,i_2,i_3,i_4}^{m,\boldsymbol{\pi}}, h_{j_1,j_2,j_3,j_4}^{m,\boldsymbol{\pi}}).
\end{align*}
Since 1) the number of cases in which $\{i_1,i_2,i_3,i_4\}$ and $\{j_1,j_2,j_3,j_4\}$ have more than one common index is at most $\asymp \sigma_m^6$ and 2) that the kernel is bounded, the last term in the above display is bounded by 
\begin{align*}
	 \frac{1}{\binom{\sigma_m}{4}^2} \sum_{A_{\mathrm{com}}^{1+}} \mathrm{Cov}(h_{i_1,i_2,i_3,i_4}^{m,\boldsymbol{\pi}}, h_{j_1,j_2,j_3,j_4}^{m,\boldsymbol{\pi}}) ~=~ O\left(\sigma_m^{-2}\right).
\end{align*}
Next consider the cases in which $\{i_1,i_2,i_3,i_4\}$ and $\{j_1,j_2,j_3,j_4\}$ have the exactly one common index. Without loss of generality, we let $i_1 = j_1$ and consider the resulting pairs
\begin{align*}
	 \{ & (X_{i_1,m},Y_{\pi_{i_1,m}}),(X_{i_2,m},Y_{\pi_{i_2,m}}),(X_{i_3,m},Y_{\pi_{i_3,m}}),(X_{i_4,m},Y_{\pi_{i_4,m}}), \\[.5em]
	&(X_{j_2,m},Y_{\pi_{j_2,m}}),(X_{j_3,m},Y_{\pi_{j_3,m}}),(X_{j_4,m},Y_{\pi_{j_4,m}})\}.
\end{align*}
For now, we assume that there is no collision among the above set. Then we can treat them as if they are i.i.d.~observations with $X \independent Y$. Conditional on this assumption of no collision, we can follow the proof of Lemma 5.1 of \cite{neykov2020minimax} and show
\begin{align} \label{Eq: goal (1 common index)}
	\mathrm{Cov}(h_{i_1,i_2,i_3,i_4}^{m,\boldsymbol{\pi}}, h_{j_1,j_2,j_3,j_4}^{m,\boldsymbol{\pi}}) ~=~ O(M^{-2}).
\end{align}
To prove this, first let $\boldsymbol{\tau} := \{\tau(1),\tau(2),\tau(3),\tau(4)\}$ be sampled uniformly from the set of all permutations of $\{i_1,i_2,i_3,i_4\}$. Similarly let $\boldsymbol{\tau}' := \{\tau'(1),\tau'(2),\tau'(3),\tau'(4)\}$ be sampled uniformly from the set of all permutations of $\{j_1,j_2,j_3,j_4\}$. Then by treating $\boldsymbol{\tau}$ and $\boldsymbol{\tau}'$ to be independent random vectors, we can re-write the expectation as
\begin{align*}
	& \mE [h_{i_1,i_2,i_3,i_4}^{m,\boldsymbol{\pi}} h_{j_1,j_2,j_3,j_4}^{m,\boldsymbol{\pi}}] \\[.5em]
	=~ & \mE \Bigg[\sum_{x \in [\ell_1], y \in [\ell_2]} \psi^{m,\boldsymbol{\pi}}_{\tau(1) \tau(2)}(x,y) \psi^{m,\boldsymbol{\pi}}_{\tau(3)\tau(4)}(x,y) \sum_{x' \in [\ell_1], y' \in [\ell_2]} \psi^{m,\boldsymbol{\pi}}_{\tau'(1) \tau'(2)}(x',y') \psi^{m,\boldsymbol{\pi}}_{\tau'(3)\tau'(4)}(x',y') \Bigg].
\end{align*}
Using the condition that there exists exactly one common index between $\{\tau(1),\tau(2),\tau(3),\tau(4)\}$ and $\{\tau'(1),\tau'(2),\tau'(3),\tau'(4)\}$ (without loss of generality, let the shared common index be between $\{\tau(1),\tau(2)\}$ and $\{\tau'(1),\tau'(2)\}$), the above expectation becomes 
\begin{align*}
	  & \sum_{x \in [\ell_1], y \in [\ell_2]}  \sum_{x' \in [\ell_1], y' \in [\ell_2]} \mE\big[\psi^{m,\boldsymbol{\pi}}_{\tau(1) \tau(2)}(x,y)  \psi^{m,\boldsymbol{\pi}}_{\tau'(1) \tau'(2)}(x',y')\big] \mE\big[ \psi^{m,\boldsymbol{\pi}}_{\tau(3)\tau(4)}(x,y)\big]  \mE\big[\psi^{m,\boldsymbol{\pi}}_{\tau'(3)\tau'(4)}(x',y')\big] \\[.5em]
	  \leq ~ & \sqrt{\sum_{x \in [\ell_1], y \in [\ell_2]}  \sum_{x' \in [\ell_1], y' \in [\ell_2]} \Big\{ \mE\big[\psi^{m,\boldsymbol{\pi}}_{\tau(1) \tau(2)}(x,y)  \psi^{m,\boldsymbol{\pi}}_{\tau'(1) \tau'(2)}(x',y')\big]\Big\}^2 } \times \\[.5em]
	  &  \sqrt{\sum_{x \in [\ell_1], y \in [\ell_2]}  \sum_{x' \in [\ell_1], y' \in [\ell_2]} \Big\{ \mE\big[ \psi^{m,\boldsymbol{\pi}}_{\tau(3)\tau(4)}(x,y)\big]  \mE\big[\psi^{m,\boldsymbol{\pi}}_{\tau'(3)\tau'(4)}(x',y')\big] \Big\}^2 },
\end{align*}
where the second step follows by Cauchy--Schwarz inequality. Now, based on the notation used in (\ref{Eq: bounding the varphi})
\begin{align} \nonumber
	\sum_{x \in [\ell_1], y \in [\ell_2]} \Big\{ \mE\big[ \psi^{m,\boldsymbol{\pi}}_{\tau(3)\tau(4)}(x,y)\big] \Big\}^2 ~\leq~ & \max_{k_1,k_2 \in [b]} \sum_{x \in [\ell_1], y \in [\ell_2]} \Big\{ p_{m,k_1}(x) q_{m,k_1}(x) -  p_{m,k_1}(x) q_{m,k_2}(x) \Big\}^2 \\[.5em] \label{Eq: TV bound}
	\leq ~ &  \max_{k_1,k_2 \in [b]} \mathcal{D}_{\mathrm{TV}}^2(q_{m,k_1},q_{m,k_2})  ~\leq~  \frac{C_1}{M^2}.
\end{align}
This bound gives 
\begin{align} \label{Eq: expectation of the product}
	 \mE [h_{i_1,i_2,i_3,i_4}^{m,\boldsymbol{\pi}} h_{j_1,j_2,j_3,j_4}^{m,\boldsymbol{\pi}}]  ~=~  O\left( M^{-2} \right).
\end{align}
Next we look at the product of the expectations
\begin{align*}
	& \mE [h_{i_1,i_2,i_3,i_4}^{m,\boldsymbol{\pi}}] \times \mE[h_{j_1,j_2,j_3,j_4}^{m,\boldsymbol{\pi}}]  \\[.5em]
	= ~ & \mE \Bigg[\sum_{x \in [\ell_1], y \in [\ell_2]} \psi^{m,\boldsymbol{\pi}}_{\tau(1) \tau(2)}(x,y) \psi^{m,\boldsymbol{\pi}}_{\tau(3)\tau(4)}(x,y) \Bigg] \times  \mE \Bigg[ \sum_{x' \in [\ell_1], y' \in [\ell_2]} \psi^{m,\boldsymbol{\pi}}_{\tau'(1) \tau'(2)}(x',y') \psi^{m,\boldsymbol{\pi}}_{\tau'(3)\tau'(4)}(x',y') \Bigg].
\end{align*}
By focusing on the first term, Cauchy--Schwarz inequality together with the previous bound~(\ref{Eq: TV bound}) yields
\begin{align*}
	& \Bigg| \mE \Bigg[\sum_{x \in [\ell_1], y \in [\ell_2]} \psi^{m,\boldsymbol{\pi}}_{\tau(1) \tau(2)}(x,y) \psi^{m,\boldsymbol{\pi}}_{\tau(3)\tau(4)}(x,y) \Bigg] \Bigg| \\[.5em]
	~ \leq ~ & \sqrt{\sum_{x \in [\ell_1], y \in [\ell_2]} \Big\{ \mE\big[ \psi^{m,\boldsymbol{\pi}}_{\tau(1)\tau(2)}(x,y)\big] \Big\}^2} \times \sqrt{\sum_{x \in [\ell_1], y \in [\ell_2]} \Big\{ \mE\big[ \psi^{m,\boldsymbol{\pi}}_{\tau(3)\tau(4)}(x,y)\big] \Big\}^2} \\[.5em]
	\leq ~ & \frac{C_2}{M^2}.
\end{align*}
The other term can be similarly analyzed and thus 
\begin{align*}
\mE [h_{i_1,i_2,i_3,i_4}^{m,\boldsymbol{\pi}}] \times \mE[h_{j_1,j_2,j_3,j_4}^{m,\boldsymbol{\pi}}]  ~=~  O\left( M^{-4} \right).
\end{align*}
Combining the above with the bound~(\ref{Eq: expectation of the product}) yields the claim~(\ref{Eq: goal (1 common index)}).

Suppose that there are at least one collision among 
\begin{align*}
	\{ & (X_{i_1,m},Y_{\pi_{i_1,m}}),(X_{i_2,m},Y_{\pi_{i_2,m}}),(X_{i_3,m},Y_{\pi_{i_3,m}}),(X_{i_4,m},Y_{\pi_{i_4,m}}), \\[.5em]
	&(X_{j_1,m},Y_{\pi_{j_1,m}}),(X_{j_2,m},Y_{\pi_{j_2,m}}),(X_{j_3,m},Y_{\pi_{j_3,m}}),(X_{j_4,m},Y_{\pi_{j_4,m}})\}.
\end{align*}
where $\{i_1,i_2,i_3,i_4,j_1,j_2,j_3,j_4\} \in A_{\mathrm{com}}^1$. The number of such cases is at most $O(\sigma_m^6)$. Therefore 
\begin{align*}
	\frac{1}{\binom{\sigma_m}{4}^2} \sum_{A_{\mathrm{com}}^{1}} \mathrm{Cov}(h_{i_1,i_2,i_3,i_4}^{m,\boldsymbol{\pi}}, h_{j_1,j_2,j_3,j_4}^{m,\boldsymbol{\pi}}) ~ = ~ O\left( \frac{1}{\sigma_m M^2} + \frac{1}{\sigma_m^2} \right).
\end{align*}

Lastly we focus on the case when $\{i_1,i_2,i_3,i_4,j_1,j_2,j_3,j_4\} \in A_{\mathrm{com}}^0$. First of all, when there is no collision, two kernels become independent and thus $\mathrm{Cov}(h_{i_1,i_2,i_3,i_4}^{m,\boldsymbol{\pi}}, h_{j_1,j_2,j_3,j_4}^{m,\boldsymbol{\pi}}) = 0$. When there is one collision, we are basically in the same situation where $\{i_1,i_2,i_3,i_4,j_1,j_2,j_3,j_4\} \in A_{\mathrm{com}}^1$ and there is no collision. Therefore we can re-use the previous bound for the covariance~(\ref{Eq: goal (1 common index)}). Furthermore, since there are at most $O(\sigma_m^7)$ cases of one collision, $O(\sigma_m^6)$ of two collisions and so on, we have
\begin{align*}
	\frac{1}{\binom{\sigma_m}{4}^2} \sum_{A_{\mathrm{com}}^{0}} \mathrm{Cov}(h_{i_1,i_2,i_3,i_4}^{m,\boldsymbol{\pi}}, h_{j_1,j_2,j_3,j_4}^{m,\boldsymbol{\pi}}) ~=~ O\left( \frac{1}{\sigma_m M^{2}} + \frac{1}{\sigma_m^2} \right).
\end{align*}
Therefore, combining the results yields
\begin{align*}
	\mathrm{Var}[U(\bW_m^{\boldsymbol{\pi}})|\sigma_m] ~=~ O\left( \frac{1}{\sigma_m M^{2}} + \frac{1}{\sigma_m^2} \right).
\end{align*}
Given the above bound, the variance of the combined test statistic is 
\begin{align*}
	\mathrm{Var}[T_{\mathrm{CI}}^{\boldsymbol{\pi}}] ~=~  \mE\big[\mathrm{Var}[T_{\mathrm{CI}}^{\boldsymbol{\pi}}|\boldsymbol{\sigma}]\big] + \mathrm{Var}\big[\mE[T_{\mathrm{CI}}^{\boldsymbol{\pi}}|\boldsymbol{\sigma}]\big].
\end{align*}
Let us upper bound these two terms separately. Starting with the first term, note that
\begin{align*}
	\mathrm{Var}[T_{\mathrm{CI}}^{\boldsymbol{\pi}}|\boldsymbol{\sigma}] ~= ~ &  \sum_{m \in [M]} \mathds{1}(\sigma_m \geq 4) \sigma_m^2 \times \mathrm{Var}[U(\bW_m^{\boldsymbol{\pi}})|\sigma_m] \\[.5em]
	= ~ & \sum_{m \in [M]} \mathds{1}(\sigma_m \geq 4) \sigma_m^2 \times O\left( \frac{1}{\sigma_m M^2} + \frac{1}{\sigma_m^2} \right),
\end{align*}
which results in
\begin{align} \label{Eq: the first bound}
	\mE\big[\mathrm{Var}[T_{\mathrm{CI}}^{\boldsymbol{\pi}}|\boldsymbol{\sigma}]\big] = O\left( \frac{n}{M^2} + M \right).
\end{align}
For the second term,
\begin{align*}
	\mE[T_{\mathrm{CI}}^{\boldsymbol{\pi}}|\boldsymbol{\sigma}] ~=~ \sum_{m \in [M]} \mathds{1}(\sigma_m \geq 4) \sigma_m \times \mE[U(\bW_m^{\boldsymbol{\pi}})| \sigma_m]. 
\end{align*}
Since $\sigma_1,\ldots,\sigma_M$ are independent,
\begin{align*}
	\mathrm{Var}\big[\mE[T_{\mathrm{CI}}^{\boldsymbol{\pi}}|\boldsymbol{\sigma}]\big] ~=~ & \sum_{m \in [M]} \mathrm{Var}\big[\mathds{1}(\sigma_m \geq 4) \sigma_m \times \mE[U(\bW_m^{\boldsymbol{\pi}})| \sigma_m]\big] \\[.5em]
	\leq ~ &  \sum_{m \in [M]} \mE\big[\mathds{1}(\sigma_m \geq 4) \sigma_m^2 \times \big\{\mE[U(\bW_m^{\boldsymbol{\pi}})| \sigma_m]\big\}^2 \big].
\end{align*}
Recall from the claim~(\ref{Eq: bounding the expectation}) that 
\begin{align} 
	\mE[U(\bW_m^{\boldsymbol{\pi}})|\sigma_m] = O\left( \frac{1}{M^2} + \frac{1}{M\sigma_m} + \frac{1}{\sigma_m^2} \right) \quad \text{for each $m \in [M]$.}
\end{align}
Therefore
\begin{align*}
	\mathrm{Var}\big[\mE[T_{\mathrm{CI}}^{\boldsymbol{\pi}}|\boldsymbol{\sigma}]\big] ~\leq~ & \frac{C_1}{M^4} \sum_{m \in [M]} \mE\big[ \mathds{1}(\sigma_m \geq 4) \sigma_m^2 \big] + \frac{C_2}{M^2} \sum_{m \in [M]} \mE\big[ \mathds{1}(\sigma_m \geq 4) \big] \\[.5em]
	+ ~& C_3 \sum_{m \in [M]}  \mE\big[ \mathds{1}(\sigma_m \geq 4) \sigma_m^{-2} \big].
\end{align*}
Now using the property that $\sigma_m \overset{\mathrm{i.i.d.}}{\sim} \text{Pois}(np_m)$ where $\sum_{m=1}^M p_m = 1$ and $p_m \geq 0$, we can upper bound 
\begin{align*}
	& \sum_{m \in [M]} \mE\big[ \mathds{1}(\sigma_m \geq 4) \sigma_m^2 \big] \leq \sum_{m \in [M]} np_m + n^2p_m^2 \leq n + n^2 \quad \text{and} \\[.5em]
	& \sum_{m \in [M]} \mE\big[ \mathds{1}(\sigma_m \geq 4) \big] \leq  \sum_{m \in [M]} \mE[\sigma_m] \leq n.
\end{align*}
Furthermore, the inequality~(\ref{Eq: poisson reciprocal}) implies that 
\begin{align*}
	\sum_{m \in [M]}  \mE\big[ \mathds{1}(\sigma_m \geq 4) \sigma_m^{-2} \big] ~\leq~ O\left( \frac{M^2}{n} \right).
\end{align*}
Putting pieces together yields
\begin{align*}
	\mathrm{Var}\big[\mE[T_{\mathrm{CI}}^{\boldsymbol{\pi}}|\boldsymbol{\sigma}]\big] ~\leq~ O\left( \frac{n^2}{M^4} + \frac{n}{M^2} + \frac{M^2}{n} \right).
\end{align*}
Combining the above with (\ref{Eq: the first bound}) concludes
\begin{align*}
	\mathrm{Var}[T_{\mathrm{CI}}^{\boldsymbol{\pi}}] ~=~ O\left( \frac{n^2}{M^4} + \frac{n}{M^2} + \frac{M^2}{n} + M\right).
\end{align*}

\vskip 1em

\noindent \textbf{Step 4. Completing the proof.} Now by taking $M = \ceil{n^{2/5}}$, the previous bounds for the mean and the variance yield
\begin{align*}
	\mE[T_{\mathrm{CI}}^{\boldsymbol{\pi}}] + \sqrt{\mathrm{Var}[T_{\mathrm{CI}}^{\boldsymbol{\pi}}]} ~=~ O\big(n^{1/5}\big).
\end{align*}
This together with the inequality~(\ref{Eq: summary}) gives
\begin{align*}
	\mP(T_{\mathrm{CI}} \leq q_{1-\alpha}) ~\leq~ \mP\big(T_{\mathrm{CI}} \lesssim n^{1/5} \big) + \frac{1}{200} + \mP(\mathcal{A}^c).
\end{align*}
Finally, we prove that $\mP(\mathcal{A}^c)$ converges to zero uniformly over a class of distributions of $(X,Y,Z) \in [\ell_1] \times [\ell_2] \times [0,1]$ where the density of $Z$ is lower bounded by some small constant $c_{\mathrm{low}}>0$. In particular, we recall the event $\mathcal{A}$ in (\ref{Eq: event A}) and define
\begin{align} \label{Eq: delta}
	\rho_n :=  \sup_{P_{X,Y,Z} \in \mathcal{P}'_{1,[0,1],\mathrm{TV}}(L)} \mP_{P_{X,Y,Z}^N,N}(\mathcal{A}^c).
\end{align}
We then show that $\rho_n \rightarrow 0$ under the conditions in Theorem~\ref{Theorem: Type II error of the double-binned test}. 

By letting $\sigma_{m,k}$ be the sample size within the bin $B_{m,k}$ and $\tau_\alpha  > \log^{-1}((1-\alpha/2)/(1-\alpha))$, let us define another event $\mathcal{E}$ in which $\min_{(m,k) \in [M] \times [b]} \sigma_{m,k} \geq \tau_{\alpha} \times (M b)$. Note that, under Poissonization, each $\sigma_{m,k} \sim \text{Pois}(n p_{m,k})$ where $np_{m,k} \geq c_{\mathrm{low}} n / (Mb)$ since the density of $Z$ is bounded below by $c_{\mathrm{low}}$. We further assume that $c_{\mathrm{low}}n/(Mb) \geq 2\tau_\alpha \times (Mb)$. Then, by Lemma~\ref{Lemma: Poisson lower bound} in Appendix~\ref{Section: Poissonization}, the event $\mathcal{E}$ holds with probability at least 
\begin{align*}
	1 - \sum_{(m,k) \in [M] \times [b]} e^{-\frac{np_{m,k}}{6}} \geq 1 - Mb \times e^{- \frac{c_{\mathrm{low}} n}{6Mb}}.
\end{align*}
This probability converges to one as $n\rightarrow \infty$ as long as $M,b$ are chosen such that $n/(Mb) \rightarrow \infty$. Therefore, we can work under the assumption that $\min_{(m,k) \in [M] \times [b]} \sigma_{m,k} \geq \tau_\alpha \times (M b)$. Under this event $\mathcal{E}$ with $\tau_\alpha(M b) \geq 2$ (this can be assumed as $Mb \rightarrow \infty$), $\mP(\mathcal{A}^c|\mathcal{E})$ is equivalent to 
\begin{align*}
	& \mP\Bigg( \frac{\prod_{(m,k) \in [M] \times [b]} \max\{\sigma_{m,k}-1,1\}}{\prod_{(m,k) \in [M] \times [b]} \max\{\sigma_{m,k},1\}} \leq \frac{1-\alpha}{1-\alpha/2} \bigg| \mathcal{E} \Bigg) \\[.5em]
	\leq ~ & \mP\Bigg( \prod_{(m,k) \in [M] \times [b]} (1 - \sigma_{m,k}^{-1}) \leq \frac{1-\alpha}{1-\alpha/2}  \bigg| \mathcal{E} \Bigg)  \\[.5em]
	\leq ~ & \mathds{1}\Bigg(  \big(1 - \tau_{\alpha}^{-1}(Mb)^{-1} \big)^{Mb}  \leq \frac{1-\alpha}{1-\alpha/2} \Bigg).
\end{align*}
Given that $\tau_\alpha  > \log((1-\alpha/2)/(1-\alpha))$ and $(1-a/n)^n \rightarrow e^{-a}$ for any constant $a$, the last display converges to zero as $Mb \rightarrow \infty$. Finally, as we argued before, we have $\mP(\mathcal{E}^c) \rightarrow 0$. Therefore
\begin{align*}
	\mP(\mathcal{A}^c)  ~=~& 	\mP(\mathcal{A}^c|\mathcal{E}) \mP(\mathcal{E})  + 	\mP(\mathcal{A}^c|\mathcal{E}^c) \mP(\mathcal{E}^c) \\[.5em]
	\leq ~ &  	\mP(\mathcal{A}^c|\mathcal{E})  +  \mP(\mathcal{E}^c)  \rightarrow 0,
\end{align*}
as desired. Since we set $\alpha = 0.01$, we can take $\tau_\alpha = 200 > \log^{-1}((1-\alpha/2)/(1-\alpha)) \approx 198.5$ and the above conditions are satisfied when $n/(Mb)^2 \geq 400 c_{\mathrm{low}}^{-1}$ and $Mb \rightarrow \infty$ as $n\rightarrow \infty$. These are the conditions imposed in Theorem~\ref{Theorem: Type II error of the double-binned test} and thus we complete the proof of Theorem~\ref{Theorem: Type II error of the double-binned test}.

\vskip 2em

\section{Auxiliary lemmas and additional results} \label{Section: Auxiliary lemmas}
This section collects some auxiliary lemmas used for the main proofs as well as additional results.

\subsection{Metrics and inequalities}
The next lemma connects the Hellinger distance between binned distributions $Q_{X,Y|\widetilde{Z}=m}$ and $\widetilde{Q}_{X,Y|\widetilde{Z}=m}$, defined in Section~\ref{Section: Validity}, to that between the underlying conditional distributions. This lemma is the basis for the proof of Theorem~\ref{Theorem: Validity of the permutation test under Hellinger smooth} and Theorem~\ref{Theorem: Validity of the permutation test under Renyi smooth}.
\begin{lemma}[A bound on the Hellinger distance] \label{Lemma: Bound on the Hellinger distance}
	For each $m=1,\ldots,M$, we let $Z_m$ denote a random variable from the conditional distribution $Z|Z \in B_m$, and we let $Z_m^{'},Z_m^{''}$ be i.i.d.~copies of $Z_m$. Then under $X \independent Y|Z$, it holds that 
	\begin{align*}
		& \mathcal{D}_{\text{\emph{H}}}^2 \big(Q_{X,Y|\widetilde{Z}=m},\widetilde{Q}_{X,Y|\widetilde{Z}=m}\big) \\[.5em] 
		\leq ~ & 6 \mE_{Z_m,Z_m^{'},Z_m^{''}}\big[ \mathcal{D}_{\text{\emph{H}}}^2\big(P_{X|Z=Z_m},P_{X|Z=Z_m^{'}}\big) \times \mathcal{D}_{\text{\emph{H}}}^2\big(P_{Y|Z=Z_m},P_{Y|Z=Z_m^{''}}\big)\big].
	\end{align*}
\end{lemma}
It is interesting to observe that Lemma~\ref{Lemma: Bound on the Hellinger distance} connects the distance between the joint conditional distributions of the binned data to the distances between the marginal conditional distributions of the unbinned data. In particular, the upper bound consists of the product between two distances, which allows us to control the type I error more tightly. The constant factor of $6$ is a convenient choice from our analysis and it is by no means tight. To prove this result, we use Jensen's inequality for a vector valued function and the proof can be found in Appendix~\ref{Section: Proof of Lemma: Bound on the Hellinger distance}.

The next lemma presents an inclusion property of the generalized Hellinger distance, which is useful for the proof of Theorem~\ref{Theorem: Validity of the permutation test under Hellinger smooth}. 
\begin{lemma}[Generalized Hellinger distance] \label{Lemma: Generalized Hellinger distance}
	For any $1 \leq \gamma_1 \leq \gamma_2$, we have 
	\begin{align} \label{Eq: Hellinger inequality 1}
		\mathcal{D}_{\gamma_2,\text{\emph{H}}}^{\gamma_2}(P,Q) \leq \mathcal{D}_{\gamma_1,\text{\emph{H}}}^{\gamma_1}(P,Q). 
	\end{align}
	On the other hand, for any $\gamma \geq 1$ and $1 \leq \alpha \leq 2$, it holds that
	\begin{align} \label{Eq: Hellinger inequality 2}
		\mathcal{D}_{\gamma,\text{\emph{H}}}(P,Q) \leq 2^{1 - 1/\gamma + 1/(\gamma \alpha)} \mathcal{D}_{\gamma \alpha, \text{\emph{H}}}(P,Q). 
	\end{align}
\end{lemma}
We note that the first inequality~(\ref{Eq: Hellinger inequality 1}) was proved in \cite{kamps1989hellinger} and, as far as we are aware, the second inequality~(\ref{Eq: Hellinger inequality 2}) is new in the literature. In particular, the inequality~(\ref{Eq: Hellinger inequality 2}) generalizes the well-known inequality between the TV distance and the Hellinger distance:
\begin{align*}
	\mathcal{D}_{\mathrm{TV}}(P,Q) \leq \sqrt{2} \mathcal{D}_{\mathrm{H}}(P,Q), 
\end{align*}
by taking $\gamma = 1$ and $\alpha = 2$. The proof of Lemma~\ref{Lemma: Generalized Hellinger distance} can be found in Appendix~\ref{Section: Proof of Lemma: Generalized Hellinger distance}. As a direct consequence of Lemma~\ref{Lemma: Generalized Hellinger distance}, we have the following corollary. 
\begin{corollary}[A bound on the Hellinger distance] \label{Corollary: Hellinger bound}
	Suppose that $\gamma > 2$ is within the range of $2^{k} < \gamma \leq 2^{k+1}$ for some positive integer $k$. Then 
	\begin{align} \label{Eq: upper bound on Hellinger}
		\mathcal{D}_{\text{\emph{H}}}(P,Q) \leq 2^{k} \mathcal{D}_{\gamma, \text{\emph{H}}}(P,Q). 
	\end{align}
\end{corollary}
Next we state a few inequalities related to R{\'e}nyi divergence, which will be useful for the proof of Theorem~\ref{Theorem: Validity of the permutation test under Renyi smooth}. Let $\mathcal{D}_{\text{KL}}(P\|Q) =\int p \log (p/q) d\mu$ and $\mathcal{D}_{\chi^2}(P\|Q) =\int (p-q)^2/q d\mu$ be the KL divergence and the $\chi^2$ divergence of $P$ from $Q$, respectively. 
Then we have the identities:
\begin{align*}
	& \mathcal{D}_{1/2 \text{,R{\'e}nyi}}(P \| Q) =  -2 \log \{ 1 - \mathcal{D}_{\mathrm{H}}^2(P,Q) /2\}, \\[.5em]
	& \mathcal{D}_{1 \text{,R{\'e}nyi}}(P \| Q) = \mathcal{D}_{\text{KL}}(P\|Q), \\[.5em]
	& \mathcal{D}_{2 \text{,R{\'e}nyi}}(P \| Q) =\log \{1 + \mathcal{D}_{\chi^2}(P\|Q)\}.
\end{align*}
Notice that $ \mathcal{D}_{\gamma, \text{R{\'e}nyi}}(P \| Q) $ is nondecreasing in $\gamma$ and thus, by the inequality $\log t \leq t - 1$ for all $t > 0$, we have
\begin{align} \label{Eq: Renyi inequality}
	\mathcal{D}_{\mathrm{H}}^2(P,Q) \leq \mathcal{D}_{1/2 \text{,R{\'e}nyi}}(P \| Q) \leq  \mathcal{D}_{\text{KL}}(P\|Q) = R_{1} (f \| g)  \leq  \mathcal{D}_{2 \text{,R{\'e}nyi}}(P \| Q) \leq \mathcal{D}_{\chi^2}(P\|Q).
\end{align} 
Other notable inequalities are
\begin{align*}
	& \frac{\gamma_1}{1 - \gamma_1} \frac{1 - \gamma_2}{\gamma_2} \mathcal{D}_{\gamma_2 \text{,R{\'e}nyi}}(P \| Q) \leq \mathcal{D}_{\gamma_1 \text{,R{\'e}nyi}}(P \| Q) \leq \mathcal{D}_{\gamma_2 \text{,R{\'e}nyi}}(P \| Q)\quad \text{for $0 < \gamma_1 < \gamma_2 < 1$}.
\end{align*}
Based on the above inequalities, we can upper bound the Hellinger distance by R{\'e}nyi divergence as follows.
\begin{lemma}[R{\'e}nyi divergence] \label{Lemma: Renyi divergence}
	For any $\gamma \geq 1/2$, we have
	\begin{align*}
		\mathcal{D}_{\text{\emph{H}}}^2(P,Q) \leq \mathcal{D}_{\gamma, \text{\emph{R{\'e}nyi}}}(P \| Q).
	\end{align*}
	On the other hand, for $0 < \gamma < 1/2$, it holds that
	\begin{align*}
		\mathcal{D}_{\text{\emph{H}}}^2(P,Q) \leq \mathcal{D}_{1/2, \text{\emph{R{\'e}nyi}}}(P \| Q) \leq \frac{1-\gamma}{\gamma} \mathcal{D}_{\gamma, \text{\emph{R{\'e}nyi}}}(P \| Q).
	\end{align*}
\end{lemma}
The proof of the above lemma can be found in \cite{van2014renyi}.

\subsection{Connection between permutation $p$-value and quantile}

The next lemma provides an alternative expression of the permutation test in terms of the quantile of the permutation distribution, which is more intuitive to analyze. 

\begin{lemma}[Quantile] \label{Lemma: quantile}
	Recall the permutation $p$-value $p_{\text{\emph{perm}}}$ in (\ref{Eq: permutation p-value}). Let $q_{1-\alpha}$ be the $1-\alpha$ quantile of the empirical distribution of $T_{\text{\emph{CI}}}^{\boldsymbol{\pi}_1},\ldots,T_{\text{\emph{CI}}}^{\boldsymbol{\pi}_K}$. Then it holds that 
	\begin{align*}
		\mathds{1}(p_{\text{\emph{perm}}} > \alpha)  = \mathds{1}(T_{\text{\emph{CI}}} \leq q_{1-\alpha}). 
	\end{align*}
	Furthermore, by letting $k = \ceil{(1-\alpha)K}$, $q_{1-\alpha}$ is the same as the $k$th order statistic $T_{\mathrm{CI}}^{\boldsymbol{\pi}_{(k)}}$ of $T_{\text{\emph{CI}}}^{\boldsymbol{\pi}_1},\ldots,T_{\text{\emph{CI}}}^{\boldsymbol{\pi}_K}$.
\end{lemma}
The proof of this result follows by the definition of the quantile function and the details can be found in Appendix~\ref{Section: Proof of Lemma: quantile}.

\subsection{Poissonization} \label{Section: Poissonization}

In the following proposition, we demonstrate that the local permutation test under Poissonization has the same validity as before. 

\begin{proposition}[Validity under Poissonization] \label{Proposition: Poissonization}
	Let $\mathcal{P}_{0}' \subset \mathcal{P}_0$ be the collection of distributions for which it holds that $\sup_{P_{X,Y,Z} \in \mathcal{P}_{0}'} \mE_{P_{X,Y,Z}^n}[\phi_{\text{\emph{perm}},n}] \leq \alpha + C_1 n^{-\epsilon}$ where $C_1 > 0$ is a universal constant, $\epsilon > 0$ and $n \geq 1$. We further let $\phi_{\text{\emph{perm}},n} = 0$ when $n=0$. Then under Poissonization where  $N \sim \text{\emph{Pois}}(n)$, the local permutation test has the type I error guarantee as
	\begin{align*}
		\sup_{P_{X,Y,Z} \in \mathcal{P}_{0}'} \mE_{P_{X,Y,Z}^N,N}[\phi_{\text{\emph{perm}},N}] \leq \alpha +  C_2(\epsilon) \cdot n^{-\epsilon} + ne^{-n/12},
	\end{align*} 
	where $C_2(\epsilon) > 0$ is a constant that only depends on $\epsilon$. 
\end{proposition}

The proof of Proposition~\ref{Proposition: Poissonization} uses the fact that a Poisson random variable is tightly concentrated around its mean, and the details can be found in Appendix~\ref{Section: Proof of Proposition: Poissonization}. Suppose we choose the maximum diameter $h$ such that the upper bounds in Theorem~\ref{Theorem: Validity of the permutation test under Hellinger smooth} and Theorem~\ref{Theorem: Validity of the permutation test under Renyi smooth} converge to $\alpha$ at a certain rate. Then Proposition~\ref{Proposition: Poissonization} essentially implies that the same convergence rate of the type I error can be obtained under Poissonization as well.

The following lemma provides a concentration bound for the minimum of Poisson random variables.
\begin{lemma}[Poisson concentration] \label{Lemma: Poisson lower bound}
	For $m=1,\ldots,M$, let $X_m$ be a Poisson random variable with parameter $\lambda_m \geq c_0$. Then 
	\begin{align*}
		\mP\Big(\min_{m \in [M]} X_m \leq c_0/2\Big) \leq  \sum_{m=1}^M e^{-\frac{\lambda_m}{6}}.
	\end{align*}
\end{lemma}

\subsection{Unbounded support of $Z$} \label{Section: Unbounded support}

In this subsection, we discuss how to deal with unbounded support of $Z$ for type I error control. In particular, we demonstrate that it is possible to obtain similar validity results as in Theorem~\ref{Theorem: Validity of the permutation test under Hellinger smooth} and Theorem~\ref{Theorem: Validity of the permutation test under Renyi smooth}, up to a vanishing term in the sample size, even when the support of the conditional variable $Z$ is unbounded.

To fix ideas, suppose $Z$ is a univariate random variable and let $\{B_1,\ldots,B_M\}$ denote a partition of $[a,b]$ for some $a,b \in \mathbb{R}$. We further denote the complement of $\cup_{i=1}^{M} B_i$ by $B_{M+1} = (-\infty,a) \cup (b,\infty)$ so that $\cup_{i=1}^{M+1} B_i = \mathbb{R}$. Under this setting, the permutation test that builds on partitions $\{B_1,\ldots,B_{M+1}\}$ has the same type I error guarantee as those in Theorem~\ref{Theorem: Validity of the permutation test under Hellinger smooth} and Theorem~\ref{Theorem: Validity of the permutation test under Renyi smooth} except that we have an additional term to control: $n q_{\widetilde{Z}}(m+1)  \times \mathcal{D}_{\mathrm{H}}^2 (Q_{X,Y|\widetilde{Z}=m+1},\widetilde{Q}_{X,Y|\widetilde{Z}=m+1})$. In particular, since the Hellinger distance is bounded by one, this extra term goes to zero as long as $n q_{\widetilde{Z}}(m+1) \rightarrow 0$ as $n \rightarrow \infty$. Therefore, the problem boils down to controlling $n q_{\widetilde{Z}}(m+1)$. We discuss two strategies to handle this depending on whether there exists prior knowledge on the distribution of $Z$. 
\begin{enumerate}
	\item \textbf{Prior information on $Z$.} Suppose that we have prior knowledge on the distribution of $Z$. To fix ideas, suppose that $Z$ has an exponential tail such that $\mP(|Z| > t) \leq Ce^{-t}$ for some $C>0$ and for all $t \geq 0$. Then by choosing $b$ sufficiently larger than $\log n$ and by letting $a=-b$, it is guaranteed that the extra term vanishes as $n q_{\widetilde{Z}}(m+1) \rightarrow 0$. 
	
	\item \textbf{Sample splitting.} When there is no information on the distribution of $Z$ a priori, we may choose $a$ and $b$ in a data-dependent way. In detail, we split the sample into two sets of size $n_1$ and $n_2$, denoted by $D_1$ and $D_2$, respectively. Let $Z_{(1)}$ and $Z_{(n_1)}$ be the minimum and the maximum order statistic of the first set $D_1$ and let $a = Z_{(1)}$ and $b= Z_{(n_1)}$. We then work with the second set $D_2$ to perform the local permutation test. Provided that $Z$ has a density with respect to the Lebesgue measure, it can be seen that $\mP(Z < Z_{(1)}) = \frac{1}{n_1+1}$ and $\mP(Z > Z_{(n_1)}) = \frac{1}{n_1+1}$ where $Z$ is independent of $D_1$. Note that $q_{\widetilde{Z}}(m+1) = \mP(Z < Z_{(1)}|D_1) + \mP(Z>Z_{(n_1)}|D_1)$ and therefore $\mE_{D_1}[q_{\widetilde{Z}}(m+1) ] = \frac{2}{n_1+1}$. In view of (\ref{Eq: Error bound in terms of the Hellinger distance}), the type I error is bounded by
	\begin{align*}
		& \mE_{P_{X,Y,Z}^n}[\phi_{\mathrm{perm},n_2}] = \mE_{D_1} [\mE_{D_2}[\phi_{\mathrm{perm},n_2|D_1}]] \\[.5em] 
		 ~\leq~ & \alpha + \mE_{D_1}\Bigg[ \Bigg\{2n_2 \sum_{m=1}^{M+1} q_{\widetilde{Z}}(m)  \times \mathcal{D}_{\mathrm{H}}^2\big(Q_{X,Y|\widetilde{Z}=m},\widetilde{Q}_{X,Y|\widetilde{Z}=m}\big) \Bigg\}^{1/2}\Bigg] \\[.5em]
		\leq ~ & \alpha +  \underbrace{\mE_{D_1}\Bigg[ \Bigg\{2n_2 \sum_{m=1}^{M} q_{\widetilde{Z}}(m)  \times \mathcal{D}_{\mathrm{H}}^2\big(Q_{X,Y|\widetilde{Z}=m},\widetilde{Q}_{X,Y|\widetilde{Z}=m}\big) \Bigg\}^{1/2}\Bigg]}_{(\mathrm{I})} \\[.5em]
		& ~~ + \underbrace{\mE_{D_1}\Bigg[ \Bigg\{2n_2 q_{\widetilde{Z}}(M+1)  \times \mathcal{D}_{\mathrm{H}}^2\big(Q_{X,Y|\widetilde{Z}=M+1},\widetilde{Q}_{X,Y|\widetilde{Z}=M+1}\big) \Bigg\}^{1/2}\Bigg]}_{(\mathrm{II})}.
	\end{align*}
	The first term $(\mathrm{I})$ can be handled similarly as in the bounded case. For the second term~$(\mathrm{II})$, it can be seen that $(\mathrm{II}) \leq \sqrt{\frac{4n_2}{n_1+1}}$ by using the fact that the Hellinger distance is bounded by one together with Jensen's inequality. Thus, the second term is asymptotically vanishing as long as $n_2/n_1 \rightarrow 0$ as $n \rightarrow \infty$. 
\end{enumerate}

\section{Definitions} \label{Section: Definitions}

We recall a H\"{o}lder smoothness class from \cite{neykov2020minimax}, which is used to analyze the local permutation test for continuous data. 

\begin{definition}[H\"{o}lder Smoothness] \label{Definition: Holder smoothness} Let $s > 0$ be a fixed real number, and let $\lfloor s \rfloor$ denote the maximum integer strictly smaller than  $s$. Denote by $\mathcal{H}^{2,s}(L)$, the class of functions $f:[0,1]^2 \mapsto \mathbb{R}$, which posses all partial derivatives up to order $\lfloor s \rfloor$ and for all $x,y,x',y' \in [0,1]$ we have
	\begin{align*} 
		\sup_{k \leq \lfloor s \rfloor}\bigg|\frac{\partial^k}{\partial x^k}\frac{\partial^{\lfloor s \rfloor - k}}{\partial y^{\lfloor s \rfloor - k}}f(x,y) - \frac{\partial^k}{\partial x^k}\frac{\partial^{\lfloor s \rfloor - k}}{\partial y^{\lfloor s \rfloor - k}}f(x',y') \bigg| \leq L((x-x')^2 + (y - y')^2))^{\frac{s - \lfloor s \rfloor}{2}},
	\end{align*}
	and in addition
	\begin{align*}
		\sup_{k \leq \lfloor s \rfloor}\bigg|\frac{\partial^k}{\partial x^k}\frac{\partial^{\lfloor s \rfloor - k}}{\partial y^{\lfloor s \rfloor - k}}f(x,y)\bigg| \leq L.
	\end{align*}
\end{definition}

The next definition describes a cyclic permutation used in Section~\ref{Section: Double-binning strategy}. 

\begin{definition}[Cyclic permutation] \label{Definition: cyclic permutation}
	For a set $\mathcal{X}$ with finite elements $a_0,a_2,\ldots,a_{n-1}$, a cyclic permutation $\pi$ is a bijective function from $\mathcal{X}$ to $\mathcal{X}$ such that for some integer $k$, $\pi(a_i) = a_{(i+k) \bmod  n}$.
\end{definition}

\vskip 2em

\section{Proofs of the auxiliary lemmas}

This section collects the proof of the lemmas in Appendix~\ref{Section: Auxiliary lemmas}.

\subsection{Proof of Lemma~\ref{Lemma: Bound on the Hellinger distance}} \label{Section: Proof of Lemma: Bound on the Hellinger distance}
To ease the notation, we suppress $(x,y|m)$ in $q_{XY|\widetilde{Z}}(x,y|m)$ and write $q_{XY|\widetilde{Z}}$. Similarly, we write $q_{X\cdot|\widetilde{Z}}(x|m)$ and $q_{\cdot Y|\widetilde{Z}}(y|m)$ by $q_{X\cdot|\widetilde{Z}}$ and $q_{\cdot Y|\widetilde{Z}}$, respectively. Using this shorthand notation, consider an alternative expression of the Hellinger distance and its upper bounds
\begin{align*}
	2 \mathcal{D}_{\mathrm{H}}^2 \big(Q_{X,Y|\widetilde{Z}=m},\widetilde{Q}_{X,Y|\widetilde{Z}=m}\big) ~=~ &   \int \int  \big(\sqrt{q_{XY|\widetilde{Z}}} - \sqrt{q_{X\cdot|\widetilde{Z}} q_{\cdot Y|\widetilde{Z}}} \big)^2 d\mu_Xd\mu_Y \\[.5em]
	~=~ & \int \int \frac{ \big( {q_{XY|\widetilde{Z}}} - {q_{X\cdot|\widetilde{Z}} q_{\cdot Y|\widetilde{Z}}} \big)^2 }{ \big( \sqrt{q_{XY|\widetilde{Z}}} + \sqrt{q_{X\cdot|\widetilde{Z}} q_{\cdot Y|\widetilde{Z}}} \big)^2 } d\mu_Xd\mu_Y \\[.5em]
	\leq~ &  \int \int  \frac{\big({q_{XY|\widetilde{Z}}} - {q_{X\cdot|\widetilde{Z}} q_{\cdot Y|\widetilde{Z}}}\big)^2}{{q_{XY|\widetilde{Z}}} + {q_{X\cdot|\widetilde{Z}} q_{\cdot Y|\widetilde{Z}}}} d\mu_Xd\mu_Y \\[.5em]
	\leq ~ & 3 \int \int  \frac{\big({q_{XY|\widetilde{Z}}} - {q_{X\cdot|\widetilde{Z}} q_{\cdot Y|\widetilde{Z}}}\big)^2}{{q_{XY|\widetilde{Z}}} + 3 {q_{X\cdot|\widetilde{Z}} q_{\cdot Y|\widetilde{Z}}}} d\mu_Xd\mu_Y.
\end{align*}
Under $X \independent Y|Z$, the term inside of the last expression is equivalent to
\begin{align*}
	\frac{ \big( q_{XY|\widetilde{Z}} - q_{X\cdot|\widetilde{Z}} q_{\cdot Y|\widetilde{Z}} \big)^2}{ q_{XY|\widetilde{Z}} + 3q_{X\cdot|\widetilde{Z}} q_{\cdot Y|\widetilde{Z}}} ~=~ & \frac{\big\{ \mE_{Z_m}[ (p_{X|Z_m} - q_{X\cdot|\widetilde{Z}}) (p_{Y|Z_m} - q_{\cdot Y|\widetilde{Z}})] \}^2}{\mE_{Z_m}[ (p_{X|Z_m} + q_{X\cdot|\widetilde{Z}}) (p_{Y|Z_m} + q_{\cdot Y|\widetilde{Z}})]},
\end{align*}
where we denote $p_{X|Z_m} = p_{X|Z}(x|Z_m)$ and $p_{Y|Z_m} = p_{Y|Z}(y|Z_m)$. Next we apply Jensen's inequality several times and obtain that
\begin{align*}
	\frac{\big\{ \mE_{Z_m}[ (p_{X|Z_m} - q_{X\cdot|\widetilde{Z}}) (p_{Y|Z_m} - q_{\cdot Y|\widetilde{Z}})] \}^2}{\mE_{Z_m}[(p_{X|Z_m} + q_{X\cdot|\widetilde{Z}}) (p_{Y|Z_m} + q_{\cdot Y|\widetilde{Z}})]} ~ \overset{(\text{i})}{\leq} ~ & \mE_{Z_m} \Bigg[ \frac{(p_{X|Z_m} - q_{X\cdot|\widetilde{Z}})^2 (p_{Y|Z_m} - q_{\cdot Y|\widetilde{Z}})^2 }{(p_{X|Z_m} + q_{X\cdot|\widetilde{Z}}) (p_{Y|Z_m} + q_{\cdot Y|\widetilde{Z}})} \Bigg] \\[.5em]
	~ \overset{(\text{ii})}{\leq} ~ & \mE_{Z_m,Z_m',Z_m''} \Bigg[ \frac{(p_{X|Z_m} - p_{X|Z_m'})^2 (p_{Y|Z_m} - p_{Y|Z_m''})]^2 }{(p_{X|Z_m} + p_{X|Z_m'}) (p_{Y|Z_m} + p_{Y|Z_m''})} \Bigg]
\end{align*}
where we detail step~(i) and step~(ii) at the end of this section. Applying Tonelli's theorem yields
\begin{align*}
	& 3 \int \int  \frac{\big({q_{XY|\widetilde{Z}}} - {q_{X\cdot|\widetilde{Z}} q_{\cdot Y|\widetilde{Z}}}\big)^2}{{q_{XY|\widetilde{Z}}} + 3 {q_{X\cdot|\widetilde{Z}} q_{\cdot Y|\widetilde{Z}}}} d\mu_Xd\mu_Y \\[.5em]
	~\leq~ & 3 \mE_{Z_m,Z_m',Z_m''} \Bigg[ \int \frac{(p_{X|Z_m} - p_{X|Z_m'})^2 }{p_{X|Z_m} + p_{X|Z_m'}} d\mu_X \times  \int \frac{(p_{Y|Z_m} - p_{Y|Z_m''})^2 }{p_{Y|Z_m} + p_{Y|Z_m''}} d\mu_Y \Bigg]  \\[.5em]
	\leq ~ & 12 \mE_{Z_m,Z_m^{'},Z_m^{''}}\big[ \mathcal{D}_{\mathrm{H}}^2\big(P_{X|Z=Z_m},P_{X|Z=Z_m^{'}}\big) \times \mathcal{D}_{\mathrm{H}}^2\big(P_{Y|Z=Z_m},P_{Y|Z=Z_m^{''}}\big)\big],
\end{align*}
where the last inequality holds since $(x+y) \leq (\sqrt{x} + \sqrt{y})^2 \leq 2(x + y)$ for $x,y \geq 0$. Hence the result follows. 

\vskip 1em 

\noindent \textbf{Details of step~(i) and step~(ii)}

\noindent Throughout this part, we assume that $q_{X\cdot|\widetilde{Z}} >0$ and $q_{\cdot Y|\widetilde{Z}} >0$. Otherwise we have $q_{XY|\widetilde{Z}} - q_{X\cdot|\widetilde{Z}}q_{\cdot Y|\widetilde{Z}} =0$ and the bound becomes trivial. Based on the identity
\begin{align*}
	& \frac{\big\{ \mE_{Z_m}[ (p_{X|Z} - q_{X\cdot|\widetilde{Z}}) (p_{Y|Z} - q_{\cdot Y|\widetilde{Z}})] \}^2}{\mE_{Z_m}[ (p_{X|Z} + q_{X\cdot|\widetilde{Z}}) (p_{Y|Z} + q_{\cdot Y|\widetilde{Z}})]} \\[.5em] 
	= ~ & \frac{\big\{ \mE_{Z_m}[ p_{X|Z_m}p_{Y|Z_m}] - q_{\cdot Y|\widetilde{Z}}  \mE_{Z_m}[ p_{X|Z_m}] - q_{X\cdot|\widetilde{Z}} \mE_{Z_m}[ p_{Y|Z_m}] + q_{X\cdot|\widetilde{Z}} q_{\cdot Y|\widetilde{Z}} \}^2}{\mE_{Z_m}[ p_{X|Z_m}p_{Y|Z_m}] + q_{\cdot Y|\widetilde{Z}} \mE_{Z_m}[p_{X|Z_m} ] + q_{X\cdot|\widetilde{Z}} \mE_{Z_m}[p_{Y|Z_m}] + q_{X\cdot|\widetilde{Z}} q_{\cdot Y|\widetilde{Z}}},
\end{align*}
it is enough to show that the function
\begin{align*}
	f(x,y,z) = \frac{(x - ay -bz + ab)^2}{x + ay + bz + ab}
\end{align*}
is convex for $x,y,z \geq 0$ and $a,b > 0$. Then the result follows by Jensen's inequality for a vector valued function. To this end, we prove that the Hessian matrix of $f(x,y,z)$ is positive semidefinite. Some calculations show that
\begin{align*}
	\frac{\partial^2 f}{\partial x^2} ~=~& \frac{8 (ay + bz)^2}{(x + bz + ay + ab)^3}, \\[.5em]
	\frac{\partial^2 f}{\partial y^2} ~=~& \frac{8 a^2(x+ab)^2}{(x + bz + ay + ab)^3}, \\[.5em]
	\frac{\partial^2 f}{\partial z^2} ~=~& \frac{8 b^2(x+ab)^2}{(x + bz + ay + ab)^3}, \\[.5em]
	\frac{\partial^2 f}{\partial x \partial y} ~=~& -\frac{8 a(x+ab)(ay+bz)}{(x + bz + ay + ab)^3}, \\[.5em]
	\frac{\partial^2 f}{\partial x \partial z} ~=~& -\frac{8 b(x+ab)(ay+bz)}{(x + bz + ay + ab)^3}, \\[.5em]
	\frac{\partial^2 f}{\partial y \partial z} ~=~& -\frac{8 ab(x+ab)^2}{(x + bz + ay + ab)^3},
\end{align*}
and thus the Hessian matrix becomes
\begin{align*}
	\textbf{H}_f ~=~ \frac{8}{(x + bz+ ay + ab)^3} 
	\begin{bmatrix}
		ay + bz \\
		-a(x+ab) \\
		-b(x+ab)
	\end{bmatrix} 
	\begin{bmatrix}
		ay + bz & -a(x+ab)  & -b(x+ab)
	\end{bmatrix}.
\end{align*}
This completes the proof of step~(i). Step~(ii) follows similarly by conditional Jensen's inequality together with the observation that $\frac{(a - x)^2}{b + x}$ is convex for $x \geq 0$ and $a,b \geq 0$.

\subsection{Proof of Lemma~\ref{Lemma: Generalized Hellinger distance}} \label{Section: Proof of Lemma: Generalized Hellinger distance}
The first inequality~(\ref{Eq: Hellinger inequality 1}) is stated in Corollary 2 of \cite{kamps1989hellinger}. Hence we focus on the second inequality~(\ref{Eq: Hellinger inequality 2}). We first claim that, for any $1 \leq \alpha \leq \beta$ such that $1/\alpha + 1/\beta = 1$,
\begin{align} \label{Eq: x-y inequality}
	|x - y| \leq |x^{1/\alpha} - y^{1/\alpha}| \times |x^{1/\beta} + y^{1/\beta}|.
\end{align}
Without loss of generality, assume that $x > y \geq 0$ (when $x=y$, there is nothing to prove). Then, by expanding the terms, the claim becomes equivalent to 
\begin{align*}
	x^{1/\alpha} y^{1/\beta} \geq x^{1/\beta} y^{1/\alpha}.
\end{align*}
By taking the log on both sides,
\begin{align*}
	\left(\frac{1}{\alpha} - \frac{1}{\beta} \right) \log x \geq \left( \frac{1}{\alpha} - \frac{1}{\beta} \right) \log y,
\end{align*}
which is true since $1 \leq \alpha \leq \beta$ and $x > y \geq 0$. This verifies the inequality~(\ref{Eq: x-y inequality}). Leveraging this inequality, observe that
\begin{align*}
	\int \big| p^{1/\gamma} - q^{1/\gamma} \big|^\gamma d\mu \leq \int \big| p^{1/(\gamma \alpha)} - q^{1/(\gamma \alpha)} \big|^\gamma \big| p^{1/(\gamma \beta)} + q^{1/(\gamma\beta)} \big|^\gamma d\mu
\end{align*}
for any $1 \leq \alpha \leq \beta$ and $1/\alpha + 1/\beta = 1$. By H\"{o}lder's inequality, we further have
\begin{align*}
	\bigg[ \int \big| p^{1/\gamma} - q^{1/\gamma} \big|^\gamma d\mu \bigg]^{1/\gamma} \leq \bigg[ \int \big| p^{1/(\gamma \alpha)} - q^{1/(\gamma \alpha)} \big|^{\gamma \alpha} d\mu \bigg]^{1/(\gamma \alpha)}  \bigg[ \int \big| p^{1/(\gamma \beta)} + q^{1/(\gamma \beta)} \big|^{\gamma \beta} d\mu \bigg]^{1/(\gamma \beta)}. 
\end{align*}
Since $p,q$ are density functions, the second term in the product is bounded above by 
\begin{align*}
	\bigg[ \int \big| p^{1/(\gamma \beta)} + q^{1/(\gamma \beta)} \big|^{\gamma \beta} d\mu \bigg]^{1/(\gamma \beta)} = ~  \bigg[  2^{\gamma\beta} \int \bigg| \frac{p^{1/(\gamma \beta)} + q^{1/(\gamma \beta)}}{2} \bigg|^{\gamma \beta} d\mu \bigg]^{1/(\gamma \beta)} 	\leq ~   2,
\end{align*}
where the last step follows by Jensen's inequality (note that $\gamma \beta \geq 1$). Therefore
\begin{align*}
	\mathcal{D}_{\gamma,\mathrm{H}} (P,Q) =  \bigg[ \frac{1}{2}\int \big| p^{1/\gamma} - q^{1/\gamma} \big|^\gamma d\mu \bigg]^{1/\gamma} \leq ~ & 2^{1 - 1/\gamma}\bigg[ \int \big| p^{1/(\gamma \alpha)} - q^{1/(\gamma \alpha)} \big|^{\gamma \alpha} d\mu \bigg]^{1/(\gamma\alpha)} \\[.5em] 
	=~ & 2^{1 - 1/\gamma + 1/(\gamma\alpha)}\mathcal{D}_{\gamma \alpha, \mathrm{H}} (P,Q).
\end{align*}
Note that by the restriction on $\alpha$ and $\beta$, the parameter $\alpha$ should be in the range of $1 \leq \alpha \leq 2$. This completes the proof of Lemma~\ref{Lemma: Generalized Hellinger distance}.

\subsection{Proof of Lemma~\ref{Lemma: quantile}} \label{Section: Proof of Lemma: quantile}
By the definition of the quantile function, $q_{1-\alpha}$ can be written as
\begin{align*}
	q_{1-\alpha} ~=~& \inf \bigg\{x \in \mathbb{R}: 1-\alpha \leq  \frac{1}{K} \sum_{\boldsymbol{\pi}_i \in \boldsymbol{\Pi}} \mathds{1} \big( T_{\mathrm{CI}}^{\boldsymbol{\pi}_i} \leq x \big) \bigg\} \\[.5em]
	= ~ &  \inf \bigg\{x \in \mathbb{R}: 1-\alpha \leq  \frac{1}{K} \sum_{\boldsymbol{\pi}_i \in \boldsymbol{\Pi}} \mathds{1} \big( T_{\mathrm{CI}}^{\boldsymbol{\pi}_i} < x \big) \bigg\} \\[.5em]
	= ~ & T_{\mathrm{CI}}^{\boldsymbol{\pi}_{(k)}},
\end{align*}
where we recall $k = \ceil{(1-\alpha)K}$. Given this representation and by noting that the permutation distribution is discrete, we see that $T_{\mathrm{CI}} \leq q_{1-\alpha}$ holds if and only if
\begin{align*}
	\frac{1}{K} \sum_{\boldsymbol{\pi}_i \in \boldsymbol{\Pi}} \mathds{1} \big( T_{\mathrm{CI}}^{\boldsymbol{\pi}_i} < T_{\mathrm{CI}} \big) < 1-\alpha,
\end{align*}
which is equivalent to $p_{\mathrm{perm}} > \alpha$. This completes the proof. 

\subsection{Proof of Proposition~\ref{Proposition: Poissonization}} \label{Section: Proof of Proposition: Poissonization}
Since we assume that $\phi_{\mathrm{perm},k} = 0$ when $k=0$, it holds that
\begin{align*}
	\sup_{P_{X,Y,Z} \in \mathcal{P}_{0}'} \bigg\{ \sum_{k=0}^\infty \mP(N = k) \mE_{P_{X,Y,Z}^k} [\phi_{\text{{perm}},k}] \bigg\} ~ = ~ & 	\sup_{P_{X,Y,Z} \in \mathcal{P}_{0}'} \bigg\{ \sum_{k=1}^\infty \mP(N = k) \mE_{P_{X,Y,Z}^k} [\phi_{\text{{perm}},k}] \bigg\} \\[.5em]
	\leq ~ & \alpha + C_1 \sum_{k=1}^\infty \mP(N = k) k^{-\epsilon}. 
\end{align*} 
Furthermore, by transforming $k$ to $k=k'+1$, we can simplify the infinite sum as
\begin{align*}
	\sum_{k=1}^\infty \mP(N = k) k^{-\epsilon} ~=~ & \sum_{k=1}^\infty \frac{n^k e^{-n}}{k!} k^{-\epsilon} =  n \sum_{k'=0}^\infty \mP(N=k') (k'+1)^{-\epsilon - 1}  \\[.5em] 
	~ = ~ & n \mE[(N+1)^{-\epsilon - 1}],
\end{align*}
where we recall that $N \sim \text{Pois}(n)$. Now decompose and upper bound the expectation by
\begin{align*}
	n \mE[(N+1)^{-\epsilon - 1}] ~=~  &n \mE[(N+1)^{-\epsilon - 1} \mathds{1}\{ N \leq n/2\} ] + n \mE[(N+1)^{-\epsilon - 1} \mathds{1}\{ N > n/2\} ] \\[.5em]
	\leq ~ & n \mP(N \leq n/2) + n \left[ \frac{1}{n/2+1}\right]^{1+\epsilon} \\[.5em]
	\leq ~ & n e^{-n/12} + 2^{1+\epsilon} n^{-\epsilon},
\end{align*}
where the last inequality uses a Chernoff bound for a Poisson random variable  \citep[e.g.][]{canonne2019note}. This completes the proof of Proposition~\ref{Proposition: Poissonization}.

\subsection{Proof of Lemma~\ref{Lemma: Poisson lower bound}} 
Let $X$ have a Poisson distribution with parameter $\lambda>0$. Then by using a Poisson tail bound~\citep[e.g.][]{canonne2019note}, we have for $x>0$,
\begin{align*}
	\mP(X \leq \lambda - x)  \leq e^{-\frac{x^2}{\lambda + x}}.
\end{align*}
Now, by setting $x = \lambda/2$, the above inequality guarantees that 
\begin{align} \label{Eq: individual poisson}
	\mP(X \leq \lambda/2) \leq e^{-\frac{\lambda}{6}}.
\end{align}
By the union bound along with (\ref{Eq: individual poisson}), we have
\begin{align*}
	\mP\Big(\min_{m \in [M]}X_m \leq c/2\Big) \leq \sum_{m=1}^M \mP (X_m \leq \lambda_m/2) \leq \sum_{m=1}^M e^{-\frac{\lambda_m}{6}},
\end{align*}
which proves the result.

\end{document}